\input amstex
\documentstyle{amsppt}
\TagsOnRight
\CenteredTagsOnSplits
\NoRunningHeads
\NoBlackBoxes

\define\la{\lambda}
\define\La{\Lambda}
\define\R{{\Bbb R}}
\define\C{{\Bbb C}}
\define\Y{{\Bbb Y}}
\define\Z{{\Bbb Z}}

\define\Om{\Omega}
\define\om{\omega}
\define\al{\alpha}
\define\be{\beta}
\define\ga{\gamma}
\define\Ga{\Gamma}
\define\de{\delta}
\define\De{\Delta}

\define\si{\sigma}

\define\Reg{\operatorname{Reg}}
\define\End{\operatorname{End}}
\define\lev{\operatorname{lev}}
\define\const{\operatorname{const}}
\define\supp{\operatorname{supp}}
\define\tr{\operatorname{tr}}
\define\Prob{\operatorname{Prob}}
\define\Ex{\operatorname{Ex}}
\define\sgn{\operatorname{sgn}}
\define\diag{\operatorname{diag}}

\define\wt{\widetilde}
\define\wh{\widehat}

\define\tht{\thetag}

\topmatter
\title
    Harmonic analysis on the infinite symmetric group
\endtitle
\author
$\boxed{\text{Sergei Kerov}}$\,, Grigori Olshanski, and Anatoly
Vershik
\endauthor

\thanks  Grigori Olshanski: Institute for Information Transmission Problems,
Bolshoy Karetny 19, Moscow 127994, GSP-4, Russia. \quad E-mail:
{\tt olsh\@online.ru} \newline\indent Anatoly Vershik: Steklov
Institute of Mathematics at St.~Petersburg, Fontanka 27,
St.~Petersburg 191011, Russia. \quad E-mail: {\tt
vershik\@pdmi.ras.ru} \endthanks

\abstract The infinite symmetric group $S(\infty)$, whose
elements are finite permutations of $\{1,2,3,\dots\}$, is a
model example of a ``big'' group. By virtue of an old result of
Murray--von Neumann, the one--sided regular representation of
$S(\infty)$ in the Hilbert space $\ell^2(S(\infty))$ generates a
type II$_1$ von Neumann factor  while the two--sided regular
representation is irreducible. This shows that the conventional
scheme of harmonic analysis is not applicable to $S(\infty)$:
for the former representation, decomposition into irreducibles
is highly non--unique, and for the latter representation, there
is no need of any decomposition at all. We start with
constructing a compactification $\frak S\supset S(\infty)$,
which we call the space of virtual permutations. Although $\frak
S$ is no longer a group, it still admits a natural two--sided
action of $S(\infty)$. Thus, $\frak S$ is a $G$--space, where
$G$ stands for the product of two copies of $S(\infty)$. On
$\frak S$, there exists a unique $G$-invariant probability
measure $\mu_1$, which has to be viewed as a ``true'' Haar
measure for $S(\infty)$. More generally, we include $\mu_1$ into
a family $\{\mu_t:\, t>0\}$ of distinguished $G$-quasiinvariant
probability measures on virtual permutations. By making use of
these measures, we construct a family $\{T_z:\, z\in\C\}$ of
unitary representations of $G$, called generalized regular
representations (each representation $T_z$ with $z\ne0$ can be
realized in the Hilbert space $L^2(\frak S, \mu_t)$, where
$t=|z|^2$). As $|z|\to\infty$, the generalized regular
representations $T_z$ approach, in a suitable sense, the
``naive'' two--sided regular representation of the group $G$ in
the space $\ell^2(S(\infty))$. In contrast with the latter
representation, the generalized regular representations $T_z$
are highly reducible and have a rich structure. We prove that
any $T_z$ admits a (unique) decomposition into a multiplicity
free continuous integral of irreducible representations of $G$.
For any two distinct (and not conjugate) complex numbers $z_1$,
$z_2$, the spectral types of the representations $T_{z_1}$ and
$T_{z_2}$ are shown to be disjoint. In the case $z\in\Z$, a
complete description of the spectral type is obtained. Further
work on the case $z\in\C\setminus\Z$ reveals a remarkable link
with stochastic point processes and random matrix theory.
\endabstract

\endtopmatter

\document

\newpage

\head Contents \endhead

0. Introduction

1. The space of virtual permutations

2. Measures on virtual permutations

3. Generalized regular representations

4. Distinguished spherical function

5. The commutant and the block decomposition

6. $K$-invariant vectors

7. The commutant of representation $T_{pq}$

8. Disjointness of spectral measures

9. Appendix

References

\head 0. Introduction \endhead

\subhead 0.1. The infinite symmetric group: characters, factor
representations, spherical representations \endsubhead The
present paper deals with harmonic analysis on the infinite
symmetric group $S(\infty)$ (the group of finite permutations of
the infinite set $\{1,2,\dots\}$). This group belongs to the
class of {\it big groups\/}, \footnote{This term, introduced in
Vershik \cite{Ver}, is somewhat vague but expressive and
convenient. About representations of various big groups, see,
e.g., Ismagilov \cite{Ism}, Neretin \cite{Ner1}, Olshanski
\cite{Ol4}, Str{\u a}til{\u a}--Voiculescu \cite{SV}.} such as
the infinite dimensional classical groups, diffeomorphism and
current groups. We present here the first example of harmonic
analysis on a big group. A short exposition of our results was
published in \cite{KOV}; here we provide the detailed statements
and proofs.

The most well--studied groups in representation theory --- compact,
abelian, and reductive Lie groups --- are {\it tame}. Irreducible
unitary representations of tame groups can be classified, and the
basic problem of harmonic analysis consists in decomposing
interesting reducible representations (e.g., the regular
representation) into irreducible components.

The group $S(\infty)$ is {\it wild}, not tame. The wild groups
do not admit any sensible classification of irreducible
representations, and reducible representations can be decomposed
into irreducibles in essentially different ways. Unlike tame
groups, the wild groups have factor representations of von
Neumann types II and III. All those facts imply that the
representation theory of wild groups should be based on the
principles distinct from those used for tame groups.
\footnote{For general facts about tame and wild groups, see,
e.g., Dixmier \cite{Dix}, Kirillov \cite{Kir}. Note that tame
groups are also called {\it type I} groups. About factor
representations, see, e.g., Dixmier \cite{Dix}, Thoma
\cite{Tho2}.}

The possibility of developing such a theory  was first
demonstrated  by E.~Thoma \cite{Tho1} exactly for the example of
$S(\infty)$. In that paper Thoma discovered that the finite type
factor representations of the group $S(\infty)$ admit an
explicit classification. More precisely, those representations
are labelled by the pairs $\omega=(\alpha;\beta)$ of
nonincreasing sequences of nonnegative numbers
$$
\alpha = (\alpha_1\ge\alpha_2\ge\ldots\ge0), \qquad
\beta = (\beta_1\ge\beta_2\ge\ldots\ge0)
$$
subject to the extra condition
$$
\sum_{k=1}^\infty \alpha_k + \sum_{k=1}^\infty \beta_k \le 1.
$$
The set of such pairs forms an infinite--dimensional simplex
$\Omega$ which we call the {\it Thoma simplex\/}. \footnote{For
more detail on Thoma's theorem, see \S9.6 below.}

For an arbitrary unitary representation of a countable group,
generating a von Neumann algebra of finite type, a decomposition
into a direct integral of factor representations always exists
and is essentially unique. Therefore, the main problem of
harmonic analysis on $S(\infty)$ may consist in the actual
decomposition of most notable representations of finite type.
This agrees with an old idea (especially developed by Pukanszky)
that as ``elementary representations'' for a wild group, one
should consider not irreducible representations but factor
representations of an appropriate class.

Instead of decompositions of representations one can speak of
decompositions of their {\it characters}. Denote by $\Cal X$ the
set of positive definite functions $\chi$ on the group
$S(\infty)$, constant on conjugacy classes and normalized at the
identity permutation by the condition $\chi(e)=1$. In the
topology of pointwise convergence, the space $\Cal X$ is a
convex compact set. It is known that the extreme points of the
set $\Cal X$ are exactly the characters of factor
representations of finite type. Moreover, the space $\Cal X$ is
a Choquet simplex, that is, any point $\chi\in \Cal X$ can be
uniquely represented as a continual convex combination of
extreme characters.

However, as was shown by one of the authors (\cite{Ol2},
\cite{Ol3}), one can develop another approach to representation
theory of the group $S(\infty)$, which makes it possible to
avoid factor representations. Recall that, by definition, the
group $S(\infty)$ consists of {\it finite\/} permutations of the
set $\{1,2,\ldots\}$, fixing all but finitely many points.
Consider now the group $\overline{S(\infty)}$ of {\it all\/}
bijections of the set $\{1,2,\ldots\}$ onto itself. Taking the
subsets
$$
\overline{S_m(\infty)}= \{g\in\overline{S(\infty)}:g(k)=k,\; k=1,\,\ldots,m\},
$$
where $m=1,2,\dots$, as a fundamental system of neighborhoods of
identity we turn $\overline{S(\infty)}$ into a topological
group.

The subgroup $S(\infty)$ is dense in $\overline{S(\infty)}$.
Furthermore, denote by $\overline G$ the group of pairs
$(g_1,g_2)\in\overline{S(\infty)}\times\overline{S(\infty)}$,
such that $g_1^{-1}g_2\in S(\infty)$, and identify
$\overline{S(\infty)}$ with the diagonal subgroup $\overline K$
in $\overline G$. We introduce a topology in $\overline G$ by
proclaiming the subgroup $\overline K$ open. The group
$S(\infty)\times S(\infty)$ is dense in $\overline G$. The
groups $\overline K$ and $\overline G$ are totally disconnected
and {\it not\/} locally compact.

Remarkably enough, the topological group $\overline G$ is tame
and its irreducible representations can be completely described,
see Olshanski \cite{Ol3}, Okounkov \cite{Ok1}, \cite{Ok2}. In
the present paper we are only interested in {\it spherical\/}
irreducible representations of the group $\overline G$ (that is,
representations containing a $\overline K$--invariant vector).
It is known that $(\overline G,\overline K)$ is a {\it Gelfand
pair\/}: in any irreducible spherical representation there is
only one, up to a scalar factor, $\overline K$--invariant vector
$\xi$.

There exists a natural bijection $T\leftrightarrow\pi$ between
irreducible spherical representations $T$ of the pair
$(\overline G,\overline K)$ and factor representations $\pi$ of
finite type of the group $S(\infty)$. Specifically, $\pi$ is the
restriction of $T$ to the subgroup
$S(\infty)\times\{e\}\subset\overline G$. The character
$\chi(s)=\operatorname{Tr}\,\pi(s)$ of the representation $\pi$
(here ``$\operatorname{Tr}$'' denotes the normalized trace on
the factor) is related to the spherical function
$\varphi(g)=(T(g)\xi,\xi)$ by the simple formula
$$
\varphi(g_1,g_2) = \chi(g_1^{-1}g_2), \qquad
(g_1,g_2) \in \overline G.
$$
Therefore, we can think of {\it finite\/} factor representations of the group
$S(\infty)$ as of irreducible spherical representations of the pair $(\overline
G,\overline K)$, thus returning to the conventional setup of representation
theory of tame groups.

\subhead 0.2. The generalized regular representations $T_z$ and the problem of
harmonic analysis  \endsubhead We shall now describe the representations of the
group $\overline G$ whose decomposition is the purpose of the present paper.
The choice of those representations is itself a nontrivial problem.
Traditionally, the object of harmonic analysis for a Gelfand pair $(\Cal G,\Cal
K)$ is the decomposition of the natural representation of $\Cal G$ in the space
$L^2(\Cal G/\Cal K)$ (for instance, for Riemannian symmetric spaces $\Cal
G/\Cal K$, the decomposition problem was studied in classical works of
Harish--Chandra and Gindikin--Karpelevich). In our situation the space
$\overline G/\overline K\cong S(\infty)$ is discrete, so that the Hilbert space
$L^2(\overline G/\overline K)=\ell^2(S(\infty))$ makes sense. However, the
corresponding representation turns out to be irreducible so that there is no
need of any further decomposition. This irreducibility effect is equivalent to
the following fact (which is probably better known): the one--sided regular
representation of the group $S(\infty)$  in the space $\ell^2(S(\infty))$
generates a  type II$_1$ factor. \footnote{More generally, for any discrete
group whose conjugacy classes, except $\{e\}$, are infinite, the one--sided
regular representation is a type II$_1$ factor representation, while the
two--sided regular representation is irreducible. See Murray--von Neumann
\cite{MvN, chapter 5},  Naimark \cite{Nai, chapter VII, \S38.5}.}

Instead of the homogeneous space $\overline G/\overline
K=S(\infty)$ we introduce a compact space $\frak S$ containing
$S(\infty)$ as a dense subset. More precisely, the space $\frak
S$ is defined as a projective limit of finite sets,
$$
S(1) \leftarrow S(2) \leftarrow \ldots \leftarrow S(n)
\leftarrow \ldots.
$$
Here $S(n)$ is the set of all permutations of $n$ objects, and
the projections are specified in \S1 below. The points of $\frak
S$ are called {\it virtual permutations\/}. Unlike $S(\infty)$,
the space $\frak S$ is {\it not} a group. However, it admits a
canonical action of the group $\overline G$, which is sufficient
for our purposes. In particular, we have a two--sided action of
the group $S(\infty)$ on $\frak S$.

On the space $\frak S\supset S(\infty)$ of virtual permutations
there is a (unique) probability measure invariant with respect
to the action of the group $\overline G$. This measure, which we
denote as $\mu_1$, is much more interesting and useful than the
counting (Haar) measure on $S(\infty)$. More generally, the
measure $\mu_1$ can be included into a one--parameter family of
probability measures $\mu_t$, $t>0$. All those measures are
invariant with respect to the subgroup $\overline K$, and
quasiinvariant with respect to $\overline G$.\; \footnote{The
idea of extending the group space (of a big group) in order to
obtain measures with good transformation properties comes from
the measure theory on infinite--dimensional linear spaces and is
rather old, see, e.g., Gelfand--Vilenkin \cite{GV, chapter IV}.}

By applying a standard general construction of producing unitary
representations from a group action with quasiinvariant measure,
we arrive at a family of representations $T_z$ of the group
$\overline G$. Here $z$ ranges over the set $\C\setminus\{0\}$,
and the representation $T_z$ acts in the space $L^2(\frak
S,\mu_t)$ where $t=|z|^2$. One can also define two more unitary
representations $T_0$, $T_\infty$ in such a way that the
resulting family is continuously parametrized by the points of
the Riemann sphere $\C\cup\{\infty\}$. The representation
$T_\infty$ coincides with the above mentioned two--sided regular
representation in the space $\ell^2(S(\infty))$. We regard the
family $\{T_z\}$ as a deformation of $T_\infty$, and we call
$T_z$ the {\it generalized regular representations}.
\footnote{Notice a certain similarity between our family
$\{T_z\}$ and Neretin's deformation of the natural
representation on $L^2(\Cal G/\Cal K)$, where $\Cal G/\Cal K$ is
a Riemannian symmetric space, see Neretin \cite{Ner3}.}

The only irreducible representation in the family $\{T_z:\,z\in\Bbb
C\cup\{\infty\}\}$ is $T_\infty$. All the representations $T_z$, $z\in\C$, are
reducible, and their structure is rather complicated. {\it We consider the
explicit decomposition of those representations as the main problem of harmonic
analysis for the group $S(\infty)$.}

Each representation $T_z$ can also be realized as an inductive
limit of the two--sided regular representations of finite groups
$S(n)\times S(n)$,
$$
\Reg^1\to \Reg^2\to \ldots \to \Reg^n\to \ldots, \tag *
$$
with very special embeddings $\Reg^n\to\Reg^{n+1}$ depending on
$z$. The existence of a finite dimensional approximation allows
one to employ the powerful combinatorial and probabilistic
techniques used in the theory of approximately
finite--dimensional ($AF$--) algebras. In case of the group
$S(\infty)$, these techniques, based on combinatorics of Young
diagrams and the theory of symmetric functions, were developed
in \cite{VK2}, \cite{VK3}.

\subhead 0.3. Main results of the paper \endsubhead These are as
follows:

1. The representations $T_z$ and $T_{\bar z}$ are equivalent, and we describe
an intertwining operator realizing their equivalence.

2. For any $z\in \C$, the representation $T_z$ can be decomposed into a
multiplicity free direct integral of irreducible spherical representations. So,
the equivalence class of $T_z$  is completely described by an equivalence class
of measures on the Thoma simplex $\Omega$. We will refer to the latter
equivalence class as to the {\it spectral type\/} of $T_z$.

3. The spectral types of two representations $T_{z_1}$ and
$T_{z_2}$ are disjoint (mutually singular) whenever $z_1$ and
$z_2$ are not equal or conjugate to each other. This means that
there exist two disjoint Borel subsets in $\Omega$ supporting
the measures from the spectral types of $T_{z_1}$ and $T_{z_2}$,
respectively. Equivalently, there is no intertwining operator
between $T_{z_1}$ and $T_{z_2}$.

4. The spectral type of $T_z$ substantially depends on whether $z$ is an
integer or not. In the present paper we focus on the case when $z\in\Z$. Then
we are able to describe the spectral type quite explicitly. Namely, let
$\Omega(p,q)$ denote the subset of those pairs $(\alpha,\beta)\in\Omega$ for
which $\alpha_1+\dots+\alpha_p+\beta_1+\dots+\beta_q=1$ (consequently, all the
remaining coordinates $\alpha_i$ and $\beta_j$ vanish); notice that
$\Omega(p,q)$ is a ($p+q-1$)--dimensional face of the Thoma simplex $\Omega$.
Our result says that when $z=0,\pm1,\pm2,\dots$, the spectral type of $T_z$ is
determined by the union of Lebesgue measures on the finite--dimensional faces
$\Omega(p,q)$ with $p-q=z$ (which agrees with the result on disjointness of
spectral types stated above). We also show that if $z\in\C\setminus\Z$, then
the spectral type of $T_z$ is concentrated ``inside'' the simplex $\Omega$.
That is, all faces $\Omega(p,q)$ are null sets with respect to the spectral
type of $T_z$.

\subhead 0.4. The case $z\in\C\setminus\Z$ \endsubhead This case
is studied in detail in a series of papers by Borodin and
Olshanski, see \cite{Ol5}, \cite{Bor2}, \cite{BO3}, the
expository papers \cite{BO1}, \cite{Ol6}, and references
therein. When $z\in\Bbb C\setminus\Z$, there exists a
distinguished probability measure in the spectral type of $T_z$
(in the present paper it is denoted as $\si_z$). The key idea is
to assign to $(\Om,\si_z)$ a {\it stochastic point process\/},
and to study its {\it correlation functions.\/}  It turns out
that the point processes thus obtained are close to those
arising in random matrix theory. The link with random matrix
theory seems to be especially interesting and promising.

The results of Borodin and Olshanski were obtained as a continuation of the
project started in \cite{KOV}. Together with the present work they provide a
description of the spectral types for all representations $T_z$.

\subhead 0.5. Organization of the paper \endsubhead

In \S1, we start with the definition of the {\it canonical projection}
$p_n:S(n)\to S(n-1)$. Using the projections $p_n$  we define the  space $\frak
S=\varprojlim S(n)$ of virtual permutations. Then we describe four different
concrete realizations of that space. In one of them, $\frak S$ turns into the
product of an infinite sequence of finite sets.  Finally, we show that $\frak
S$ is a $G$--space, where $G=S(\infty)\times S(\infty)$, and we introduce an
additive $\Z$--valued 1--cocycle for the action of $G$ on $\frak S$.

In \S2, we study the family $\{\mu_t\}_{t>0}$ of probability
measures on the space $\frak S$. This is a deformation of the
unique $G$--invariant probability measure (in our notation,
$\mu_1$). The measures $\mu_t$ are central, i.e., invariant with
respect to the action of the diagonal subgroup $K\subset G$. We
show that they turn into product measures in one of the
realizations of $\frak S$. Moreover, they are essentially the
only central measures with this property. Applying Kakutani's
theorem we show that the measures $\mu_t$ are pairwise disjoint.
Then we show that any $\mu_t$ is quasiinvariant with respect to
the action of the group $G$, i.e., for any $g\in G$, the shift
$\mu_t^g$ of $\mu_t$ by $g$ is a measure in the equivalence
class of $\mu_t$. We also calculate the Radon--Nikodym
derivative $\mu_t^g/\mu_t$, which we will need later.

In \S3, we construct the representations $T_z$ in two different
ways. First, we realize $T_z$, where $z\in\C^*$, in the Hilbert
space $L^2(\frak S,\mu_t)$, where $t=|z|^2$. Here we use a
multiplicative 1--cocycle $\frak S\times G\to \C^*$ which
depends on $z$ and is defined via the additive cocycle from \S2.
Second, we realize the same representation $T_z$ as the
inductive limit corresponding to a chain of embeddings \tht{*},
see above. Using the latter realization, it is easy to complete
the family $\{T_z\}$ by two limit representations, $T_0$ and
$T_\infty$. Finally, we present a transparent interpretation of
the embeddings in \tht{*}.

In \S4, we define and study a distinguished matrix coefficient
of $T_z$. The representation $T_z$ has a distinguished
$K$--invariant vector $\xi_0$: in the first realization, this is
the identity function $f_0\equiv 1$ on the space $\frak S$.
Though the space of $K$--invariant vectors is infinite
dimensional, $\xi_0$ is the only $K$--invariant vector  which is
seen at once: constructing other examples of $K$--invariant
vectors is already a nontrivial task (we do this in section 6).
To the vector $\xi_0$ one associates a spherical function on
$G$, $\varphi_z(g_1,g_2)=(T(g_1,g_2)\xi_0,\xi_0)$, and a
character $\chi_z(s)=\varphi(s,e)$ of the group $S(\infty)$. We
do not dispose of a simple expression for the values of the
function $\chi_z$ on conjugacy classes of the group $S(\infty)$;
instead of this, we find a very nice formula for the
coefficients $M_z(\lambda)$ in the expansions
$$
\chi_z\mid_{S(n)}=\sum_{\lambda\in\Bbb Y_n}M_z(\lambda)\widetilde\chi^\lambda,
\qquad n=1,2,\dots, \tag**
$$
where $\Bbb Y_n$ denotes the set of Young diagrams with $n$ boxes, and
$\widetilde\chi^\lambda$ is the normalized irreducible character of $S(n)$,
indexed by $\lambda$. In other words, we get the ``Fourier coefficients'' of
$\chi_z$ (see Theorem 4.2.1).

This explicit formula for $M_z(\lambda)$ plays a key role in the
present paper. \footnote{For alternative derivations of the
formula and generalizations, see Kerov \cite{Ker3}, Borodin's
appendix in Olshanski \cite{Ol5}, Borodin--Olshanski
\cite{BO2}.} We derive from it the following results:

First, we prove that $\xi_0$ is a cyclic vector if (and only if) $z\in\Bbb
C\setminus\Z$. This means that for nonintegral values of $z$, the spectral type
of $T_z$ is entirely determined by the decomposition of the matrix coefficient
associated with $\xi_0$.

Second, we prove the equivalence $T_z\sim T_{\bar z}$, which is not
evident from the
construction of the representations. We also exhibit an isometric
operator intertwining $T_z$ and $T_{\bar z}$.

Finally, following the general philosophy of \cite{VK2},
\cite{VK3}, we assign to the family $M_z=\{M_z^{(n)}\}$ the
so--called transition probabilities $p_z(\lambda,\nu)$. Given
$\lambda\in\Bbb Y_n$, the numbers $p_z(\lambda,\nu)$ form a
probability distribution on the set of those diagrams $\nu$ that
can be obtained from $\lambda$ by adding a single box.

In \S5, we deal with the realization of $T_z$ as an inductive limit of regular
representations of the finite groups $S(n)\times S(n)$. To any such inductive
limit we assign a ``transition function'' defined on the edges of the Young
graph. This leads to a description of the commutant of the representation. We
explicitly calculate the transition function for the representations $T_z$.
Using this, we get, for integral values of the parameter $z$, a decomposition
of $T_z$ into a direct sum of subrepresentations $T_{pq}$ which we call the
{\it blocks\/} of $T_z$.

In \S6, we get a convenient realization of the subspace $V_z$ of $K$--invariant
vectors in the Hilbert space of $T_z$. In terms of this realization we
construct, for any $z\in\Z$, a $K$--invariant vector in each block of $T_z$,
and we calculate the spectral decomposition of the corresponding matrix
coefficient.

In \S7, we prove that all the vectors constructed in section 6 are cyclic
vectors in the corresponding blocks. Together with the results of section 6
this provides us with a complete description of the spectral type of the
representation $T_z$ in case of integral $z$.  The argument of this section
turns out to be rather long. At the end we give an example illustrating the
origin of one of the difficulties.

In \S8, we deal with arbitrary values of the parameter $z$. Here we prove that
the spectral types of the representations $T_z$ are pairwise disjoint (except
the equivalence $T_z\sim T_{\bar z}$).

In \S9 (Appendix), we collected a number of necessary definitions and facts of
general nature.

\subhead 0.6. Concluding remarks \endsubhead The measures
$\mu_t$ involved in the construction of the representations
$T_z$ are very interesting objects in their own right. Each
measure $\mu_t$ can be written as a projective limit
$\varprojlim\mu^n_t$, as $n\to\infty$, where $\mu^n_t$ is a
remarkable probability measure on the finite symmetric group
$S(n)$. These measures $\mu^n_t$ were first discovered in the
context of population genetics and were considered in many
subsequent works (see, e.g., the encyclopedic article
Tavar\'e--Ewens \cite{TE} and references therein). We will call
the measures $\mu_t$ the {\it Ewens measures.} The probability
space $(\frak S, \mu_t)$ is closely related to the {\it Chinese
Restaurant Process\/} construction, see Aldous \cite{Ald},
Pitman \cite{Pit, \S3.1}.

There is a similarity between the spectral decomposition of the
characters $\chi_z$ attached to the representations $T_z$, and
the decomposition of the measures $\mu_t$ into $K$--invariant
ergodic components. In a certain sense, these two problems are
dual to each other: the latter refers to the ``group level''
while the former refers to the ``group dual level''. Moreover,
there exists a general scheme unifying both problems and
providing an interpolation between them, see Kerov \cite{Ker3},
Kerov--Okounkov--Olshanski \cite{KOO}, Borodin--Olshanski
\cite{BO2}. The decomposition of the measures $\mu_t$ is
governed by the {\it Poisson--Dirichlet distributions,\/} see
Kingman \cite{Kin1}, \cite{Kin2}, and also Olshanski \cite{Ol5}.

There is a deep analogy between the infinite symmetric group
$S(\infty)$ and the infinite--dimensional unitary group
$U(\infty)$. This analogy becomes apparent when one compares the
description of characters of both groups, given in the
fundamental papers Thoma \cite{Tho1} and Voiculescu \cite{Voi}.
The theory of harmonic analysis for $S(\infty)$, as developed in
the present paper and the papers of Borodin and Olshanski
mentioned above, also has a counterpart for the group
$U(\infty)$, see Olshanski \cite{Ol6} and Borodin--Olshanski
\cite{BO5}.

In particular, the counterparts of the Ewens measures $\mu_t$
are the so--called {\it Hua--Pickrell measures,\/} see
Borodin--Olshanski \cite{BO4}. More generally, similar measures
can be associated to all 10 infinite series of classical
Riemannian symmetric spaces of compact type, see Neretin
\cite{Ner2}. A pioneer work in this direction is that of
Pickrell \cite{Pic}; our construction of the space $\frak S$ of
virtual permutations was largely influenced by that paper.

As was first observed by Borodin, the expression \tht{**} for
the characters $\chi_z$ can be analytically continued to provide
a {\it complementary series\/} of characters. In the papers of
Borodin and Olshanski, the characters of the complementary
series are considered together with the characters $\chi_z$, the
latter being viewed as the {\it principal series\/} (for a
justification of this terminology, see Okounkov \cite{Ok3}; note
that there also exists a {\it degenerate series}). A natural
question is whether these series of characters (and the
corresponding representations) exhaust all ``reasonable''
objects of harmonic analysis for $S(\infty)$. In this direction
little is known. Using the idea of Borodin \cite{Bor},
Rozhkovskaya \cite{Rozh} obtained an elegant combinatorial
characterization of the characters $\chi_z$ and their analytic
continuation. Kerov \cite{Ker1} considered the so--called {\it
Ewens--Pitman measures\/} generalizing the Ewens measures
$\mu_t$, as possible candidates for an extension of the basic
construction of the representations $T_z$. It would be
interesting to pursue further the study of this question.

\head \S1. The space of virtual permutations \endhead

\subhead 1.1. Canonical projections \endsubhead Let $S(n)$ be
the group of permutations of the finite set $\{1,\,\ldots,n\}$,
the {\it symmetric group of degree} $n$. We identify $S(n)$ with
the subgroup of permutations $s\in S(n+1)$ preserving the last
element $n+1$, i.e., $s(n+1)=n+1$. The inductive limit of groups
$S(n)$ with respect to these embeddings (i.e., the union of
these groups) will be denoted as $S(\infty)=\varinjlim S(n)$.
The elements of $S(\infty)$ are {\it finite} permutations of the
set $\{1,2,\,\ldots\}$, fixing all but finitely many natural
numbers. We call $S(\infty)$ the {\it infinite symmetric group.}

Given a permutation $\tilde s \in S(n+1)$, $n=1,2,\,\ldots$, we define its {\it
derivative permutation\/} (we borrow the term from dynamical systems theory,
see \cite{CFS, chapter 1, \S5}) $s = \tilde s' \in S(n)$ as follows
$$
s(i) = \cases
\tilde s(i),& \text{if } \tilde s(i) \le n,\\
\tilde s(n+1),& \text{if } \tilde s(i) = n+1,
\endcases
$$
where $i=1,\,\ldots,n$. The map $\tilde s \mapsto s$, denoted
$p_{n,n+1}$, will be referred to as the {\it canonical
projection} of $S(n+1)$ onto $S(n)$. Here is an alternative
description of the canonical projection in terms of the cycle
structure of permutations. Depending on the position of the
element $n+1$ in the cycles of the permutation $\tilde s$ we
distinguish between two cases: $n+1$ belongs to a cycle
$(\ldots \to i \to n+1 \to j \to \ldots)$ of length $\ge 2$, or
$n+1$ is a fixed point of $\tilde s$. In the former
case we remove $n+1$ out of its cycle, i.e., we replace this
cycle with $(\ldots \to i \to j \to \ldots)$. In the latter case
$\tilde s$ already belongs to the subgroup $S(n)\subset S(n+1)$, and we
set $s=\tilde s$.

It is clear from the definition that the preimage of a permutation $s \in S(n)$
with respect to the canonical projection contains $n+1$ permutations in
$S(n+1)$. In fact, in order to obtain a permutation $\tilde s \in
p_{n,n+1}^{-1}(s)$ one should insert $n+1$ in a cycle of $s$ right before one
of the elements $j=1,\,\ldots,n$, or take $n+1$ as a new 1--cycle. Note that in
the latter case $\tilde s = s$.

Yet another useful description of the canonical projection
employs the following simple operation on graphs. It will be
convenient to identify permutations with bipartite graphs. We
associate with a permutation $s \in S(n)$ a graph with the vertex
set $\{1,\,\ldots,n;\;1',\,\ldots,n'\}$. Its edges are couples of the
form $(i,j')$, where  $s(i)=j$. The projection
$p_{n,n+1}: \tilde s\mapsto s$ can now be described as follows.
Take the graph of $\tilde s$, and add an extra edge connecting
the vertices $(n+1)$ and $(n+1)'$. Then the graph of the
derivative permutation $s$ arises if one takes for the edges the
paths connecting the vertices in $1,\,\ldots,n$ with the vertices
in $1',\,\ldots,n'$.

Note that the group $S(n)$ acts by left and right multiplications
on both $S(n)$ and $S(n+1)$.

\proclaim{Proposition 1.1.1}
The canonical projection $p_{n,n+1}$ is equivariant with respect
to two--sided action of the group $S(n)$. If $n\ge4$, this is the
only map $S(n+1) \to S(n)$ with this property.
\endproclaim

\demo{Proof} The first claim is immediate from the description of the canonical
projection in terms of bipartite graphs, and the obvious graphical
interpretation of the product of permutations. Assume now that the a map $p:
S(n+1) \to S(n)$ is equivariant. Then, for each permutation $s \in S(n)$, we
have $s^{-1}p(e)s=p(s^{-1}es)=p(e)$, where $e$ is the identity permutation.
Since $e$ is the only central element in $S(n)$ for $n\ge3$, this implies
$p(e)=e$. By the same token, $s^{-1}p((n,n+1))s= p((n,n+1))$ for every
permutation $s \in S(n-1)$. If $n\ge4$, then $e$ is the only element of $S(n)$
commuting with all permutations in $S(n-1)$ (for $n=2,3$ the transposition
$(12)$ also shares this property). Therefore, $n\ge4$ implies $p((n,n+1)) = e$.
Since the group $S(n+1)$ is made of just two double $S(n)$--cosets (the group
$S(n)$ itself and the class containing $(n,n+1)$), it follows that $p =
p_{n,n+1}$. \qed\enddemo

\example{Remark 1.1.2} As it is clear from the proof, there are non--canonical
projections $p:S(n+1)\to S(n)$ for $n=2,3$. They are determined by the
equalities $p((23))=(12)$ and $p((34))=(12)$ respectively. Note that there are
lots of maps $S(n+1) \to S(n)$ which are equivariant with respect to a
one--sided (left or right) action of the group $S(n)$.
\endexample

\subhead 1.2. Virtual permutations \endsubhead
Consider the sequence
$$
S(1) \leftarrow \dots \leftarrow S(n) \leftarrow S(n+1) \leftarrow \dots
$$
of canonical projections, and let
$$
\frak S = \varprojlim S(n)
$$
denote the projective limit of the sets $S(n)$. By definition,
the elements of $\frak S$ are arbitrary sequences $x = (x_n \in
S(n))$, such that $p_{n,n+1}(x_{n+1})=x_n$ for all
$n=1,2,\ldots$. The set $\frak S$ is a closed subset of the
compact space of all sequences $(x_n)$, therefore it is a compact
space itself. Let $p_n: \frak S \to S(n)$ denote the natural
projection, $n=1,2,\ldots$. The group $S(\infty)$ can be
identified with a subset of $\frak S$ via the map $s \mapsto
x=(x_n)$, $x_n=s$ for sufficiently large $n$. In other words,
$S(\infty)\subset\frak S$ consists of the stable sequences $(x_n)$.

The subset $S(\infty)$ is dense in $\frak S$, because $p_n(S(\infty))=p_n(\frak
S)=S(n)$ for any $n$. Hence, the space $\frak S$ is a compactification of the
discrete space $S(\infty)$. The elements of $\frak S$ will be called {\it
virtual permutations} of the set $\{1,2,\ldots\}$.

The definition of the space $\frak S$ does not change if we remove from the
projective limit any finite number of first terms, i.e., if we start from
$S(n)$ instead of $S(1)$, where $n$ is chosen arbitrarily. This simple
observation will be tacitly used in what follows.

In particular, it implies that changing the canonical projections $p_{n,n+1}$
for $n=2,3$ by noncanonical ones (see Remark 1.1.2) does not affect the
construction of the space $\frak S$.

\subhead 1.3.  Realizations of the space of virtual permutations
\endsubhead
We shall use four concrete realizations of the space $\frak S$ of virtual
permutations: \roster \item as of the infinite product $\prod_{n=1}^\infty
\{0,1,\ldots,n-1\}$; \item as of the space of growing trees; \item as of the
space of decreasing maps of the set $\{0,1,\ldots\}$ in itself; \item as of the
space of cyclic structures on the set $\{1,2,\ldots\}$.
\endroster

All the constructions use appropriate realizations of finite sets
$S(n)$, and the canonical projections preserve the specific
structures of those sets. In this sense, the above realizations
are {\it natural}.

\proclaim{Proposition 1.3.1}
There exists a natural homeomorphism $x=(x_n)\mapsto i=(i_n)$
between the space $\frak S$ and the infinite product
$$
I = I_1 \times I_2 \times \ldots \qquad
\text{where } I_n = \{0,1,\,\ldots,n-1\}.
$$
\endproclaim

\demo{Proof}
Given an element $x=(x_n) \in \frak S$, we define the sequence
$i=(i_n) \in I$ as follows. Set $i_1=0$. For every $n=1,2,\ldots$
the coordinate $i_{n+1}$ encodes the relation of $s=x_n$ and
$\tilde s=x_{n+1}$. Specifically, $i_{n+1}=0$ means that $s=\tilde s$,
and $i_{n+1} = j \in \{1,\,\ldots,n\}$ means that the element $n+1$
is inserted in a cycle of $s$ immediately before $j$. One can
easily check that this correspondence is indeed a homeomorphism
$\frak S \to I$, where $I$ is equipped with the product topology.
\qed
\enddemo

Likewise, we get a bijection $S(n)\to I_1\times\dots\times I_n$.
The vector $i(x)=(i_1,i_2,\ldots,i_n)$ is called the {\it code}
of a permutation $x\in S(n)$. We shall apply this terminology to
virtual permutations, too.

We define a {\it finite increasing tree} as a rooted labelled tree with $n+1$
vertices, labelled by the numbers $0,1,\,\ldots,n$ in such a way that the root
has label $0$ and the labels increase along every path leading off the root
(cf. \cite{Sta}, \S1.3). {\it Countable} increasing trees are defined in a
similar way.

\proclaim{Proposition 1.3.2}
There exists a natural bijection between virtual permutations and
countable increasing trees.
\endproclaim

\demo{Proof} Let $\tau$ be a countable increasing tree.  Removing vertices with
the labels $n+1,n+2,\ldots$ we obtain a finite increasing tree which we denote
as $\tau_n$. Thus, a tree $\tau$ can be considered as the union of a chain
$\tau_1 \subset \tau_2 \subset \ldots$ of finite increasing trees, where
$\tau_i$ has $i+1$ vertices. Given a tree $\tau_n$, there are $n+1$ options for
the next tree $\tau_{n+1}$, which are naturally indexed by the numbers
$0,1,\,\ldots,n$. In fact, we can join the vertex $(n+1)$ to every one of the
vertices $0,\,\ldots,n$. Therefore, each countable increasing tree $\tau$ is
determined by the sequence $i=(i_n) \in I$. One can easily check that that the
map $\tau \mapsto (i_n)$ provides a bijection onto $I$. From Proposition 1.3.1
we get a bijection between countable increasing trees, and virtual
permutations. The above construction is quite similar to that of the bijection
between permutations in $S(n)$, and increasing trees with $n+1$ vertices, cf.
\cite{Sta}, \S1.3. \qed
\enddemo

Note that the bijection $\frak S \to I$ identifies finite
permutations $s\in S(\infty) \subset \frak S$ with the sequences
$i=(i_n)$ with only finitely many nonzero coordinates. There is
also an obvious interpretation in terms of increasing trees.

The first two realizations of virtual permutations are actually
very close to each other. The third realization is a version of
the first one, too. We say that a map $\varphi$ of the set
$\{0,1,\,\ldots\}$ into itself is {\it decreasing} if
$\varphi(0)=0$ and $\varphi(n)<n$ for all $n>0$. Setting
$\varphi(n)=i_n$ for $n\ge1$ we obtain a bijective correspondence
between decreasing maps and the points of $I$.

In the next subsection we describe the forth realization of virtual
permutations; its nature is rather different from the first
three.

\subhead 1.4. Cyclic structures \endsubhead Let $J$ be an arbitrary set. We
define a {\it cyclic order} on $J$ as a family of subsets $[i,j)$, called {\it
arcs}, and labelled by ordered pairs $i,j$ of distinct points in $J$. We assume
that the following four axioms hold:
\roster
\item for each pair $i,j$ (where
$i \ne j$) the arcs $[i,j)$, $[j,i)$ do not intersect, and their union is $J$;
\item for each arc $[i,j)$, $i \in [i,j)$ (and hence $j \notin [i,j)$);
\item
if $i,j,k$ are pairwise distinct, then one of the arcs $[i,j)$, $[i,k)$ is a
strict subset of the other;
\item if the arc $[i,j)$ is a strict subset of
$[i,k)$, then $[j,k)=[i,k)\setminus [i,j)$.
\endroster
If $|J|=1$, there is only one cyclic structure on $J$ with the
empty family of arcs.

Note that cyclic orders on a finite set $J$ are in a bijective
correspondence with cyclic permutations of the elements of $J$.
In fact, given a cyclic permutation $s$ of the set $J$, we define
the corresponding cyclic order as follows: for each pair $i\ne j$
the arc $[i,j)$ consists of the points
$i,s(i),\,\ldots,s^{m-1}(i)$ where $m$ is the least number with
$s^m(i)=j$. On the contrary, for infinite sets there is no
natural relation between cyclic orders and permutations.

More generally, we define a {\it cyclic structure} on a set $J$
as a partition of $J$ into nonempty disjoint subsets (called {\it
cycles\/}) with a specified cyclic order on each cycle.
Let $CS(J)$ denote the set of all cyclic structures on $J$. If
$J$ is finite, there is a natural bijection between $CS(J)$ and
the set of all permutations of $J$. Specifically, we  associate with a
permutation $s$ its cycle partition with the cyclic order
on each cycle as defined above.

For every subset $J'\subset J$ there is a natural projection
$CS(J)\to CS(J')$. Moreover, if $J$ is a union of an increasing
sequence of subsets $(J_n)$, then $CS(J)$ can be naturally
identified with the projective limit of the sets $CS(J_n)$. This
observation leads to the following result.

\proclaim{Proposition 1.4.1}
There exists a natural bijection between $\frak S$ and the set
$CS(\{1,2,\,\ldots\})$ of cyclic structures on
$\{1,2,\,\ldots\}$.
\endproclaim

Note that the bijection of Proposition 1.4.1 identifies $p_n: \frak S \to S(n)$
with the projection $CS(\{1,2,\,\ldots\}) \to CS(\{1,\,\ldots,n\})=S(n)$.

We define the {\it cycles of a virtual permutation\/}
$x\in \frak S$ as the cycles of the corresponding cycle structure
on the set $\{1,2,\,\ldots\}$.

In the second realization of $\frak S$ (see \S1.3), the cycles of $x$
correspond to those subtrees of the rooted tree $\tau$ associated
with $x$ which are one edge apart from the root of $\tau$.

In the third realization of $\frak S$, the cycles
correspond to the orbits of the associated map $\varphi$ of the
set $\{1,2,\,\ldots\}$ to itself. Here,  by definition, two numbers
$i,j$ belong to the same $\varphi$--orbit if there exist $k$, $l$,
such that $\varphi^k(i) = \varphi^l(j)$.

\subhead 1.5. The groups $G$ and $K$, and their action on the space $\frak S$
\endsubhead Set $G=S(\infty)\times S(\infty)$. We call $G$ the infinite {\it
bisymmetric\/} group. Let $K=\diag S(\infty)$ denote the diagonal subgroup
$\{(s,s)\in G:\,s\in S(\infty)\}$ in $G$. We shall also use parallel notation
$$
G(n)=S(n)\times S(n), \qquad K(n)=\diag S(n)\subset G(n).
$$

Clearly, the group $K$ is isomorphic to $S(\infty)$. We consider
$S(\infty)$ as a right homogeneous space $K\backslash G$, where
the action of $G$ is defined as follows
$$
s \cdot g = g_2^{-1} s\, g_1, \qquad s \in S(\infty),\; g=(g_1,g_2) \in
G=S(\infty) \times S(\infty).
$$
Since the canonical projections $p_{n,n+1}$ are equivariant, this
action can be naturally extended to the action $\frak S\times
G\to\frak S$ defined as $x\cdot g=y$, where
$$
y_n = x_n \cdot g \quad \text{for all $n$ large enough.}
$$
Specifically, this equality holds whenever $n$ is so large that $g\in
G$ already lies in $G(n)$.

Note that $G$ acts on $\frak S$ by homeomorphisms. Thus, this is the
(only) continuous extension of the action of $G$ on $K\backslash G$ to
the space $\frak S$.

Unfortunately, in all our realizations the description of this action
of $G$ is rather awkward. Only the action of the subgroup $K$ can be
easily described in terms of the forth realization of $\frak S$. Indeed,
the action of $K$ on $\frak S$ is the only continuous extension of  its
action on $S(\infty)$ by conjugations. Under the identification of the
space $\frak S$ with the space of cyclic structures on
$\{1,2,\ldots\}$, this turns into the natural action of $K$ (as a
group of permutations of  $\{1,2,\dots\})$ on the set
$CS(\{1,2,\dots\})$.

\subhead 1.6. The fundamental cocycle \endsubhead
Given $s\in S(n)$, denote by $[s]$ the number of cycles in $s$.
It is important for what follows that the difference $[x\cdot
g]-[x]$ can be defined correctly for all $x \in \frak S$ and $g
\in G$.

\proclaim{Proposition 1.6.1} {\rm(i)} There exists an integer valued function
$c(x,g)$ on $\frak S \times G$, uniquely defined by the following property: if
$n$ is large enough, so that $g \in S(n) \times S(n)$, then
$$
c(x,g) = [p_n(x \cdot g)] - [p_n(x)] =
[p_n(x) \cdot g] - [p_n(x)].
$$

{\rm(ii)} This function is an additive cocycle:
$$
c(x, g_1\,g_2) = c(x,g_1) + c(x \cdot g_1, g_2).
$$

{\rm (iii)} If $g \in K$, then $c(\,\cdot\,,g) \equiv 0$.
\endproclaim

\demo{Proof} (i) It suffices to prove the claim for $g$ of the form $(s,e)$ or
$(e,s)$, where $s$ is a transposition $(ij) \in S(\infty)$. In this case (i)
can be easily checked directly. In fact, let $1\le i<j\le n$. If $i,j$ are both
in the same cycle of $x\in S(n)$, then the multiplication by $(ij)$ from the
left or from the right splits this cycle into two; otherwise the two cycles of
$x$ containing the elements $i$ and $j$ merge into a single cycle of the
product $x\cdot(ij)$ or $(ij)\cdot x$. On the other hand, if
$x=p_{n,n+1}(\widetilde{x})$, where $\widetilde x\in S(n+1)$, then $i,j$ belong
to one and the same cycle of $x$ if and only if they belong to one and the same
cycle of $\widetilde{x}$. Therefore, $[p_n(x)\cdot(ij)]-[p_n(x)]=[(ij)\cdot
p_n(x)]-[p_n(x)]=\pm1$, and this number does not change when $n$ is replaced by
$n+1$.

The claim (ii) follows from (i), and (iii) is obvious.
\qed\enddemo

We call $c(x,g)$ the {\it fundamental cocycle} of the dynamical
system $(\frak S,G)$.

\example{Remark 1.6.2} Let $\overline G\supset G$ and $\overline K\supset K$ be
topological groups as defined above in \S0.1. One can prove that the action of
$G$ on $\frak S$ can be extended to a continuous action $\frak S\times\overline
G\to\frak S$. In particular, the subgroup $\overline K\subset\overline G$ also
acts on $\frak S$. Moreover, all claims of Proposition 1.6.1 hold when $G$ and
$K$ are replaced by $\overline G$ and $\overline K$, respectively.
\endexample

\head 2. Quasiinvariant measures \endhead

\subhead 2.1. The $G$--invariant measure on $\frak S$ \endsubhead Recall that
$\frak S = \varprojlim S(n)$ is the space of virtual permutations, and that we
have defined an action of the bisymmetric group $G=S(\infty)\times S(\infty)$
on $\frak S$. In what follows, all measures are Borel measures.

\proclaim{Proposition 2.1.1} There exists a unique $G$--invariant probability
measure $\mu_1$ on $\frak S$.
\endproclaim

\demo{Proof}
Let $\mu_1^n$ denote the normalized Haar measure on $S(n)$.
Clearly, $\mu_1^n$ coincides with the image of $\mu_1^{n+1}$
under the canonical projection $p_{n, n+1} \colon S(n+1) \to
S(n)$. As $n\to\infty$, the projective limit measure
$$
\mu_1 = \varprojlim \mu_1^n
$$
on $\frak S$ is well defined. It is $G$--invariant, because the measures
$\mu_1^k$, $k \ge n$, are $G(n)$--invariant for all $n=1,2,\ldots$.

Conversely, let $\mu$ be a $G$--invariant probability measure on $\frak S$. For
any $n$, the push--forward of $\mu$ under the projection $p_n\,:\,\frak S \to
S(n)$ should coincide with $\mu_1^n$, the only $G(n)$--invariant probability
measure on $S(n)$. Therefore, $\mu = \mu_1$. \qed\enddemo

Note that the assumption of $G$--invariance in Proposition 2.1.1 can be
replaced by a
weaker assumption of the invariance under the subgroups
$S(\infty) \times \{e \} \subset G$ or
$\{e\} \times S(\infty) \subset G$.

We think of $\mu_1$ as of a substitute of Haar measure for the
group $S(\infty)$.

\subhead 2.2. Ewens measures $\mu_t$ \endsubhead We shall now
include the measure $\mu_1$ into a one--parameter family of
probability measures $\{\mu_t\}_{t>0}$ on $\frak S$.

For $t > 0$ and arbitrary $n=1,2,\ldots$ we define a measure
$\mu_t^n$ on $S(n)$ by the formula
$$
\mu_t^n (\{x\}) = \frac {t^{[x]}} {t(t+1) \ldots (t +n - 1)},
\qquad x \in S(n).
$$
Here, as in \S1.6, $[x]$ stands for the number of cycles in $x$. The measure
$\mu_t^n$ is a probability distribution, as it follows from a well--known
identity (see Stanley \cite{Sta, Proposition 1.3.4})
$$
\sum_{k=1}^n c(n, k) t^k = t (t+1) \ldots (t+n-1).
$$
Here $c(n, k) = \# \{x \in S(n) \colon\; [x]=k\}$ is the absolute
value of the Stirling number of the first kind. Another proof of
this fact follows from

\proclaim{Proposition 2.2.1} Under the  bijection $S(n)\cong
I_1\times\ldots\times I_n$ of\/ {\rm \S1.3},  the measure $\mu_t^n$ turns into
a product measure $\bar \mu_t\,^1\times \bar \mu_t\,^2 \times\, \ldots \times
\bar \mu_t\,^n$. Here $\bar \mu_t\,^m$ is the measure on the set $I_m = \{0,
1,\ldots, m -1 \}$ defined as follows
$$
\bar \mu_t\,^m(\{k\}) = \cases
{1\over t+m-1},& \text{\rm if  $k =1,\,\ldots, m-1$},\\
{t \over t+m-1},& \text {\rm if $k=0$}.
\endcases
$$
\endproclaim

\demo{Proof} Let $i(x)=(i_1,\ldots,i_n)\in I_1 \times\ldots\times I_n$ be the
code of a permutation $x\in S(n)$, and let $l$ denotes the number of 0's among
the coordinates $i_1,\dots,i_n$. By the very definition,
$$
(\bar \mu_t\,^1\times \bar \mu_t\,^2 \times\, \ldots \times \bar
\mu_t\,^n)(i(x))=\frac{t^l}{t(t+1)\dots(t+n-1)}\,.
$$
On the other hand, a coordinate $i_k$
vanishes if and only if the element $k$ creates a new cycle of the
permutation $p_k(x)$, hence the number of cycles $[x]$ in $x$
coincides with the number $l$ of zeros in the vector $i(x)$. This
concludes the proof.
\qed\enddemo

\proclaim{Corollary 2.2.2}
Given $t>0$, the canonical projections $p_{n-1,n}$ preserve
the measures $\mu_t^n$, hence the measure
$$
\mu_t = \varprojlim \mu_t^n
$$
on $\frak S$ is correctly defined. Under the identification $\frak S
\cong I$ of \S1.3, the measure $\mu_t$ looks as the product measure
$$
\bar\mu_t = \bar\mu_t^1\times\bar\mu_t^2\times\,\ldots
$$
\endproclaim

\demo{Proof} Indeed, the canonical projection $p_{n-1,n}$ corresponds
to deleting the last entry $i_n$ of the code $(i_1,i_2,\ldots,i_n)\in
I_1\times\,\ldots\times I_n$. This immediately implies the both
claims. \qed
\enddemo

The measures $\mu^n_t$ on the groups $S(n)$ are known as {\it Ewens measures\/}
(see \S0.6). We will use the same name for the measures $\mu_t$ on the space
$\frak S$, which are built from the measures $\mu^n_t$.

\proclaim{Proposition 2.2.3} For any $t>0$, the Ewens measure $\mu_t$ is
invariant under the action of the group $K$ on $\frak S$.
\endproclaim

\demo{Proof} Indeed, it suffices to prove that for any $n$, the
measure $\mu_t^n$ on $S(n)$ is invariant with respect to the action of
the subgroup
$K(n)\subset K$, isomorphic to $S(n)$. The action under question is
simply the action of $S(n)$ on itself by conjugations. Since
the measure $\mu_t^n$ on $S(n)$ has constant weights on conjugacy
classes, it is invariant.
\qed\enddemo

Since the measures $\mu_t$, $0<t<\infty$, live on a
compact space, it is natural to ask for their limits as
$t$ goes to $0$ or to $+\infty$.

\proclaim{Proposition 2.2.4} There exist weak limits
$$
\underset {t \to 0}\to {\operatorname{w-lim }} \mu_t = \mu_0,\qquad \underset
{t \to +\infty}\to {\operatorname{w-lim }} \mu_t = \mu_\infty.
$$
Here $\mu_0$ is supported by the subset of virtual permutations with a
single cycle while $\mu_\infty$ is the Dirac measure at the point
$e\in S(\infty)\subset\frak S$. Under the identification $\frak S \cong
I$, both $\mu_0$ and $\mu_\infty$ become product measures.
\endproclaim

\demo{Proof} Let us deal with the realization $\frak S\cong I$.
Then it suffices to examine the limit behavior of the measure
$\bar\mu_t^n$ on $I_n=\{0,\dots,n-1\}$, where $n$
is fixed and  $t$ goes to 0 or $\infty$.

When  $t\to0$, the limit exists and is the measure $\bar\mu_0^n$ such
that
$$
\bar\mu_0^n(0)=0, \qquad
\bar\mu_0^n(1)=\dots=\bar\mu_0^n(n-1)=\frac1{n-1}\,.
$$
This means that the finite product measure
$\bar\mu_0^1\times\dots\times\bar\mu_0^n$, being transferred to
$S(n)$, lives on
maximal cycles in $S(n)$. Therefore, the infinite product measure
$\bar\mu_0^1\times\bar\mu_0^2\times\dots$ on $I$, being transferred to
$\frak S$, lives on the virtual permutations with a single cycle in the
sense of \S1.4.

When $t\to\infty$, the measure $\bar\mu_t^n$ tends to the Dirac
measure at $0\in I_n$, which we denote as $\bar\mu_\infty^n$. This
means that the measure on $\frak S$ corresponding to the infinite product
$\bar\mu_\infty^1\times\times\bar\mu_\infty^2\times\dots$ is simply the
Dirac measure at the point $x$ such that $i(x)=(0,0,\dots)$. This
point is just $e$.
\qed
\enddemo

\subhead 2.3. $K$--invariant product measures on $I$ \endsubhead

We can characterize the family of measures $\mu_t$ as follows.

\proclaim{Proposition 2.3.1} The measures $\mu_t$, $0\le t\le\infty$, are
precisely those probability measures on $\frak S$ that are both product
measures (with respect to the identification $\frak S\cong I$ of subsection
1.3) and invariant under $K$.
\endproclaim

\demo{Proof} Let $\mu$ be a $K$--invariant product measure on $\frak S$.
Clearly, $\mu=\varprojlim\mu^n$ with $\mu^n=p_n(\mu)$. We have to show that
$\mu$ coincides with one of the measures $\mu_t$, $0\le t\le\infty$.

Consider the measure $\mu^2$ on $S(2)$. It coincides with some $\mu_t^2$, $0\le
t\le\infty$, and the parameter $t$ is determined uniquely. In fact,
$\mu_\infty^2$ is the Dirac measure at $e\in S(2)$, $\mu_0^2$ is the Dirac
measure at the involution $(1,2)\in S(2)$, and the measure $\mu_t^2$ with
$0<t<\infty$ has weights ${1 \over 1+t}$ and ${t \over 1 + t}$ at the elements
$(1,2)$ and $e=(1)(2)$ respectively. We shall prove by induction in $n\ge2$
that $\mu^n =\mu_t^n$ for all $n$.

We start by considering the induction step $n\to n+1$ in the degenerated cases
of $t=0$ and $t=\infty$. Here we need not the assumption that $\mu$ is a
product measure. In the first case we assume that $\mu^n$ is the uniform
distribution on $n$--cycles in $S(n)$. The preimage of $n$--cycles in $S(n+1)$
under the canonical projection is the set $A\cup B \subset S(n+1)$ where $A$
consists of all $n$--cycles in $S(n)\subset S(n+1)$, and $B$ is the set of all
$(n+1)$--cycles. Each of the two sets is a $\diag S(n)$--orbit, and $B$ (but
not $A$) is also a $\diag S(n+1)$--orbit. Since $\mu^{n+1}$ is $\diag
S(n+1)$--invariant, it is supported by $B$ alone and uniform on $B$.

The case $t=\infty$ is similar. We assume that $\mu^n$ is supported by the
point $\{e\}$. Then the measure $\mu^{n+1}$ is supported by $A\cup B\subset
S(n+1)$ where $A$ consists of transpositions $(n+1,j)$, $1\le j\le n$, and
$B=\{e\}$. Since $A$ is not $S(n+1)$--invariant, $\mu^{n+1}$ is supported by
$B$.

Let us now assume that $\mu^n=\mu_t^n$ where $t\ne0,\infty$. We have to show
that $\mu^{n+1}=\mu_t^{n+1}$. We shall write an arbitrary permutation $x\in
S(n+1)$ as a pair $\{y,\,j\}\in S(n)\times I_{n+1}$ where $y=p_{n,n+1}(x)$ and
$$
j = \cases
0& \text{if $x(n+1) = n+1$};\\
x(n+1)& \text{if $x(n+1) \ne n+1$}.
\endcases
$$
Since $\mu^{n+1}$ is a product measure,
$$
\mu^{n+1}( \{ x \}) = \frac {t^{[y]}\, \nu(j)} {t (t+1)\, \ldots (t + n - 1)},
\tag2-3-1
$$
where $\nu(0),\,\ldots,\nu(n)$ are the weights of a probability
measure on $I_{n+1}$.

Applying a conjugation by appropriate permutation in $S(n)\subset S(n+1)$, one
can replace any $j\ne0$ with any other $j'\ne0$ leaving the cycle structure of
$y$ intact. It follows that $\nu(1)=\,\ldots=\nu(n)$, and we only have to find
out the relation of $\nu(0)$ and $\nu(1)=\dots=\nu(n)$. To this end, choose a
permutation $y\in S(n)$ such that the cycle containing 1 is of length $\ge2$,
i.e., $y(1)=j\ne1$. Then
$$
(1, n+1)\, \cdot \{y,\,0\} \cdot \, (1, n+1) = \{y',\,j\}.
$$
Here $y'$ is obtained from $y$ by removing $j$ from its cycle and forming an
additional trivial cycle $(j)$. An important point is that the initial cycle
containing $j$ does not disappear completely. It follows that $[y']=[y]+1$.
Since, by the invariance assumption, $\mu^{n+1}(\{y,0\})=\mu^{n+1}(\{y',j\})$,
we obtain from \tht{2-3-1}
$$
\nu(0) = t\,\nu(j)=t\,\nu(1).
$$
Since
$$
\nu(0) + n\,\nu(1) = 1,
$$
we obtain
$$
\nu(0) = {t \over t+n}, \qquad \nu(1) = {1 \over t+n},
$$
and the desired identity $\mu^{n+1}=\mu_t^{n+1}$ follows.
\qed\enddemo

\subhead 2.4. Disjointness of the measures $\mu_t$ \endsubhead

\proclaim{Proposition 2.4.1} The measures $\mu_t$, $0\le
t\le\infty$, are disjoint {\rm(}mutually singular{\rm)}.
\endproclaim

\demo{Proof}
We can replace the measures $\mu_t$ on $\frak S$ with the
corresponding product measures
$$
\bar \mu_t = \prod_{n=1}^\infty \bar \mu_t^n
$$
on $I=I_1\times I_2\times\,\ldots$. Note that $\bar\mu_\infty$
is the Dirac measure at the point $i=(0,0,\,\ldots)\in I$, and
$\bar\mu_0$ is supported by the sequences $i=(i_1,i_2,\,\ldots)$
with nonzero coordinates: $i_n\ne0$ for $n\ge2$. Obviously, the
measures $\bar \mu_0$, $\bar \mu_\infty$ and $\bar \mu_t$ are
pairwise disjoint, for every $t\in(0,\infty)$.

Assume now that $s,t\in(0,\infty)$ and $s\ne t$. We shall show
that $\mu_s$ and $\mu_t$ are disjoint, $\mu_s \perp \mu_t$, by
applying the well--known Kakutani criterion \cite{Kak}. To this
end we check that the infinite product $\prod_{n=1}^\infty a_n$,
where
$$
a_n =
\sum_{i=0}^{n-1} \sqrt{\bar \mu_s\,^n (i)\, \bar \mu_t\,^n (i)},
$$
diverges. Set $u=\sqrt{st}$; then
$$
\align
a_n = & \frac {u+n-1} {\sqrt{(s+n-1)\, (t+n-1)}} =\\
=& \frac{1 + {u \over n-1}}
{\sqrt{(1 +{s \over n-1})
(1 + {t \over n-1})}}=\\
= & 1 + \frac {u - {1 \over 2}(s +t)}
 {n-1} + O\left({1 \over (n-1)^2}\right).
\endalign
$$
Now note that
$$
u - {1 \over 2} (s+t) = \sqrt{st} - {1 \over 2}(s + t) \ne 0,
$$
since $s \ne t$. Therefore, the product $\prod_{n=1}^\infty a_n$
is indeed divergent.
\qed\enddemo

\subhead 2.5. Quasiinvariance \endsubhead

\proclaim{Proposition 2.5.1}
Each of the measures $\mu_t$, $0<t<\infty$ is quasiinvariant
with respect to the action of $G$ on the space $\frak S$. More
precisely,
$$
\frac{\mu_t( dx \cdot g)} {\mu_t (dx)} = t^{c(x,g)}; \qquad x \in \frak S, \, g
\in G \tag 2-5-1
$$
where $c(x,g)$ is the fundamental additive cocycle of subsection 1.6. The
measures $\mu_0$ and $\mu_\infty$ are not quasiinvariant.
\endproclaim

\demo{Proof}
It suffices to check that
$$
\mu_t(V \cdot g) = \int_V t^{c(x,g)} \mu_t(dx), \quad g \in G \tag 2-5-2
$$
for every Borel subset $V\subseteq\frak S$. This would imply \tht{2-5-1}, hence
the quasiinvariance of $\mu_t$.

Fix $g\in G$ and and choose $m$ so large that $g\in G(m)$. For arbitrary $n\ge
m$ and $y\in S(n)$, let $V_n(y)\subset\frak S$ denote the preimage of the point
$y\in S(n)$ under the canonical projection $p_n\colon\frak S\to S(n)$; this is
a cylinder set. It suffices to check \tht{2-5-2} for $V=V_n(y)$.

Note that $V_n(y)\cdot g=V_n(y\cdot g)$ and $\mu_t(V_n(y))=\mu_t^n(\{y\})$,
hence
$$
\mu_t(V_n(y)\,  \cdot \, g) = \mu_t^n(\{y  \cdot g \}).
$$
On the other hand,
$$
c(x,g) = [p_n(x  \cdot g)] - [p_n(x)] =
[y  \cdot g] -[y]; \qquad x \in V(y),
$$
so that
$$
t^{c(x,y)} \equiv t^{[y  \cdot g] -[y]}; \qquad x \in V_n(y).
$$

The equation \tht{2-5-2} takes the form
$$
\mu_t^n (\{y  \cdot g\}) =
t^{[y  \cdot g] - [y]}\, \mu_t^n (\{y\}),
$$
which is immediate from the definition of the measure $\mu_t^n$.
\qed\enddemo

Using a well--known trick from ergodic theory, one can replace quasiinvariant
measures $\mu_t$ by invariant, though infinite, measures. In order to do this,
consider the space $\widetilde{\frak S}=\frak S\times\Z$ and define an action
of $G$ on $\widetilde{\frak S}$ as
$$
(x,k) \cdot g = (x \cdot g, \, k + c(x,g)), \qquad x \in \frak S,\quad k \in
\Z, \quad g \in G. \tag 2-5-3
$$
Define the infinite measure $\nu_t$ on $\Z$ as
$$
\nu_t (\{k\}) = t^{-k}, \quad k \in \Z,
$$
and introduce the infinite measure $\tilde\mu_t=\mu_t\times\nu_t$
on $\widetilde{\frak S}$.

\proclaim{Proposition 2.5.2} For every $0<t<\infty$, the measure $\tilde\mu_t$
on the space $\widetilde{\frak S}$ is invariant with respect to the action
\tht{2-5-3} of the group $G$.
\endproclaim

\demo{Proof} This is immediate from \tht{2-5-1} and the definition of the
measure $\nu_t$. \qed\enddemo

This construction will be used below, see \S3.1.

\example{Remark 2.5.3} One can prove that for every $t>0$, the action of each
of the subgroups $S(\infty)\times \{e\}$ and $\{e\}\times S(\infty)$ on the
space $\frak S$ with the measure $\mu_t$ is ergodic and topologically minimal
(the latter means that every orbit is dense in $\frak S$). These claims will
not be used in the sequel.
\endexample

\head 3. Generalized regular representations \endhead

In this Section we introduce a family $\{ T_z \}$ of unitary representations of
the group $G=S(\infty)\times S(\infty)$ parameterized by points
$z\in\C\cup\{\infty\}$ of the Riemann sphere. First we assume $z\in\Bbb
C^*=\C\setminus\{0\}$.

\subhead 3.1. The representations $T_z$ \endsubhead We shall always assume
below that the parameters $t>0$ and $z\in\C^*$ are related as $t=z\bar z$. By
virtue of Proposition 2-5-1,
$$
\frac {\mu_t\, (dx \cdot g)} {\mu_t\, (dx)} = t^{c(x,g)} = |z^{c(x,g)}|^2,
\qquad x \in \frak S, \quad g \in G.
$$
Recall that $c(x,g)$ is an additive cocycle, so that $z^{c(x,g)}$
is a multiplicative one. Therefore, the following formula allows
one to define a unitary representation $T_z$ of the group $G$ in
the Hilbert space $\Cal H=L^2\,(\frak S,\mu_t)$,
$$
\big( T_z(g) \,f \big) (x) = f(x\cdot g)\, z^{c(x,g)},\qquad g \in G,\quad
x\in\frak S, \quad f\in \Cal H.
$$
Note two cases when the multiplier $z^{c(x,g)}$ is equal to 1: \roster \item if
$z=1$ (then the measure $\mu_t=\mu_1$ is invariant); \item if $z$ is arbitrary
but $g\in K\subset G$ (since $c(x,\cdot)\equiv0$ on the subgroup $K$).
\endroster

By the reasons to be made clear later on we call representations
$T_z$ the {\it generalized regular} representations of the group
$G$.

\example{Remark 3.1.1} The above construction is nothing but a specialization
of a well--known general construction. Indeed, to any triple $(\frak X, \Cal G,
\mu)$, where $\Cal G$ is a group acting on the right on a space $\frak X$ with
a quasiinvariant measure $\mu$, one associates a one--parameter family of
unitary representations acting in $L^2(\frak X,\mu)$ according to the formula
$$
(T(g)f)(x)=f(x\cdot g)\left(\frac{\mu(dx\cdot
g)}{\mu(dx)}\right)^{\frac{1+is}2}, \qquad g\in\Cal G, \quad x\in\frak X,\quad
f\in L^2(\frak X,\mu),
$$
where $s\in\R$ is a parameter. When $(\frak X, \Cal G,\mu)=(\frak S, G,\mu_t)$,
we obtain
$$
\left(\frac{\mu(dx\cdot
g)}{\mu(dx)}\right)^{\frac{1+is}2}=\left(t^{c(x,g)}\right)^{\frac{1+is}2}=z^{c(x,g)},
$$
where $z=t^{\frac{1+is}2}$. Hence, we return to the definition of $T_z$.
\endexample

One can introduce representations $T_z$ in a slightly different
way, using the action
$$
g:(x,k)\quad \longmapsto \quad
(x,k)\cdot g = (x\cdot g, \, k + c(x,g))
$$
of the group $G$ on the space $\widetilde{\frak S}=\frak S\times\Z$ with the
infinite invariant measure $\tilde\mu_t=\mu_t\times\nu_t$, see Proposition
2.5.2. There is a natural unitary representation $\widetilde T_t$ of the group
$G$ in the Hilbert space $L^2(\widetilde{\frak S},\tilde\mu_t)$,
$$
(\widetilde T_t(g)f)(x,k)=f((x,k)\cdot g), \qquad g\in G.
$$
Let $\Bbb T$ denote the unit circle $\{\xi\in\C:\, |\xi|=1\}$.

\proclaim{Proposition 3.1.2} For any $t>0$, the representation $\widetilde T_t$
of the group $G$ is unitary equivalent to the direct integral
$$
\int_{\xi\in\Bbb T}  T_{\xi\sqrt{t}} \; d\theta, \qquad \xi=e^{2\pi i\theta},
$$
of the generalized regular representations $T_z$ with $|z|^2=t$.
\endproclaim

\demo{Proof} Define a unitary representation of the group $\Z$ acting in the
space of $\widetilde T_t$ and commuting with the latter representation:
$$
(\widetilde T(l)\,f)(x,k) = f(x,k+l)\,t^{-l/2}, \qquad l\in\Z.
$$
We claim that the decomposition in question is determined by this commuting
representation of $\Z$. To see this, we shall pass to a slightly different
realization of $\widetilde T_t$.

Consider the Hilbert space $L^2(\frak S\times \Bbb T,\mu_t\times d\theta)$,
where $d\theta$ is the normalized Lebesgue measure on the circle $\Bbb T$ (we
again write $\xi=e^{2\pi i\theta}$). In this space, we introduce two commuting
unitary representations of the groups $G$ and $\Z$, as follows
$$
\gathered \big(\widehat T_t(g)\,\widehat f\,\,\big)(x,\xi) =
\widehat f(x\cdot g,\xi)\; (\xi\sqrt{t})^{c(x,g)}, \qquad g\in G,\\
\big(\widehat T(l)\,\widehat f\,\,\big)(x,\xi) = \widehat f(x,\xi)\; \xi^l,
\qquad l\in\Z,
\endgathered
$$
where $\widehat f$ ranges over $L^2(\frak S\times \Bbb T,\mu_t\times d\theta)$.

Clearly, the representation $\widehat T_t$ admits the required decomposition,
which is determined by the action of $\Z$. Hence, it suffices to check that the
representation $\widetilde T_t\times\widetilde T$ in the space
$L^2(\widetilde{\frak S},\widetilde{\mu}_t)$ is unitary equivalent to the
representation $\widehat T_t\times\widehat T$ of the group $G\times\Z$ in the
space $L^2(\frak S\times \Bbb T,\mu_t\times d\theta)$.

The desired unitary equivalence is provided by the transform $\Cal F:\,
f\mapsto\widehat f$,
$$
\widehat{f} (x,\xi) = \sum_{k = - \infty}^\infty f(x,k) t^{-k/2} \xi^{-k}, \tag
3-1-1
$$
which is, in essence, the Fourier transform with respect to the second
argument.

Clearly, $\Cal F$ is an isometry. Let us check that it intertwines both
representations of $G$:
$$ \align \Cal F\big(\widetilde T_t(g)f\big)(x,\xi) &=
\sum_{k = -\infty}^\infty
(\widetilde T_t (g)\, f) \, (x,k)\, t^{-{k \over 2}} \xi^{-k} = \\
&= \sum_{k = -\infty}^\infty f(x \cdot g, \, k + c(x,g))
\,(\xi \sqrt{t})^{-k}=\\
&= \sum_{j = -\infty}^\infty f(x \cdot g, j)
\, (\xi \sqrt{t})^{-j + c(x,g)}= \\
&= \widehat T_t\big(\Cal F(f)\big)(x,\xi).
\endalign
$$

The intertwining property for the action of the group $\Z$ can be checked in a
similar way.

 \qed\enddemo

\subhead 3.2. Admissibility \endsubhead The definition of admissible
representations is given in \S9.9.

\proclaim{Proposition 3.2.1} All the representations $T_z$, $z\in\C^*$, are
admissible representations of the pair $(G,K)$.
\endproclaim

\demo{Proof}
Given $n=1,2,\,\ldots$, consider the canonical projection
$p_n\colon\frak S\to S(n)$. A function $F\circ p_n$, where $F$ is
any function on $S(n)$, will be called a {\it cylinder function}
of {\it level} $n$ on the space $\frak S$. We denote the space of
such functions by $\operatorname{Cyl}^n(\frak S)$, and we call
$$
\operatorname{Cyl}(\frak S) = \underset {n\ge1}\to\cup
\operatorname{Cyl}^n(\frak S)
$$
the {\it space of cylinder functions} on $\frak S$. Clearly,
$$
\operatorname{Cyl}(\frak S) \subset \Cal H=L^2(\frak S,\mu_t)
$$
for every $t>0$.

Note that the canonical projection $S(n)\to S(m)$ is invariant
with respect to conjugations with the elements of the subgroup
$S_m(n)\subset S(n)$, for all $m<n$. It follows that
$$
p_m(x \cdot u) = p_m(x), \qquad x \in \frak S, \quad u \in K_m.
$$
Since the factor $z^{c(x,g)}$ is trivial on the group $K\subset
G$, we derive that
$$
\operatorname{Cyl}^m(\frak S) \subseteq \Cal H_m.
$$
The space $\operatorname{Cyl}(\frak S)$ is clearly dense in $\Cal H=L^2(\frak
S,\mu_t)$, hence $\Cal H_\infty$ is also dense in $\Cal H$, and the
representation $T_z$ is admissible. \qed\enddemo

\example{Remark 3.2.2} The space $\operatorname{Cyl}^m (\frak S)$ is strictly
smaller than $\Cal H_m$. Indeed, the former space is finite--dimensional for
any $m$, while the latter space (as we shall see later on) has infinite
dimension even for $m = 0$.
\endexample

\example{Remark 3.2.3} According to a general result (see \S9.9), Proposition
3.2.1 implies that the representations $T_z$ can be continued to the
topological group $\overline G$. This can be verified directly by making use of
Remark 1.6.2.
\endexample

\subhead 3.3. Approximations by regular representations
\endsubhead
For every $n=1,2,\,\ldots$ we denote by $H^n$ the finite dimensional space $L^2
(S(n))$ defined by the normalized Haar measure $\mu^n$ on $S(n)$, and by
$\Reg^n$ the two--sided regular representation of the group $G(n)=S(n)\times
S(n)$ in this space:
$$
(\Reg^n(g)\, f) (x) = f(g_2^{-1}\, x\, g_1)
$$
where
$$
g = (g_1, g_2) \in G(n), \quad x \in S(n), \quad f \in H^n.
$$
We shall show that every generalized regular representation $T_z$ can be
obtained as an inductive limit of the representations $\Reg^n$ determined by an
appropriate family of isometric embeddings $L_z^n\colon H^n\to H^{n+1}$,
$n=1,2,\ldots$, depending on $z\in\C^*$.

We define the operators $L_z^n\colon H^n\to H^{n+1}$ as follows:
if $f\in H^n$ and $x\in S(n+1)$,
$$
(L_z^n f)(x) = \cases z\; \sqrt{\frac {n+1} {t+n}}\; f(x)
&\text{if $ x \in S(n) \subset S(n+1)$}; \\
\sqrt {\frac {n+1} {t+n}}\; f(p_n(x)) &\text{if $ x \in S(n+1) \setminus
S(n)$}.
\endcases
\tag 3-3-1
$$
Here and below we assume that $t=z\bar z$.

\proclaim{Proposition 3.3.1} For any $z\in\C^*$ the operator $L_z^n$ provides
an isometric embedding $H^n\to H^{n+1}$ which intertwines the
$G(n)$--representations $\Reg^n$ and $\Reg^{n+1}\bigm|_{G(n)}$. Let $T'_z$
denote the inductive limit of the representations $\Reg^n$ with respect to the
embeddings
$$
H^1 \overset{L_z^1} \to\longrightarrow\,
H^2 \overset{L_z^2}  \to\longrightarrow \,
H^3 \overset{L_z^3} \to\longrightarrow\,
\ldots.
$$
Then the representations $T'_z$ and $T_z$ are equivalent.
\endproclaim

\demo{Proof} Note that for every $n=1,2,\,\ldots$ the subspace
$\operatorname{Cyl}^n \subset\Cal H=L^2(\frak S,\mu_t)$ of cylinder functions
of level $n$ is invariant with respect to the operators $T_z(g)$ where $g\in
G(n)$. This follows from the definition of the representation $T_z$, and the
fact that for all $g\in G(n)$ the function $x\mapsto c(x,g)$ is a cylinder
function of level $n$:
$$
c(x, g) = [p_n(x) \cdot g] - [p_n(x)], \qquad x \in \frak S, \quad g \in G(n).
$$

Since the image of the measure $\mu_t$ with respect to the canonical projection
$p_n\colon\frak S\to S(n)$ coincides with $\mu_t^n$, we can identify the
Hilbert spaces $\operatorname{Cyl}^n \subset\Cal H$ and $L^2(S(n),\,\mu_t^n)$.
The operators $T_z(g)\bigm|_{\operatorname{Cyl}^n}$, where $g\in G(n)$, take
the form
$$
(T_z(g)\, f)(x) = f(x \cdot g)\,z^{[x \cdot g] - [x]} \tag 3-3-2
$$
(here $ x \in S(n)$,  $f \in L^2(S(n),\mu_t^n)$.)

Define a function $F_z^n$ on the group $S(n)$ by the formula
$$
F_z^n(x) = \left( \frac {n!} {t(t+1) \, \ldots (t+n -1)} \right)^{1/2}\,
z^{[x]}, \qquad x \in S(n).\tag3-3-3
$$
Then \tht{3-3-2} can be written in the form
$$
(T_z(g)\, f)(x) =
f(x \cdot g)\,\frac {F_z^n(x \cdot g)} {F_z^n(x)}.
$$
Note that the function $\vert F_z^n(x)\vert^2$ coincides with the density of
the measure $\mu_t^n$ with respect to the Haar measure $\mu_1^n$. It follows
that the operator of multiplication by the function $F_z^n$ defines an isometry
$$
\operatorname{Cyl}^n = L^2(S(n),\,\mu_t^n)\, \longrightarrow \,
L^2(S(n),\,\mu_1^n) = H^n
$$
intertwining the representations $T_z\bigm|_{\operatorname{Cyl}^n}$ and
$\Reg^n$ of the group $G(n)$.

Consider now the commutative diagram
$$
\CD
L^2(S(n),\,\mu_t^n)=\C^n  @>>>  \C^{n+1}=L^2(S(n+1),\,\mu_t^{n+1}) \\
  @VVV                                         @VVV \\
   H^n               @>\widetilde L_z^n>>         H^{n+1}
\endCD
$$
where the top arrow denotes the natural embedding (lifting of functions via the
projection $p_{n,n+1}$), the vertical arrows correspond to multiplication by
$F_z^n$ and $F_z^{n+1}$, respectively, and the bottom arrow $\widetilde L_z^n$
is defined by the commutativity requirement. It follows that if $f\in H^n$ and
$x\in S(n+1)$, then
$$
\align
(\widetilde L_z^n\, f) (x)&= F_z^{n+1}(x) \, (F_z^n(p_n(x)))^{-1}\,
f(p_n(x)) =\\
&= \sqrt {\frac {n+1} {t+n}}\, z^{[x] - [p_n(x)]}\, f(p_n(x)).
\endalign
$$
Since
$$
[x] = \cases
[p_n(x)] + 1& \text{if $x\in S(n) \subset S(n + 1)$};\\
[p_n(x)]& \text{if $x \in S(n+1) \setminus S(n)$},
\endcases
$$
we conclude that $\widetilde L_z^n=L_z^n$.

By the very definition, $\widetilde L_z^n$ is an isometric embedding which
commutes with the action of $G(n)$. Hence, the inductive limit representation
$T'_z$ is well defined. From the commutative diagram above we conclude that
$T_z$ and $T'_z$ are equivalent. \qed\enddemo

\subhead 3.4. Representations $T_0$ and $T_\infty$ \endsubhead

\proclaim{Proposition 3.4.1}

{\rm(i)} For every $n=1,2,\,\ldots$, the isometry $L_z^n\colon H^n\to H^{n+1}$
admits a continuous continuation, with respect to the parameter $z\in\C^*$, to
the points $z=0$ and $z=\infty$ of the Riemann sphere $\C\cup\{\infty\}$.
Therefore, the definition of the inductive limit representation $T'_z$ also
makes sense for the values $z=0$ and $z=\infty$.

{\rm(ii)} The representation $T'_\infty$ is equivalent to the natural
two--sided regular representation of the group $G=S(\infty)\times S(\infty)$ on
the Hilbert space $l^2(S(\infty))$.
\endproclaim

\demo{Proof} (i)\; Since $t=z\bar z$, \tht{3-3-1} implies that the limits
$$
L_0^n = \lim_{z \to 0} L_z^n \qquad \text{and} \qquad L_\infty^n = \lim_{|z|
\to \infty} L_z^n
$$
do exist, and have the form
$$
\aligned
\left( L_0^n f \right) (x) &={\cases 0
,&\text{if $x \in S(n) \subset S(n+1)$}, \\
\sqrt{ \frac {n+1} {n}}\, f(p_n(x))
,&\text{if $x \in S(n+1) \setminus S(n)$}, \endcases}\\
\left( L_\infty^n f \right) (x) &={\cases \sqrt {n+1}\, f(p_n(x)) ,&\text{if $x
\in S(n) \subset S(n+1)$}, \\ 0 ,&\text{if $x \in S(n+1) \setminus
S(n)$}.\endcases}
\endaligned
\tag 3-4-1
$$
By continuity, $L_0^n$ and $L_\infty^n$ determine isometric embeddings $H^n\to
H^{n+1}$ commuting with the action of the group $G(n)$. Thus, we can use them
to construct inductive limits $T'_0$, $T'_\infty$ of the two--sided regular
representations of the groups $S(n)\times S(n)$.

(ii)\; For every $n=1,2,\,\ldots$ consider the map $l^2(S(n))\to H^n$ defined
as multiplication by the scalar $\sqrt{n!}$. This is an isometry, since the
counting measure on $S(n)$ equals $n!\,\mu_1^n$. Under identification of both
spaces by this map, the embedding $L_\infty^n$ turns into the natural embedding
$l^2(S(n))\to l^2(S(n+1))$. This completes the proof. \qed\enddemo

Using the identification $T_z=T'_z$ we now may extend the definition of
representations $T_z$ to the values $z=0$ and $z=\infty$ of the parameter $z$.
In this way we obtain a family of representations $T_z$ parametrized by the
points of the Riemann sphere $\C\cup\{\infty\}$. We have shown that our family
forms a continuous deformation of the standard two--sided regular
representation of $G$ in $l^2(S(\infty))$. This is a justification of the term
``generalized regular representation''.

\subhead 3.5. A construction of $T_z$ via representations $T_0$ and $T_\infty$
\endsubhead Let us discuss now the formula for the isometric embeddings
$L^n_z:\,  H^n\to H^{n+1}$. We derived this formula from the initial definition
of the representations $T_z$ with $z\in\C^*$. Then, taking a limit transition
in the formula, we completed the construction of the representations at the
points $z=0$ and $z=\infty$.

Here we aim to show that these two steps can be realized in opposite order. We
start with the definition \tht{3-4-1} of the embeddings $L_0^n$ and
$L_\infty^n$, and then pass to the general operators $L_z^n$. First, we have to
check that \tht{3-4-1} indeed defines isometric embeddings, equivariant with
respect to $G(n)$. For $L^n_\infty$, this immediately follows from the basic
property of the canonical projection. As for $L^n_0$, we observe that up to a
scalar multiple, $L^n_0$ can be defined as lifting along the fibers of the
canonical projection $p_{n,n+1}: S(n+1)\to S(n)$, omitting the natural section
$S(n)\to S(n+1)$. This implies equivariance. To prove the isometry property, we
use the fact that the fiber over any point of $S(n)$ consists just of $n$
points, except a single point belonging to the section.

Next, we note that the spaces $L_0^n(H^n)$ and $L_\infty^n(H^n)$ are mutually
orthogonal subspaces of $H^{n+1}$: the functions in the second space are
supported by the subgroup $S(n)\subset S(n+1)$, and those in the first space
are supported by $S(n+1)\setminus S(n)$.

Comparing \tht{3-3-1} and \tht{3-4-1} we see that for $z\in\C^*$, $L_z^n$ is a
linear combination of $L_0^n$ and $L_\infty^n$:
$$
L_z^n = \frac {\sqrt n} {\sqrt{t+n}} L_0^n + \frac {z} {\sqrt{t+n}} L_\infty^n.
\tag 3-5-1
$$
Moreover, the coefficients of $L_0^n$ and $L_\infty^n$ satisfy
the relation
$$
\left| \frac {\sqrt n} {\sqrt{t+n}} \right|^2 + \left| \frac {z}
{\sqrt{t+n}} \right|^2 = \frac {n + z \bar z} {t + n} = 1.
$$
It follows at once that the operator $L_z^n$ {\it defined\/} by the formula
\tht{3-5-1} is a $G(n)$--equivariant isometry.

Formula \tht{3-5-1} looks very simple and natural. This is an argument in favor
of ``naturalness'' of the representations $T_z$.

\subhead 3.6. Connection between $T_z$ and $T_{-z}$ \endsubhead For any $s\in
S(\infty)$, the number $\operatorname{inv}(s)$ of inversions in $s$ is finite.
Let
$$
\operatorname{sgn}(s)=(-1)^{\operatorname{inv}(s)}=\pm1.
$$
Then $\operatorname{sgn}:\, S(\infty)\to\{\pm 1\}$ is a (unique) nontrivial
one--dimensional representation of the group $S(\infty)$.

\proclaim{Proposition 3.6.1} For every $z\in\C\cup\{\infty\}$, $T_{-z}$ is
equivalent to $T_z\times(\operatorname{sgn}\times \operatorname{sgn})$.
\endproclaim

\demo{Proof} Given $x\in\frak S$ and $g=(g_1,g_2)\in G$, let $n$ be so large
that $g\in G(n)$. Then
$$
c(x,g) = [g_2^{-1}p_n(x)g_1] - [p_n(x)].
$$
It follows that $c(x,g)\in 2\Z$ if the permutation $g_1g_2^{-1}$ is even, and
$c(x,g)\in 2\Z+1$ if $g_1g_2^{-1}$ is odd. Using the initial definition of
$T_z$ (for $z\ne0,\infty$) we derive that
$$
T_{-z}(g) = \operatorname{sgn}(g_1g_2^{-1})\, T_z(g).
$$

When $z=0,\infty$, Proposition 3.6.1 claims the equivalence
$$
T_0 \cong T_0 \otimes (\operatorname{sgn} \times \operatorname{sgn}),\qquad
T_\infty \cong T_\infty \otimes (\operatorname{sgn} \times \operatorname{sgn}).
$$
Such an equivalence is indeed provided by the operator of multiplication by the
function $\operatorname{sgn}(\cdot)$ (we use the realization of the
representations as inductive limits, see subsection 3.4). \qed\enddemo

\head 4. The distinguished spherical function \endhead

\subhead 4.1. The distinguished vector and the coherent system $M_z$ (case
$z\ne0,\infty)$ \endsubhead Assume $z\in\C\setminus\{0\}$ and consider the
generalized regular representation $T_z$. According to the initial construction
of $T_z$ in Hilbert space $\Cal H=L^2(\frak S,\mu_t)$ (see \S3.1), $T_z$ comes
with a {\it distinguished vector\/} $\xi_0$: this vector is simply the function
$f_0\equiv1$ on the space $\frak S$. Clearly, $\xi_0$ is $K$--invariant and has
norm 1.

In the inductive limit realization of $T_z$ as described in \S3.3, the same
vector $\xi_0$ is represented by the functions $F^n_z\in L^2(S(n),\mu^n_1)$
defined in \tht{3-3-3}.

Note that the whole space of $K$--invariant vectors in $\Cal H$ is
infinite--dimensional. However, explicitly constructing invariant vectors other
than the distinguished one is a nontrivial task.

Let $\varphi_z$ denote the spherical function on the group $G$ corresponding to
the distinguished vector $\xi_0$, let $\chi_z$ be the related character of the
group $S(\infty)$, and let $M_z$ be the corresponding coherent system. We will
derive a nice expression for $M_z$. As for $\chi_z$, it seems that it does not
admit a simple explicit expression as a function on the symmetric groups
$S(n)$. In other words, the Fourier coefficients of the functions
$\left.\chi_z\right|_{S(n)}$ are simple whereas the functions themselves are
complex.

Recall a standard notation related to Young diagrams. For a particular box
$b\in\la$ with coordinates $(i,j)$, the number
$$
c(b) = j - i \quad \text{and} \quad
$$
is called the {\it content} of $b$.

\proclaim{Theorem 4.1.1} Let $z\in\C^*$ and $t=z\bar z$. Consider the coherent
system $M_z=\{M^{(n)}_z\}$ as defined above.  For any Young diagram
$\lambda\vdash n$,
$$
M_z(\lambda) =\frac {\prod_{b \in \lambda} |z + c(b)|^2} {t(t+1) \, \ldots (t +
n - 1)} \, \frac {\dim^2 \lambda} {n!} \tag4-1-1
$$
where we abbreviate
$$
M_z(\la)=M^{(n)}_z(\la), \qquad \la\in\Y_n\,.
$$
\endproclaim

\demo{Proof} It will be convenient to identify $T_z$ with the inductive limit
$T'_z$ of regular representations $\Reg^n$. Recall that in this realization the
representation space is defined as the Hilbert completion $H$ of the inductive
limit of finite dimensional Hilbert spaces $H^n = L^2(S(n),\, \mu_1^n)$. The
distinguished vector $f_0$ belongs to $H^1$, hence to all of $H^n$. As an
element of $H^n$ it coincides with the function $F_z^n$ introduced in \S3.
Therefore, for $s\in S(n)$
$$
\left.\chi_z\right|_{S(n)}(s) = \big(\Reg^n(s,e)\, F_z^n, F_z^n \big) =
\frac{1}{n!} \sum_{s_1\in S(n)} F_z^n(s_1\,s)\, \overline{F_z^n(s_1)}.
$$
This can be rewritten as
$$
\left.\chi_z \right|_{S(n)} = (F_z^n)^* \, * F_z^n,
$$
where $f^*(s)=\overline{f(s^{-1})}$ denotes the standard involution on the
group algebra $\C[S(n)]$, and ``$*$'' is the convolution product taken with
respect to normalized Haar measure $\mu_1^n$.

Note that $F_z^n$ is a central function on $S(n)$, hence it can
be decomposed as a sum of characters $\chi^\lambda$,
$$
F_z^n = \sum_{\lambda \vdash n} a(\lambda)\, \chi^\lambda, \tag4-1-2
$$
where $a(\lambda)$ are appropriate complex coefficients. By virtue of the
orthogonality relations,
$$
\left( \chi^\lambda \right)^* \, * \chi^\mu = \delta_{\lambda\mu}
\frac {\chi^\lambda} {\dim \lambda}, \qquad \lambda,\mu \vdash n,
$$
hence
$$
M_z(\lambda) = |a(\lambda)|^2. \tag4-1-3
$$

Recall that
$$
F_z^n(x) = \left( \frac {n!} {t(t+1)\, \ldots (t+n-1)} \right)^{1/2}\, z^{[x]},
\tag4-1-4
$$
where $[x]$ denotes the number of cycles of a permutation $x\in S(n)$, see
\S3.3. We are interested in the decomposition of the central function
$z^{[x]}$.

\proclaim{Lemma 4.1.2} Given $n=1,2,\,\ldots$ and $z\in\C^*$, the decomposition
of the central function $x\mapsto z^{[x]}$ on the group $S(n)$ along the
characters $\chi^\lambda$, $\lambda\vdash n$, can be written in the form
$$
\align
z^{[x]} =& \sum_{\lambda \vdash n} \prod_{b \in \lambda} (z +
c(b)) \cdot \frac {\dim \lambda} {n!}\, \chi^\lambda(x) =\\
& = \sum_{\lambda \vdash n} \prod_{b \in \lambda} \frac {z + c(b)} {h(b)}\,
\chi^\lambda(x). \tag 4-1-5
\endalign
$$
\endproclaim

Keeping together \tht{4-1-2}, \tht{4-1-3}, \tht{4-1-4}, \tht{4-1-5} we get the
desired formula \tht{4-1-1}. Thus, it remains to prove the lemma.

\demo{Proof of the lemma} We switch from central functions on the group $S(n)$
to symmetric functions. This is done using the classical {\it characteristic
map} ``$\operatorname{ch}$'' establishing a bijection between central functions
on $S(n)$ and homogeneous symmetric functions of degree $n$, see \cite{Mac,
1.7}. It is well known that $\operatorname{ch}(\chi^\lambda)=s_\lambda$, hence
we have to prove the formula
$$
\operatorname{ch}(z^{[\cdot]}) = \sum_{\lambda \vdash n}\, \prod_{b \in
\lambda} \frac {z+c(b)} {h(b)}\, \cdot \, s_\lambda, \tag 4-1-6
$$
where $s_\lambda$ are the Schur functions. Let us recall the definition of
$\operatorname{ch}$. If $F$ is a central function on $S(n)$,
$\rho=(1^{k_1}\,2^{k_2}\,\ldots)$ is a partition of $n$, and $x_\rho\in S(n)$
is a permutation of cycle type $\rho$, then
$$
\operatorname{ch} F =
\sum_{\rho \vdash n} z_\rho^{-1} \, F(x_\rho)\, p_\rho.
$$
Here $z_\rho$ is the order of the centralizer of $x_\rho$,
$$
z_\rho =
\frac {n!} {1^{k_1}\,2^{k_2}\,\ldots k_1!\, k_2!\,\ldots},
$$
and $p_\rho=p_1^{k_1}\,p_2^{k_2}\ldots$ are the monomials in the
power sums $p_1,p_2\,\ldots$. Note that
$[x_\rho]=k_1+k_2+\,\ldots$ and $k_1+2k_2+3k_3+\,\ldots=n$.

Denote by $y_1,y_2,\,\ldots$ a sequence of formal variables of
symmetric functions, and let $u$ be still another formal
variable. One can write
$$
\align
1 + \sum_{n \ge 1}\operatorname{ch}(z^{[\cdot]})\, u^n
&= \sum_{(1^{k_1}\, 2^{k_2}\, \ldots)}
\frac {z^{k_1 + k_2 + \, \ldots}\quad u^{1k_1 + 2k_2 + \, \ldots}}
{1^{k_1}\, 2^{k_2} \, \ldots \, k_1!\, k_2! \, \ldots}\, p_1^{k_1}\,
p_2^{k_2} \, \ldots=\\
&=\prod_{n=1}^\infty \sum_{k=0}^\infty
\frac{z^k\, p^k\, u^{nk}}{n^k k!}=
\exp z \left( \frac {up_1} {1} + \frac {u^2p_2} {2} + \,
\ldots \right) =\\
&=\exp z \sum_{i=1}^\infty \left( \frac {uy_i} {1} + \frac
{(uy_i)^2} {2} + \, \ldots \right) =\\
&=\exp \left( -z \sum_{i =1}^\infty \ln (1 - uy_i) \right) =
\prod_{i=1}^\infty (1 - uy_i)^{-z}.
\endalign
$$
The formula \tht{4-1-6} takes the form
$$
\prod_{i =1}^\infty (1 - uy_i)^{-z} = \sum_\lambda \prod_{b \in
\lambda} \frac {z + c(b)} {h(b)} \, \cdot s_\lambda(uy_1, \, uy_2,
\, \ldots),
$$
where $\lambda$ in the right hand side runs over all Young
diagrams. Replacing $uy_i$ with $y_i$, we arrive at the identity
$$
\prod_{i =1}^\infty (1 - y_i)^{-z} =  \sum_\lambda \prod_{b \in \lambda} \frac
{z + c(b)} {h(b)} \, \cdot s_\lambda(y_1, \, y_2, \, \ldots). \tag 4-1-7
$$

Recall that the coefficients of Schur functions in the right hand side are the
polynomials in $z$, hence it suffices to prove \tht{4-1-7} for
$z=N=1,2,\,\ldots$.

It is well known (\cite{Mac, I.3, Example 4}) that
$$
\prod_{b \in \lambda} \frac {N + c(b)} {h(b)} = s_\lambda
(\underbrace{1, \, \ldots, 1}_N)
$$
(this is the dimension of the irreducible representation of the group $GL(\Bbb
C)$ with the highest weight $(\lambda_1,\,\ldots,\lambda_N)$ if
$\lambda_{N+1}=\lambda_{N+2}=\,\ldots=0$, and $0$ otherwise). The formula
\tht{4-1-3} takes the form
$$
\prod_{i=1}^\infty (1 - y_i)^{-N} = \sum_\lambda s_\lambda
\underbrace{(1, \, \ldots, 1)}_N\, s_\lambda (y_1, y_2, \, \ldots).
$$
This is a specialization of a more general identity (\cite{Mac,
Ch.~I, (4.3)})
$$
\prod_{j=1}^\infty \prod_{i=1}^\infty (1 - u_j y_i)^{-1} =
\sum_\lambda s_\lambda(u_1, u_2, \, \ldots)\, s_\lambda(y_1, y_2,
\, \ldots)
$$
where we put $u_1=\ldots=u_N=1$ and $u_{N+1}=u_{N+2}=\,\ldots=0$. This
completes the proof of Lemma 4.1.2 and Theorem 4.1.1. \qed\enddemo
\enddemo

\subhead 4.2. The limit coherent systems $M_0$ and $M_\infty$
\endsubhead
In Theorem 4.1.1 we did not consider the parameter values $z=0$ and $z=\infty$.
However, one can see from \tht{4-1-1} that there exist the limits
$$
M_0 = \lim_{z \to 0} M_z, \quad M_\infty = \lim_{z\to\infty} M_z,
$$
which are also coherent systems on the Young lattice. The coherent system $M_0$
is supported by hook diagrams only (see Proposition 4.3.1 (iii) below), and
$$
M^{(n)}_\infty(\lambda) = \frac {\dim^2 \lambda} {n!}, \qquad \la\in\Y_n, \tag
4-2-1
$$
is the so--called {\it Plancherel measure}.

According to Proposition 9.5.1 these limiting coherent systems
give rise to certain characters $\chi_0$ and $\chi_\infty$ of
the group $S(\infty)$, and to certain spherical functions
$\varphi_0$ and $\varphi_\infty$.

\proclaim{Proposition 4.2.1} The functions $\varphi_0$, $\varphi_\infty$ are
spherical functions of the representations $T_0=T'_0$ and $T_\infty=T'_\infty$
respectively. That is, they coincide with matrix coefficients of certain
$K$--invariant vectors of the representations in question.
\endproclaim

\demo{Proof} In order to see this, examine the behavior of the distinguished
vector $\xi_0$ (recall that as an element of the space $H^n$ it coincides with
the function $F_z^n$) as long as $z\to0$ or $z\to\infty$.

Set $z=r\zeta$ where $r>0$ and $\zeta$ is a point of the unit circle
$S^1\subset\C^*$. If $\zeta$ is fixed, the limits
$$
F_0^n = \lim_{r \to 0} F_{r\zeta}^n, \quad
F_\infty^n = \lim_{r \to \infty} F_{r\zeta}^n, \quad
n = 1, 2, \, \ldots
$$
exist and have the form
$$
\align F_0^n(x)&= \cases \sqrt n \,\zeta &\text{if $[x] = 1$, i.e., if $x$ is a
cycle in $S(n)$ of maximal length $n$},\\ 0 &\text{otherwise};
\endcases\\
F_\infty^n(x)&= \cases
\sqrt {n!}\, \zeta^n &\text {\rm if } x = e, \\
0 &\text{\rm if } x \ne e.
\endcases
\endalign
$$
In these expressions, $\zeta$ enters as a scalar factor only. Hence the
corresponding spherical functions do not depend on the choice of $\zeta$.
Clearly, they coincide with $\varphi_0$ and $\varphi_\infty$, respectively.
\qed
\enddemo

Thus, our definition of the distinguished vector $\xi_0$ (\S4.1) can be
extended to the limit cases $z=0$ and $z=\infty$ --- at least, up to an
unessential scalar factor. Note that the representation $T_\infty$ can be
realized in the Hilbert space $\ell^2(S(\infty)$, and then $\xi_0$ can be
identified with the delta function at $e\in S(\infty)$.

Note also that $\varphi_\infty$ is simply the characteristic function of the
subgroup $K\subset G$, and $\chi_\infty$ is the delta function at the identity
element of $S(\infty)$.

We proceed to analysis of the formula \tht{4-1-1}.

\subhead 4.3. Support of $M_z$
\endsubhead
Here we consider the coherent system $M_z$ with arbitrary
$z\in\C\cup\{\infty\}$. The definition of the support of a coherent system is
given in \S9.4.

\proclaim{Proposition 4.3.1} {\rm(i)} If $z\notin\Z$, then
$\operatorname{supp}(M_z)$ is the whole set $\Y$.

{\rm(ii)} If $\la$ is a nonzero integer, $\la=k$ or $\la=-k$, where
$k=1,2,\,\ldots$, then $\operatorname{supp}(M_z)$ consists of those Young
diagrams $\la$ that have no more than $k$ rows, or, respectively, no more than
$k$ columns.

{\rm(iii)} If $z=0$, then $\operatorname{supp}(M_z)$ is the set of hook
diagrams, i.e., Young diagrams contained inside the union of the first row and
the first column.
\endproclaim

\demo{Proof} (i) If $z\notin\Z\cup\{\infty\}$, then the numerator in
\tht{4-1-1} does not vanish for all $\lambda$, hence $M_z(\lambda)\ne0$. If
$z=\infty$, then it follows from \tht{4-2-1} that $M_\infty(\lambda)\ne0$.

(ii) Assume that $z=k$. Then the zeros in the numerator of \tht{4-1-1}
correspond to the boxes $b=(i,j)\in\lambda$ such that $c(b)=j-i=-k$. These
boxes lie on a diagonal of $\lambda$ passing through the box $(k+1,1)$ in the
first column. The lack of such boxes in $\lambda$ is clearly equivalent to the
fact that $\lambda$ contains $k$ or less rows. In a similar fashion, if $z=-k$
then $-k+c(b)\ne0$ for all $b\in\lambda$ if and only if  $\lambda$ contains no
more than $k$ columns.

(iii) Assume that $z=0$ and consider the limit $M_0=\lim_{z\to0}M_z$. When
$z\to0$, the zero factor $t=z\bar z$ in the denominator of \tht{4-3-1} cancels
with the factor in the numerator corresponding to the box $b=(1,1)$ (every
nonempty diagram $\lambda\ne\varnothing$ contains this box). Other zero factors
in the numerator correspond to the boxes $(2,2),(3,3),\,\ldots$ on the main
diagonal. The absence of such boxes just means that $\lambda$ is a hook.
\qed\enddemo

If a hook diagram $\la$ has {\it arm length} $a=\lambda_1-1$ and {\it leg
length} $l=\lambda'_1-1$, then
$$
M_0(\lambda) =
\frac{|z-l|^2\ldots|z-1|^2|z+1|^2\ldots|z+a|^2}
{(t+1)(t+2)(t+3)\ldots(t+n-1)}\;
\frac{\dim^2\lambda}{n!}.
$$

\subhead 4.4. When the distinguished vector $\xi_0$ is cyclic
\endsubhead

\proclaim{Proposition 4.4.1} Assume  $z\notin\Z$, then the distinguished vector
$\xi_0$ is a cyclic vector of the representation $T_z$.
\endproclaim

\demo{Proof} When $z=\infty$, this is evident from the realization in the space
$\ell^2(S(\infty)$. Assume now $z\in\C\setminus\Z$.

Since $T_z$ is the inductive limit of regular representations $\Reg^n$ and
$\xi_0$ belongs to the spaces $H^n$ of all those representations, it suffices
to check that $\xi_0$ is cyclic in $\Reg^n$ for all $n$. Recall that
$$
\Reg^n \cong \bigoplus_{\lambda \vdash n} (\pi^\lambda \times \pi^\lambda).
$$
Each of the irreducible representations $\pi^\lambda\times\pi^\lambda =
\pi^\lambda\times(\pi^\lambda)^*$ of the group $G(n)=S(n)\times S(n)$ is
spherical with respect to diagonal subgroup $K(n)$, and the corresponding
spherical vector is the function $\chi^\lambda$. Since $\xi_0=F_z^n\in H^n$ is
a $K(n)$--invariant vector, too, it is cyclic if and only if all of its
coefficients in the decomposition in functions $\chi^\lambda$ are nonzero. But
$M_z(\lambda)$ is the square modulus of the coefficient of $\chi^\lambda$,
hence the claim follows from Proposition 4.3.1(i). \qed\enddemo

We shall see below that $\xi_0$ is {\it not} cyclic in $T_z$ if $z\in\Z$.

\subhead 4.5. The equivalence of representations $T_z$ and $T_{\bar z}$
\endsubhead

\proclaim{Proposition 4.5.1} The representations $T_z$ and $T_{\bar z}$ are
unitarily equivalent for all $z\in\C\cup\{\infty\}$.
\endproclaim

\demo{Proof} If $z\in\R$ or $z=\infty$, there is nothing to prove. Hence, we
may assume that $z\notin\R\cup\{\infty\}$, in particular, $z\notin\Z$. By
virtue of Proposition 4.6.1, it suffices to check that $\varphi_z=\varphi_{\bar
z}$ which is equivalent to $\chi_z=\chi_{\bar z}$ and to $M_z=M_{\bar z}$. But
this last formula follows directly from \tht{4-1-1}. \qed\enddemo

Note that it is not immediate from the definition that $T_z$ and
$T_{\bar z}$ are equivalent.

\proclaim{Proposition 4.5.2} Assume that $z\in\C\setminus\Z$.

{\rm(i)} There exists a unique operator $A_z$ intertwining representation $T_z$
and $T_{\bar z}$, and identifying their distinguished spherical vectors.

{\rm(ii)} Realize the representations $T_z$, $T_{\bar z}$ as inductive limits
of representations $\Reg^n$. Then the operator $A_z$ preserves the subspaces
$H^n$, hence determines an operator $A_{n,z}\colon H^n\to H^n$ commuting with
the representation $\Reg^n$, for all $n=1,2,\,\ldots$.

{\rm(iii)} The operator $A_{n,z}$ on the space $H^n=L^2(S(n),\mu_1^n)$ is the
convolution operator with the central function
$$
\Theta_{n,z} = \sum_{\lambda \vdash n} \left( \prod_{b \in
\lambda} \frac {\bar  z + c(b)} {z + c(b)} \right) \, \dim
\lambda \cdot \chi^\lambda.
$$
\endproclaim

\demo{Proof} (i) Follows from the fact that the distinguished vectors are
cyclic (for $z\notin \Z$), and from the coincidence of the corresponding
spherical functions $\varphi_z$, $\varphi_{\bar z}$.

(ii) Follows from (i) and the fact that the distinguished vector (for
$z\notin\Z$) is a $G(n)$--cyclic vector in the representation $\Reg^n$, for all
$n=1,2,\,\ldots$.

(iii) We have to find the operator in $L^2(S(n),\mu_1^n)$ that commutes with
the regular representation $\Reg^n$, and transforms the function
$$
F_z^n = C \sum_{\lambda \vdash n}
\left( \prod_{b \in \lambda} \frac {z + c(b)} {h(b)}
\right) \, \chi^\lambda
$$
into
$$
F_{\bar z}^n = C \sum_{\lambda \vdash n}
\left( \prod_{b \in \lambda} \frac {\bar
z + c(b)} {h(b)} \right) \, \chi^\lambda,
$$
with the same factor $C\ne 0$ (we do not need its precise form at the moment).
Every operator commuting with $\Reg^n$ is a convolution operator with some
central function
$$
\Theta = \sum_{\lambda \vdash n} \theta(\lambda)\, \chi^\lambda.
$$
Note that the convolution operator with the function $\dim\lambda
\cdot\chi^\lambda$ is the projection onto the irreducible component
$\pi^\lambda\times\pi^\lambda$ of the representation $\Reg^n$. It follows that
$\Theta=\Theta_{n,z}$. \qed\enddemo

\subhead 4.6. Reducibility of representations $T_z$ \endsubhead For
$z\in\C\cup\{\infty\}$, let $\widetilde{T}_z$ denote the subrepresentation in
$T_z$ realized in the cyclic span of the distinguished vector. In other words,
$\widetilde{T}_z$ is the cyclic unitary representation of the group $G$
generated by the positive definite function $\varphi_z$. If $z\notin\Z$, then
$\widetilde{T}_z$ coincides with $T_z$ by Proposition 4.4.1; we shall see below
that for $z\in\Z$ it is a proper subrepresentation of $T_z$. Note that
$$
\varphi_1(g) \equiv 1,\quad \varphi_{-1}(g) \equiv \operatorname{sgn}
(g_1g_2^{-1}), \qquad g = (g_1, g_2) \in G,
$$
so that for $z=\pm1$ our representation $\widetilde{T}_z$ is one--dimensional
(more precisely, trivial for $z=1$ and equivalent to
$\operatorname{sgn}\times\operatorname{sgn}$ for $z=-1$). Moreover, for
$z=\infty$ the representation $\widetilde{T}_\infty$ is irreducible, since the
two--sided regular representation $T_\infty$ of the group $G=S\times S$ in
$l^2(S)$ is irreducible. We shall show now that in all other cases our cyclic
representation $\widetilde{T}_z$ (hence the entire representation $T_z$) is
reducible.

\proclaim{Proposition 4.6.1} For every $z\in\C\setminus\{\pm1\}$ the cyclic
representation $\widetilde{T}_z\subseteq T_z$ generated  by the distinguished
vector is reducible.
\endproclaim

\demo{Proof} Let $\chi_z$ denote the character corresponding to the spherical
function $\varphi_z$. Irreducibility of the representation $\widetilde{T}_z$
would imply that $\chi_z$ is an extreme character and hence coincides with a
certain character $\chi^{\alpha\beta}$ from the Thoma list (see the Appendix).
In this case we would also have the equality $M_z=M^{(\alpha,\beta)}$ of the
corresponding coherent systems on the Young lattice. We shall compare the
values of $M_z$ and $M^{(\alpha,\beta)}$ on one--row diagrams $\lambda =(n)$,
$n=1,2,\,\ldots$ and derive that the equality only holds for $z=\pm 1,\infty$.

Indeed, it follows from \tht{4-3-1} that the generating function for
$M_z((n))$, $z\ne0,\infty$, is
$$
1 + \sum_{n \ge 1} M_z((n))\,w^n =
{_2F_1}(z, \bar z; z \bar z; w),
$$
where $w$ is a parameter. On the other hand, the Thoma formula implies that
$$
1 + \sum_{n \ge 1} M^{\alpha\beta} ((n)) w^n = e^{\gamma w} \prod_i \frac {1
+\beta_i w} {1 - \alpha_i w}\,, \qquad
\gamma=1-\sum_{k=1}^\infty(\alpha_k+\beta_k).
$$
Hence, we are led to study the possibility of the equality
$$
_2F_1(z, \bar z; z\bar z; w) = e^{\gamma w} \prod_i \frac {1 +\beta_i w} {1 -
\alpha_i w}. \tag 4-6-1
$$
This formula would also imply that
$$
_2F_1(-z, -\bar z; z \bar z;w) = e^{\gamma w} \prod_i \frac{1 + \alpha_i w} {1
- \beta_i w}, \tag 4-6-2
$$
since the tensor multiplication by the nontrivial one--dimensional
representation $\operatorname{sgn}\times\operatorname{sgn}$ switches $z$ to
$-z$ (by Proposition 3.6.1) and replaces $\chi^{\alpha\beta}$ with
$\chi^{\beta\alpha}$.

It is now easy to see that the equalities \tht{4-8-1}, \tht{4-6-2} are only
possible if $z=\pm1$. In fact, the hypergeometric series converges absolutely
in the open disk $|w|<1$, hence the right hand side of \tht{4-6-1}, \tht{4-6-2}
cannot have poles in this disk. Recall that all the Thoma parameters are
positive, so that there are no cancellations between numerators and
denominators. Therefore, $\alpha_1=1$ or $\beta_1=1$ or $\gamma=1$, and all
other parameters vanish. The first case corresponds to $z=1$, the second one to
$z=-1$, and the last one corresponds to $z=\infty$ and cannot occur for a
finite $z$.

It remains to consider the case $z=0$. In the limit $z\to0$ the generating
function takes the form $1+\sum_{n\ge1}\frac 1nw^n$ and the equalities
\tht{4-6-1}, \tht{4-6-2} cannot hold in this case, too. \qed\enddemo

\subhead 4.7. Transition probabilities \endsubhead Recall from Proposition
4.5.1 that in case $z\notin\Z$ the support of $M_z$ is the entire Young graph.
If $z=\pm k$ where $k=1,2,\,\ldots$, the support is made of the Young diagrams
with $k$ or less rows (columns). If $z=0$, then $M_z$ is supported by the hook
diagrams.

Let $p_z(\lambda,\nu)$ denote the the transition probabilities of the coherent
system $M_z$, see \S9... These quantities are defined for any
$\la\in\supp(M_z)$.

\proclaim{Proposition 4.7.1} Let $\la\in\Y_n$, $\nu\in\Y_{n+1}$,
$\la\nearrow\nu$. We have
$$
\gather p_\infty(\lambda,\nu) = \frac {\dim \nu} {(n +1) \dim \lambda}, \\
p_z(\lambda,\nu) = \frac {|z + c_{\lambda\nu}|^2} {z \bar z +
n}\, p_\infty(\lambda,\nu), \qquad \lambda \in
\Y_n\cap\operatorname{supp}(M_z),\qquad z \ne \infty
\endgather
$$
where $c_{\lambda\nu}=c(\nu\setminus\lambda)$ is the content of the box
$\nu\setminus\lambda$.
\endproclaim

\demo{Proof} Follows immediately from \tht{4-1-1} and the definition of
transition probabilities. \qed\enddemo

\head 5. The commutant and block decomposition \endhead

\subhead 5.1. Simplicity of spectrum \endsubhead Let $T_z$,
$z\in\C\cup\{\infty\}$ be a generalized regular representation. We shall work
with the realization of $T_z$ as inductive limit of the two--sided regular
representation $\Reg^n$ of the group $G(n)=S(n)\times S(n)$. As before, we
denote the space of the representation $\Reg^n$ by $H^n$. The representation
$T_z$ acts in the Hilbert completion $H$ of the space $\bigcup_{n\ge 1} H^n$,
where the maps $L_z^n\colon H^n\to H^{n+1}$ (depending on $z$) were introduced
in \S3. It will be important that $L_z^n$ depends continuously on the parameter
$z$ ranging over the Riemann sphere. Recall that $\Reg^n$ is the direct sum of
irreducible representations $\pi^\lambda\times\pi^\lambda$, $\lambda\in\Bbb
Y_n$ of the group $G(n)$. There are no multiple components in this
decomposition.

Denote by $P_n\colon H\to H^n$ the orthogonal projection onto $H^n$, by
$H(\lambda)\subset H^n$ the space of the representation
$\pi^\lambda\times\pi^\lambda$, and by $P(\lambda)\colon H\to H(\lambda)$ the
orthogonal projection onto $H(\lambda)$. Note that the projectors $P(\lambda)$,
$\lambda\in\Bbb Y_n$, are pairwise orthogonal, and their sum equals $P_n$.

Let $\Cal A$ be the commutant of $T_z$, i.e., the algebra of all
bounded operators in $H$ commuting with the representation
$T_z$. We know that for $z\notin\Z$ the representation $T_z$
admits a cyclic $K$--invariant vector (the distinguished
vector). On the other hand, $(G,K)$ is a Gelfand pair \cite{Ol3,
\S1}. It follows that for $z\notin\Z$, the algebra $\Cal A$ is
isomorphic to the commutant of a commutative operator
$*$--algebra admitting a cyclic vector, whence $\Cal A$ is
commutative. We shall presently give another proof of this fact,
applicable for all $z$.

\proclaim{Proposition 5.1.1} For every $z\in\C\cup\{\infty\}$ the commutant
$\Cal A$ of the representation $T_z$ is a commutative algebra.
\endproclaim

\demo{Proof} We have to prove that $AB = BA$ for arbitrary $A,B\in\Cal A$. For
every $n$, the operators $P_nAP_n$ and $P_nBP_n$ viewed as operators in the
space $H^n$ commute with the representation $\Reg^n$. Since $\Reg^n$
multiplicity free, its commutant is commutative. Therefore,
$$
P_n\, A\, P_n\, B\, P_n = P_n\, B\, P_n\, A\, P_n \qquad \text{for all $n =
1,2,\ldots$}. \tag 5-1-1
$$
Since $H^n\subset H$ form an increasing chain of subspaces and
their union is dense in $H$, the projectors $P_n$ converge to
$1$ strongly as $n\to\infty$. Moreover, the multiplication
operation is continuous in the strong operator topology on every
operator ball. Since the norms of all operators in \tht{5-1-1}
do not exceed the maximum of the numbers $1$, $\Vert A \Vert$,
$\Vert B\Vert$, we can pass to the limit in \tht{5-1-1}, which
gives $AB=BA$. \qed\enddemo

\proclaim{Corollary 5.1.2} For every $z\in\C$, the representation $T_z$ is
decomposable in a multiplicity free direct integral of admissible irreducible
representations of the group $G$.
\endproclaim

\demo{Proof} The existence of a decomposition into a
multiplicity free integral of irreducible representations
follows from the fact that the group $G$ is countable and the
commutant is commutative. The admissibility is easily checked as
in \cite{Ol1, Theorem 3.6}. \qed\enddemo

\subhead 5.2. The transition function of $T_z$
\endsubhead
In order to simplify the notation, we identify the isometric embedding
$L_z^n\colon H^n\to H^{n+1}$ with the partially isometric operator $L_z^nP_n$
in the space $H$. For every $\lambda\in\Bbb Y_n$, $\nu\in\Bbb Y_{n+1}$ (where
$\lambda\nearrow\nu$) fix an isometric embedding $E(\lambda,\nu)\colon
H(\lambda)\to H(\nu)$ commuting with the action of the group $G(n)$. The choice
of $E(\lambda,\nu)$ is unique, up to a complex factor of  modulus 1. We
identify $E(\lambda,\nu)$ with the partially isometric operator
$E(\lambda,\nu)P(\lambda)$ acting in the whole space $H$.

For each pair of Young diagrams $\lambda\in\Bbb Y_n$, $\nu\in\Bbb Y_{n+1}$
consider the operator $P(\nu)\,L_z^n\,P(\lambda)$. This operator intertwines
representations $\pi^\lambda\times\pi^\lambda$ and
$(\pi_\nu\times\pi_\nu)\big\vert_{G(n)}$, hence is $0$ unless
$\lambda\nearrow\nu$. In the latter case it is proportional to
$E(\lambda,\nu)$:
$$
P(\nu)\,L_z^n\,P(\lambda) = \alpha_z(\lambda, \nu) E(\lambda, \nu), \qquad
\lambda\nearrow \nu.
$$
Set
$$
\tilde p_z(\lambda, \nu) = \vert \alpha_z(\lambda, \nu) \vert^2, \qquad \lambda
\nearrow \nu.
$$
Clearly, this function does not depend on the choice of $E(\lambda,\nu)$. It is
also clear that for any $\xi\in H(\lambda)$
$$
||P(\nu)\, \xi||^2 = \cases \tilde p_z(\lambda, \nu)\, ||\xi||^2,&
\text{if $\lambda \nearrow \nu$} \\
0,& \text{otherwise}.
\endcases
$$
Since the projections $P(\nu)$, $\nu\in\Bbb Y_{n+1}$ are pairwise orthogonal
and sum up to $P_{n+1}$, it follows that
$$
\sum_{ \nu\searrow \lambda} \tilde p_z(\lambda,\nu) = 1 \qquad \forall \lambda
\in \Bbb Y_n, \quad n = 1, 2, \, \ldots.
$$
We shall call $\tilde p_z(\lambda,\nu)$ the {\it transition function of the
representation $T_z$.\/} In Theorem 5.5.1 below we show that it coincides with
the transition probabilities of the coherent system $M_z$.

Let us emphasize that the transition function $\tilde p_z(\la,\nu)$ is defined
on the edges of the graph $\Y\setminus\{\varnothing\}$, not the whole Young
graph. That is, we {\it do not\/} attempt to define the value of this function
when $\la$ is the empty diagram $\varnothing$ and $\nu$ is the one--box
diagram.

\subhead 5.3. The commutant in terms of the transition function
\endsubhead
We shall presently show that the transition function of the
representation $T_z$ determines its commutant completely.

Let $\widetilde{\Cal A}$ denote the space of all bounded complex functions
$A(\lambda)$ on the set of vertices of the graph $\Y\setminus\{\varnothing\}$
satisfying the condition
$$
A(\lambda) = \sum_{ \nu \searrow \lambda} \tilde p_z(\lambda,\nu)\, A(\nu),
\qquad \text{for all $\lambda \in \Y$, $\la\ne\varnothing$}, \tag 5-3-1
$$
where $\nu\searrow\la$ means $\la\nearrow\nu$. We consider $\widetilde{\Cal A}$
as a Banach space with the norm $||A||=\sup_\lambda|A(\lambda)|$.

\proclaim{Proposition 5.3.1}
The commutant $\Cal A$ of the generalized regular representation
$T_z$ considered as a Banach space with the ordinary operator
norm is naturally isometric to the space $\widetilde{\Cal A}$.
\endproclaim

\demo{Proof}
For every operator $A\in\Cal A$ and every $n=1,2,\ldots$ we have
$$
P_n\,A\,P_n =
\sum_{\lambda, \mu \in \Bbb Y_n} P(\lambda)\, A\, P(\mu).
$$
But $P(\lambda)\,A\,P(\mu)$ intertwines the representations
$\pi^\mu\times\pi^\mu$ and $\pi^\lambda\times\pi^\lambda$, hence can be nonzero
only if $\lambda=\mu$. In this case the operator $P(\lambda)\,A\,P(\lambda)$
has to be proportional to $P(\lambda)$. Denoting the coefficient by
$A(\lambda)$ we obtain
$$
P_n\,A\, P_n = \sum_{\lambda \in \Bbb Y_n} A(\lambda) P(\lambda).
$$
It is clear that
$$
\Vert P_nAP_n\Vert= \sup_{\lambda \in \Bbb Y_n}\vert A(\lambda) \vert
$$
which implies that
$$
\Vert A\Vert= \sup_n \Vert P_n A P_n \Vert = \sup_{\lambda \in \Bbb Y} \vert
A(\lambda) \vert = \Vert A(\,\cdot\,)\Vert,
$$
where $\Vert A(\,\cdot\,)\Vert$ denotes the sup--norm of the function $A(\la)$.

Let us check now that for any $\lambda\in\Bbb Y_n$, $\nu\in\Bbb Y_{n+1}$,
$$
P(\lambda)\, P(\nu)\, P(\lambda) = \cases \tilde p_z(\lambda, \nu)
P(\lambda),& \text {if $\lambda \nearrow \nu$}, \\
0,& \text {otherwise}.
\endcases
$$
Indeed, the operator $P(\nu)$ commutes with the action of $G(n)\subset G(n+1)$,
hence the operator $P(\lambda)\,P(\nu)\,P(\lambda)$ which commutes with the
irreducible representation $\pi^\lambda\times\pi^\lambda$ of the group $G(n)$
should be proportional to $P(\lambda)$. In order to find the coefficient we
remark that for every $\xi\in H(\lambda)$ we have $P(\lambda)\,\xi=\xi$, so
that
$$
\gathered \left(P(\lambda)\, P(\nu)\, P(\lambda)\, \xi, \xi \right) = (P(\nu)\,
\xi,\xi) = (P(\nu)\, \xi, P(\nu)\, \xi) =\\= \cases \tilde p_z (\lambda, \nu)\,
(\xi,\xi),
& \text {if $\lambda\nearrow \nu$}, \\
0,& \text{otherwise}
\endcases
\endgathered
$$
by the definition of the transition function. On the other hand,
$(P(\lambda)\,\xi,\xi)=(\xi,\xi)$, so that the above coefficient equals $\tilde
p_z(\lambda,\nu)$ if $\lambda\nearrow\nu$, and $0$ otherwise.

Note now that $P(\lambda)P(\nu)P(\mu)=0$ if $\lambda,\mu\in\Bbb Y_n$,
$\mu\ne\lambda$. Therefore,
$$
P_n\, P(\nu)\, P_n = \sum_{\lambda \nearrow \nu} \tilde p_z(\lambda, \nu)
P(\lambda) \qquad \forall \nu \in \Bbb Y_{n+1}.
$$
Setting
$$
A_n = P_n\, A\, P_n, \qquad n = 1, 2, \, \ldots
$$
we see that the property \tht{5-3-1} of the function $A(\cdot)$ simply means
that
$$
A_n = P_n\, A_{n+1}\, P_n, \qquad n = 1, 2, \, \ldots.
$$

Thus, we have constructed above an isometric embedding of $\Cal A$ into the
space $\widetilde{\Cal A}$. In the opposite direction, we shall show that every
function $A(\,\cdot\,)\in\widetilde{\Cal A}$ stems from some operator $A\in\Cal
A$. To this end we set
$$
A_n = \sum _{\lambda \in \Bbb Y_n} A(\lambda)\, P(\lambda),
\qquad n = 1, 2, \, \ldots.
$$
The condition \tht{5-3-1} then implies that
$$
A_n= P_n\, A_{n+1}\, P_n, \qquad n = 1, 2, \, \ldots,
$$
and the condition $\Vert A(\,\cdot\,)\Vert<\infty$ implies
$$
\sup_n \Vert A_n\Vert = \Vert A(\cdot)\Vert < \infty.
$$
It follows that there exists a bounded operator
$$
A = \underset {n \to \infty}\to{\operatorname{w-lim }} A_n
$$
where $\operatorname{w-lim}$ denotes the limit in the weak operator topology.
Since $A_n$ commutes with the action of the group $G(n)$, the operator $A$
belongs to the commutant. It is clear that $A_n$ coincides with $P_n\,A\,P_n$
for all $n$, hence our function $A(\,\cdot\,)$ corresponds to this very
operator. \qed\enddemo

\subhead 5.4. The multiplication in the space $\widetilde{A}$
\endsubhead
Let us denote by $(A\circ B)(\,\cdot\,)$ the multiplication operation of
functions $A(\,\cdot\,)$, $B(\,\cdot\,)$ in $\widetilde{\Cal A}$ corresponding
to the operator multiplication in $\Cal A$. Unfortunately there is no simple
formula for this operation. Nevertheless, it can be described in terms of the
transition function using an appropriate limit procedure.

It will be convenient to extend the definition of the transition function
$\widetilde{p}_z(\lambda,\nu)$ to all pairs of Young diagrams $\lambda$, $\nu$.
Given $\lambda\in\Bbb Y_n$, $\nu\in\Bbb Y_N$, we denote by $\Cal
T(\lambda,\nu)$ the set of  paths
$$
\tau = (\tau_n\nearrow\tau_{n+1}\nearrow\ldots \nearrow\tau_N), \qquad
\tau_n=\la, \quad \tau_N=\nu,
$$
from $\lambda$ to $\nu$ in the Young graph (such paths are commonly called {\it
skew Young tableaux} of shape $\nu/\lambda$). Let
$$
\widetilde{p}_z(\tau) = \prod_{k=n+1}^N \widetilde{p}_z(\tau_{k-1},\tau_k)
$$
denote the transition probability along the path $\tau$, and let
$$
\widetilde{p}_z(\lambda,\nu) = \sum_{\tau\in \Cal T(\lambda,\nu)}
\widetilde{p}_z(\tau)
$$
be the total probability of the transition from $\lambda$ to $\nu$. If
$\lambda$ is not contained in $\nu$, then the set $\Cal T(\lambda,\nu)$ is
empty and $\widetilde{p}_z(\lambda,\nu)=0$. If $\lambda\nearrow\nu$, then the
definition of $\widetilde{p}_z(\lambda,\nu)$ does not change.

Let $C(\Bbb Y_n)$ denote the algebra of functions on the finite set $\Y_n$,
with pointwise multiplication. Given $n<N$, define a linear map
$$
\alpha_{N,n}:\; C(\Bbb Y_N) \to C(\Bbb Y_n)
$$
as follows: if $A(\,\cdot\,)\in C(\Bbb Y_N)$, then
$$
\left( \alpha_{N,n} (A) \right) (\lambda) = \sum_{\nu\in\Y_N}
\widetilde{p}_z(\lambda,\nu)\; A(\nu).
$$

\proclaim{Proposition 5.4.1} Assume that
$A(\,\cdot\,),B(\,\cdot\,)\in\widetilde{\Cal A}$, and let
$$
A_N(\,\cdot\,)=A(\,\cdot\,)\mid_{\Y_N},\qquad
B_N(\,\cdot\,)=B(\,\cdot\,)\mid_{\Y_N}
$$
denote their restrictions to $\Bbb Y_N$. Then
$$
\left( A \circ B \right) (\lambda) = \lim_{N \to \infty} \alpha_{N,n}
(A_N(\,\cdot\,)\, B_N(\,\cdot\,)) (\lambda), \qquad \lambda \in \Bbb Y_n, \tag
5-4-1
$$
\endproclaim

\demo{Proof} Denote by $A,B\in\Cal A$ the operators corresponding to
$A(\,\cdot\,)$, $B(\,\cdot\,)$ and set $A_N=P_NAP_N$, $B_N=P_NBP_N$. Then
$$
A_N\, B_N = \sum_{\nu \in \Bbb Y_N} A_N(\nu)\,B_N(\nu)\,P(\nu),
$$
hence
$$
P(\lambda) A_N B_N P(\lambda) = \alpha_{N,n} (A_N(\,\cdot\,)
B_N(\,\cdot\,))(\lambda)
$$
for all $\lambda\in\Bbb Y_n$. On the other hand, $A_NB_N$
converges strongly to $AB$ as $N\to\infty$, so that the left hand
side converges strongly to
$$
P(\lambda)\, A\, B\, P(\lambda) =
(A \circ B)(\lambda)\, P(\lambda),
$$
which proves \tht{5-4-1}. \qed\enddemo

\subhead 5.5. The identity of the functions $p_z$ and
$\widetilde{p}_z$ \endsubhead

\proclaim{Theorem 5.5.1} The transition function of the representation $T_z$,
$z\in\Bbb C\cup\{\infty\}$, is given by the same expression as the transition
probabilities of the distribution $M_z$, i.e.,
$$
\tilde p_z (\lambda, \nu) = \frac {\vert z + c_{\lambda\nu} \vert^2} {|z|^2 +
n} \cdot \frac {\dim \nu} {(n+1) \dim \lambda}, \tag 5-5-1
$$
where $\lambda\in\Bbb Y_n$, $\nu\in\Bbb Y_{n+1}$, $\lambda\nearrow\nu$ and
$c_{\lambda\nu}=c(\nu\setminus\lambda)$.
\endproclaim

\example{Remark 5.5.2} If $z\in\Z$, the function $p_z(\lambda,\nu)$ is formally
defined on the edges of the proper subgraph $\supp(M_z)\subset\Y$ only, whereas
the function $\tilde p_z (\lambda, \nu)$ is always defined on the whole graph
$\Y$. However, the expression for $p_z(\lambda,\nu)$ given in \S4.7 makes sense
for all couples $\la\nearrow\nu$. This makes it possible to say that both
functions coincide even for $z\in\Z$.
\endexample

\demo{Proof} a) Let us introduce some notation which will also
be used in the sequel. Let $\xi_\lambda$, $\lambda\in\Bbb Y_n$,
be the vector in $H(\lambda)\subset H^n\subset H$ corresponding
to $\chi^\lambda$ under the identification $H^n\cong
L^2(S(n),\mu_1^n)$. It has unit length and is $K(n)$--invariant.
The distinguished vector of $T_z$ will be denoted as $\xi_0$ as
before. As a vector in $L^2(S(n),\mu_1^n)$ it is given by the
function
$$
F_z^n(x) =
\left(\frac{n!}{t(t+1)\,\ldots\,(t+n-1)} \right)^{1/2}\,
z^{[x]}; \qquad x\in S(n),
$$
and can be decomposed as
$$
\xi_0 = \sum_{\lambda \in \Bbb Y_n} a_z(\lambda)\, \xi_\lambda,
$$
where
$$
a_z(\lambda) =
\left( \frac {1} {t(t+1) \, \ldots\, (t+n-1)}
\right)^{1/2} \prod_{b \in \lambda} (z + c(b)) \cdot \frac {\dim
\lambda} {\sqrt {n!}}.
$$

b) We now show that, for every $\nu\in\Bbb Y_{n+1}$, the vector
$\xi_\nu$ can be decomposed as
$$
\xi_\nu = \sum_{\lambda \nearrow \nu} \xi_{\lambda \nu},
$$
where the vectors $\xi_{\lambda\nu}$ are pairwise orthogonal,
$K(n)$--invariant,
$$
\left( \xi_{\lambda\nu}, \xi_{\lambda\nu} \right) =
\frac{\dim\lambda}{\dim\nu}\,,
$$
and $\xi_{\lambda\nu}$ generates, under the action of $G(n)$, the
representation $\pi^\lambda\times\pi^\lambda$.

Indeed, let $H(\pi^\nu)$ be the space of an irreducible representation
$\pi^\nu$ of the group $S(n+1)$, and let $\End\,H(\pi^\nu)$ be the algebra of
operators on this space. Endow $\End\,H(\pi^\nu)$ with the inner product
$$
(A, B) = \frac {\tr (AB^*)} {\dim \nu}
$$
and define an action of the group $G(n+1)$ on $\End\,H(\pi^\nu)$ by
$$
g \cdot A = \pi^\nu(g_1)\, A\, \pi^\nu (g_2)^{-1},\qquad A \in \End\,
H(\pi^\nu), \quad g = (g_1, g_2) \in G(n+1).
$$
It is convenient to identify the vector spaces $\End H(\pi^\nu)$ and $H(\nu)
\subset H^{n+1} = L^2 (S(n+1),\mu^{n+1}_1)$ as follows: to an operator
$A\in\End H(\pi^\nu)$ we assign the function $\widehat
A(s)=\tr(A\pi^\nu(s^{-1})$. The map $A\mapsto\widehat A$ preserves the inner
product and commutes with the action of the group $G(n+1)$.

Let $1_\nu$ denote the identity operator in the space $H(\pi^\nu)$. Its image
under the correspondence $A\mapsto\widehat A$ coincides with the vector
$\xi_\nu$. Further, for any $\la\nearrow\nu$, let $1_{\la\nu}\in\End
H(\pi^\nu)$ denote the orthogonal projection onto the subspace of vectors that
transform according to the representation $\pi^\lambda$ under the action of
$S(n)\subset S(n+1)$. Define $\xi_{\lambda\nu}$ as the image of the operator
$1_{\la\nu}$ under the correspondence $A\mapsto\widehat A$. Clearly, the
vectors $\xi_{\lambda\nu}$ satisfy all the required properties.

c) Now let us remark that it suffices to prove \tht{5-5-1} when
$z\notin\Z$. In fact, the right--hand side of \tht{5-5-1} is
continuous in the parameter $z$ ranging over the Riemann sphere.
It also follows from the definition \tht{5-3-1} of $\tilde
p_z(\lambda,\nu)$ that this function is continuous in $z$, since
so is the map $L_z^n$.

The assumption $z\notin\Z$ implies that $a_z(\la)\ne0$ for all
$\la$, which will be used in the computation below.

d) Equating the decompositions of $\xi_0$ in $\xi_\lambda$, $\lambda\in\Bbb
Y_n$, and in $\xi_\nu$, $\nu\in\Bbb Y_{n+1}$, we conclude that
$$
\xi_0 = \sum_{\lambda \in \Bbb Y_n} a_z(\lambda)\, \xi_\lambda = \sum_{\nu \in
\Bbb Y_{n +1}} a_z(\nu)\, \xi_\nu.
$$
Substituting the decomposition
$\xi_\nu=\sum_{\lambda\nearrow\nu}\xi_{\lambda\nu}$, we arrive at
$$
\sum_{\lambda \in \Bbb Y_n} a_z(\lambda)\, \xi_\lambda = \sum_{\lambda \in \Bbb
Y_n} \sum_{\nu \searrow \lambda} a_z(\nu)\, \xi_{\lambda \nu}.
$$
Comparing the components in both sides that transform according to a given
irreducible representation of $G(n)\subset G(n+1)$ we see that
$$
a_z(\lambda)\, \xi_\lambda = \sum_{\nu \searrow \lambda} a_z(\nu)\,
\xi_{\lambda \nu} \qquad \text{for any $\lambda \in \Bbb Y_n$}.
$$
This implies that
$$
P_\nu\, \xi_\lambda = \frac {a_z(\nu)} {a_z(\lambda)} \xi_{\lambda\nu}, \qquad
\nu \searrow \lambda,
$$
whence
$$
\tilde p_z (\lambda, \nu) = \left|\frac {a_z(\nu)} {a_z(\lambda)}\right|^2
\Vert\xi_{\lambda\nu}\Vert^2 = \left|\frac {a_z(\nu)} {a_z(\lambda)}\right|^2
\frac {\dim \lambda} {\dim \nu}
$$
by the definition of the transition function. This last equation, along with
the explicit formula for the coefficients $a_z(\,\cdot\,)$, implies the desired
formula for $\tilde p_z(\lambda,\nu)$. \qed\enddemo

\subhead 5.6. The subgraphs $\Y(p,q)$ and levels of Young diagrams \endsubhead
Given a couple $(p,q)$ of nonnegative integers, we define a subset
$\Y(p,q)\subset\Y$ as follows
$$
\gather
\Y(p,q) = \{ \lambda \in \Y  \mid (p,q) \in
\lambda, \, (p+1, q+1) \notin \lambda \}, \qquad p,q\ge1\\
\Y(p,0)=\{\la\in\Y\mid \la_{p+1}=\la_{p+2}=\dots=0\}, \qquad p\ge1\\
\Y(0,q)=\{\la\in\Y\mid (\la')_{q+1}=(\la')_{q+2}=\dots=0\}, \qquad q\ge1\\
\Y(0,0)=\{\varnothing\}
\endgather
$$
In other words, the set $\Y(p,q)$, where  $p,q\ge1$, consists of all diagrams
containing the rectangle of shape $p\times q$ but not the box $(p+1,q+1)$. The
set $\Y(p,0)$ consists of all diagrams with at most $p$ rows, and the set
$\Y(0,q)$ consists of all diagrams with at most $q$ columns.

Each $\Y(p,q)$ may be viewed as a connected subgraph of the Young graph.

\proclaim{Proposition 5.6.1} Fix an arbitrary integer $k$ and remove from the
Young graph all edges $\la\nearrow\nu$ such that the content of the box
$\nu\setminus\la$ equals $-k$. Then we obtain a subgraph in $\Y$ whose
connected components are exactly the $\Y(p,q)$'s with $p-q=k$.
\endproclaim

\demo{Proof} This follows from the three claims which are readily checked.
First, the sets $\Y(p,q)$ with fixed $p-q=k$ form a partition of the set of all
Young diagrams. Second, if $\la\nearrow\nu$ and $c(\nu\setminus\la)\ne k$ then
$\la$ and $\nu$ belong to one and the same part $\Y(p,q)$ of that partition.
Third, if $\la\nearrow\nu$ and $c(\nu\setminus\la)=k$ then $\la$ and $\nu$
belong to different parts: specifically, if $\la\in\Y(p,q)$ then
$\nu\in\Y(p+1,q+1)$. \qed
\enddemo

For an arbitrary $k\in\Z$, we define the $k$--{\it level} of a Young diagram
$\lambda$ as follows
$$
\lev_k(\lambda) = \# \{ b \in \lambda \mid c(b) = -k \}.
$$
In other words,  $\lev_0(\lambda)$ equals the length of the main diagonal in
$\la$, and $\lev_k(\lambda)$ is the number of boxes on the diagonal shifted
(with respect to the main diagonal) $k$ boxes downwards if $k\ge0$, and $|k|$
boxes upwards, if $k\le0$.

\proclaim{Proposition 5.6.2} Fix an arbitrary integer $k$. The partition of the
set $\Y$ into disjoint union of the sets $\Y(p,q)$ with $p-q=k$ coincides with
the partition according the value of the $k$--level. Specifically,a given part
$\Y(p,q)$ with $p-q=k$ is exactly the set of diagrams with
$\lev_k(\,\cdot\,)=l$, where $l=\min(p,q)$.
\endproclaim

\demo{Proof} This is evident. \qed
\enddemo

\subhead 5.7. The decomposition into blocks \endsubhead

The knowledge of the transition function $\tilde p_z(\lambda,\nu)$, and the
results of Propositions 5.3.1 and 5.4.1 lead us to a preliminary decomposition
of representations $T_z$ for $z\in\Bbb Z$.

Define a function $A_{pq}(\,\cdot\,)$ on nonempty diagrams as follows. If
$p,q=1,2,\dots$ then this is the characteristic function of the set $\Y(p,q)$.
If $p=1,2,\dots$ and $q=0$ then this is the characteristic function of
$\Y(p,0)\setminus\{\varnothing\}$. Similarly, if $p=0$ and $q=1,2,\dots$ then
this is the characteristic function of $\Y(0,q)\setminus\{\varnothing\}$.

\proclaim{Theorem 5.7.1} Fix $k\in\Z$ and let $(p,q)$ be a couple of
nonnegative integers, not vanishing simultaneously, and such that $p-q=k$.

{\rm(i)} The function $A_{pq}(\,\cdot\,)$ satisfies the condition \tht{5-3-1}
involving the transition function $\tilde p_k(\lambda,\nu)$. Therefore, it
determines an operator $A_{pq}$ in the commutant $\Cal A$ of the representation
$T_k$.

{\rm(ii)} Any operator $A_{pq}$ is an orthogonal projection onto a subspace
$H_{pq}\subset H$. The subspaces $H_{pq}$ are pairwise orthogonal, and their
direct sum is the whole $H$. Thus, they determine a decomposition of the
representation $T_k$ into a direct sum of subrepresentations,
$$
T_k = \underset {p - q = k}\to \oplus T_{pq}\,, \qquad (p,q) \in
\Z_+^2\setminus\{(0,0)\}\,.
$$

{\rm(iii)} Denote by $\Reg_{pq}^n$, where $n$ is large enough, the
subrepresentation of the regular representation $\Reg^n$, which is the union of
the components $\pi^\lambda\times\pi^\lambda$ such that $\lambda\in\Bbb
Y(p,q)$. Let $H_{pq}^n$ be the corresponding subspace in $H^n$. Then the
isometric embedding $L_k^n: H^n\to H^{n+1}$ maps $H_{pq}^n$ into
$H_{pq}^{n+1}$, and the representation $T_{pq}\subset T_k$ coincides with the
inductive limit of the representations $\Reg_{pq}^n\subset\Reg^n$ as
$n\to\infty$.
\endproclaim

\demo{Proof} (i) For $\lambda\in\Bbb Y_n\cap\Bbb Y(p,q)$ we have a chain of
identities
$$
A_{pq} (\lambda) = 1 = \sum _{\nu \searrow \lambda}\tilde
p_k(\lambda,\nu) = \sum _{\nu \searrow \lambda,\,\nu \in
\Y(p,q)}\tilde p_k(\lambda,\nu) = \sum _{\nu \searrow
\lambda}\tilde p_k(\lambda,\nu) A_{pq}(\nu).
$$
Therefore, the function $A_{pq}(\,\cdot\,)$ satisfies the condition
\tht{5-3-1}. Since this function is bounded, it determines, according to
Proposition 5.3.1, an operator $A_{pq}$ in the commutant $\Cal A$ of the
representation $T_k$.

(ii) From Proposition 5.3.1 and the definition of the function
$A_{pq}(\,\cdot\,)$ it follows that for all $n=1,2,\ldots$
$$
P_n A_{pq} P_n =  \sum _{\lambda \in \Bbb Y_n \cap \Bbb Y(p, q)} P(\lambda).
\tag 5-6-2
$$
Since $P_n A_{pq}P_n$ is an orthoprojector for any $n$, so is the operator
$A_{pq}$. A similar argument shows that the projectors $A_{pq}$ are pairwise
orthogonal and sum up to the identity operator. Note that here we could also
use Proposition 5.4.1.

(iii) Since all three operators $P_n$, $A_{pq}$, $P_nA_{pq}P_n$ are
orthoprojectors, it follows that $P_n$ and $A_{pq}$ commute. Along with
\tht{5-6-2} this implies that the operator $P_nA_{pq}=A_{pq}P_n$ projects $H^n$
onto the invariant subspace $H_{pq}^n$. Since $P_{n+1}$ majorizes $P_n$, we see
that $H_{pq}^n$ is a subspace of $H_{pq}^{n+1}$. \qed\enddemo

We shall call the subrepresentations $T_{pq}\subset T_k$ the {\it blocks} of
the representation $T_k$. We shall show below that the blocks are themselves
reducible, and find their decomposition into a direct integral of irreducible
representations. The only exceptions are the one--dimensional blocks
$T_{10}\subset T_1$ and $T_{01}\subset T_{-1}$: each one is generated by the
distinguished spherical vector. The block $T_{10}$ is the identity
representation, and $T_{01}$ is $\operatorname{sgn}\times\operatorname{sgn}$.

\subhead 5.8. Another approach to decomposition into blocks
\endsubhead
It is worth noting that the decomposition of the representation $T_k$, $k \in
\Z$, into blocks can be obtained in a different way -- using the operator $A_z$
intertwining representations $T_z$ and $T_{\bar z}$.

Recall that the operator $A_z$ was described in \S4.7 as inductive limit, as
$n\to\infty$, of certain operators $A_{n,z}\colon H^n\to H^n$ commuting with
the representation $\Reg^n$. In other words, $A_z$ preserves $H^n$ for all
$n=1,2,\ldots$, and the restriction of $A_z$ to $H^n$ is $A_{n,z}$. The
operator $A_{n,z}$ acts in the space $H^n=L^2(S(n),\mu_1^n)$ as the operator of
convolution with the central function
$$
\Theta_{n,z} = \sum_{\lambda \vdash n} \left( \prod_{b \in
\lambda} \frac {\bar z + c(b)} {z + c(b)} \right) \dim \lambda
\cdot \chi^\lambda
$$
on the group $S(n)$ (see Proposition 4.7.1).

This function is correctly defined for all nonintegral values of the parameter
$z$. If $z\in\R\setminus\Z$, then the function $\Theta_{n,z}$ is just
$\delta_e=\sum_{\lambda\vdash n}\dim\lambda\cdot\chi^\lambda$, the delta
function at the identity element of the group $S(n)$. Therefore, the associated
operator is simply the identity operator in $H^n$. Let us see what happens to
$\Theta_{n,z}$ when the parameter $z\in\C\setminus\R$ approaches an integer
point $k\in\Z$. It turns out that there exists a nontrivial limit of
$\Theta_{n,z}$, depending this time on the direction in which $z$ approaches
$k$.

In order to see this, set $z=k+\varepsilon w$, where $w\ne0$ is a fixed complex
number and $\varepsilon$ is a real number going to $0$. Set $u=\bar w\,/\,w$.
Then one can easily see that
$$
\lim_{\varepsilon \to 0} \Theta_{n,z} = \sum_{\lambda \vdash n}
u^{\lev_k(\lambda)} \dim \lambda \cdot \chi^\lambda.
$$
Denote by $A_{n,k}(u)$ the convolution operator with the latter function.
Clearly, for each fixed $u$, $|u|=1$, the sequence of operators $(A_{n,k}(u))$,
$n=1,2,\ldots$, is consistent with the embeddings $L_k^n: H^n\to H^{n+1}$ and
hence determines an operator $A_k(u)$ in the space $H$ of the representation
$T_k$, commuting with this representation.

On the other hand,
$$
A_k(u) \big\vert_{H_{pq}^n} = u^l \cdot1, \qquad l = \min(p,q),
$$
which implies that
$$
A_k(u) = \sum_{p-q=k} u^l A_{pq},
$$
i.e., $A_k(u)$ can be considered as a generating function of the projectors
$A_{pq}$. One can derive from this fact another proof of Theorem 5.6.1, not
relying on Proposition 5.3.1.

\head 6. The invariant vectors \endhead

Let us outline the contents of the present Section.

We start with arbitrary $z$ and describe a convenient
realization of the space $V_z\subset H(T_z)$ formed by
$K$--invariant vectors. Specifically, we show that $V_z$ is
isomorphic to the space $\Cal F_z$ of functions on $\Bbb Y$ that
satisfying two conditions: a harmonicity type condition and a
Hardy type condition. In these terms, if $z$ is an integer, the
splitting of the space $V_z$ induced by the block decomposition
$T_z=\bigoplus_{p-q=z}T_{pq}$ takes especially nice form.

Then we focus on the case when $z$ is an integer. We prove two main results:
Theorems 6.2.1 and 6.2.2. In Theorem 6.2.1 we construct, for any block
$T_{pq}$, a certain $K$--invariant vector $v_{pq}$. In Theorem 6.2.2 we compute
the spectral decomposition of the corresponding spherical function: we show
that the spectral measure lives on a finite dimensional face $\Omega(p,q)$ of
the simplex $\Omega$.

Later on in \S7 we shall show that, for any couple $(p,q)$,  the vector
$v_{pq}$ is a cyclic vector in $T_{pq}$. This allows us to completely
understand the spectral decomposition of $T_z$ at the integer points $z$.

\subhead 6.1. The space $\Cal F_z$
\endsubhead Recall that the vectors $\xi_\la$ were introduced in \S5.5 (see the
proof of Theorem 5.5.1, part a)).

\proclaim{Proposition 6.1.1} For any $n=1,2,\dots$ and any $\lambda\in\Bbb
Y_n$, $\nu\in\Bbb Y_{n+1}$,
$$
(\xi_\lambda, \xi_\nu) = \cases (z+c_{\lambda\nu})/\sqrt{(t+n)(n+1)},
& \text{if $\lambda\nearrow\nu$} \\
0, & \text{ otherwise},
\endcases
$$
where, as before, $c_{\lambda\nu}$ is the content of the box $\nu/\la$,
$t=|z|^2$.
\endproclaim

\demo{Proof} The inner product $(\xi_\lambda,\xi_\nu)$ is a
continuous function in $z$, hence it suffices to prove the
formula under the assumption $z\notin\Z$. In this case
$(\xi_0,\xi_\la)\ne0$ for all $\la$. Indeed, this follows from
the equality  $(\xi_0,\xi_\la)=a_z(\la)$ and the explicit
expression for $a_z(\la)$, see the proof of Theorem 5.5.1, part
a).

Denote by $Q_\lambda$ the projection operator from the Hilbert
space $H$ of the representation $T_z$ to the subspace of all
vectors that transform, under the action of the subgroup $G(n)$,
according to the representation $\pi^\lambda\otimes\pi^\lambda$
(that is, the range of $Q_\la$ is the isotypical component of
$\pi^\lambda\otimes\pi^\lambda$ in $T_z\mid_{G(n)}$). Notice
that $P(\la)\le Q_\la$ and $Q_\la\xi_\la=\xi_\la$. It follows
from the proof of Theorem 5.5.1, part b), that
$$
(Q_\lambda\xi_\nu,Q_\lambda\xi_\nu) = \cases \dim\lambda/\dim\nu,
& \text{if  $\lambda\nearrow\nu$} \\
0, & \text{otherwise}.
\endcases
$$
In particular, if the condition $\lambda\nearrow\nu$ does not hold then
$Q_\lambda\xi_\nu=0$ and $(\xi_\lambda,\xi_\nu)=0$.

Consider the decompositions
$$
\xi_0 = \sum_{\lambda\in\Bbb Y_n} (\xi_0,\xi_\lambda) \xi_\lambda =
\sum_{\nu\in\Bbb Y_{n+1}} (\xi_0,\xi_\nu) \xi_\nu.
$$
Applying the operator $Q_\lambda$ we derive
$$
Q_\lambda\xi_0 = (\xi_0,\xi_\lambda) \xi_\lambda = \sum_{\nu\searrow\lambda}
(\xi_0,\xi_\nu) Q_\lambda\xi_\nu.
$$
Fix a diagram $\nu$ such that $\nu\searrow\lambda$. Taking the inner product
with $\xi_\nu$ we obtain
$$
(\xi_0,\xi_\lambda)\,(\xi_\lambda,\xi_\nu) =
(\xi_0,\xi_\nu)\,(Q_\lambda\xi_\nu,\xi_\nu).
$$
We have used here the relation $(Q_\lambda\xi_{\nu'},\xi_\nu)=0$ for any
$\nu'\in\Y_{n+1}\setminus\{\nu\}$, which in turn follows from the fact that
$Q_\lambda H(\nu)\subset H(\nu)$ for any $\nu\in\Y_{n+1}$.

Now remark that
$$
(Q_\lambda\xi_\nu,\xi_\nu) = (Q_\lambda\xi_\nu,Q_\lambda\xi_\nu) = {\dim\lambda
\over \dim\nu},
$$
and hence
$$
(\xi_0,\xi_\lambda)\,(\xi_\lambda,\xi_\nu) = (\xi_0,\xi_\nu) {\dim\lambda \over
\dim\nu}.
$$
Since $(\xi_0,\xi_\lambda)\ne0$ this implies
$$
(\xi_\lambda,\xi_\nu) = \frac{(\xi_0,\xi_\nu)}
{(\xi_0,\xi_\lambda)}\,
\frac{\dim\la}{\dim\nu}=\frac{a_z(\nu)}{a_z(\la)}\,\frac{\dim\la}{\dim\nu}.
$$
Substituting the explicit expression for $a_z(\,\cdot\,)$ (see
the proof of Theorem 5.5.1 part a)) concludes the proof.
\qed\enddemo

\example{Definition 6.1.2} Denote by $\Cal F_z$ the space of complex--valued
functions $f(\lambda)$ on the vertices $\lambda\ne\varnothing$ of the Young
graph, satisfying the following two conditions.

(i) {\it Pseudoharmonicity\/}: for any $\lambda\in\Bbb Y_n$, $n=1,2,\dots$,
$$
f(\lambda) = \sum_{\nu\searrow\lambda} f(\nu)\, (\xi_\nu,\xi_\lambda) =
\sum_{\nu\searrow\lambda} f(\nu)\, \frac{\bar z +
c_{\lambda\nu}}{\sqrt{(t+n)(n+1)}}
$$
(we have used here the formula of Proposition 6.1.1).
\footnote{Cf. the definition of harmonic functions on $\Y$, see
the end of \S9.3}

(ii) {\it Hardy type condition\/}: for any $n=1,2,\dots$
$$
\Vert f\Vert^2 := \sup_n \sum_{\lambda\in\Bbb Y_n} |f(\lambda)|^2 < \infty.
$$
\endexample

It is worth noting that for any $f$ satisfying the
pseudoharmonicity condition, the sum $\sum_{\lambda\in\Bbb
Y_n}|f(\lambda)|^2 $ does not decrease, as $n\to\infty$ (this
follows from the proof of Proposition 6.1.3 below). This shows
that we could equally well define the norm $\Vert f\Vert$ by the
formula
$$
\Vert f\Vert^2 = \lim_{n\to\infty} \sum_{\lambda\in\Bbb Y_n} |f(\lambda)|^2.
$$

\proclaim{Proposition 6.1.3}
Let $V_z$ be the subspace of $K$--invariant vectors of the
representation $T_z$. Then the map
$$
v \mapsto f_v,\qquad f_v(\lambda) = (v,\xi_\lambda)
$$
provides a linear isomorphism $V_z\to \Cal F_z$, preserving the norm.
\endproclaim

\demo{Proof} Recall that by $P_n$ we denote the orthogonal
projection from $H=H(T_z)$ onto $H^n$. Denote by $V^n$ the
subspace of $K(n)$--invariant vectors in $H^n$. Note that $V^n$
is not contained in $V_z$, but $P_nV_z\subseteq V^n$ and
$P_nV^{n+1}\subseteq V^n$. The vectors $\xi_\lambda$, where
$\lambda\in\Bbb Y_n$, form an orthonormal basis in $V^n$. Given
two vectors
$$
v_n = \sum_{\lambda\in\Bbb Y_n} a(\lambda) \xi_\lambda\in V^n
\qquad \text{ and } \qquad v_{n+1} = \sum_{\nu\in\Bbb Y_{n+1}}
b(\nu) \xi_\nu\in V^{n+1},
$$
we have
$$
v_n = P_n v_{n+1} \quad\Longleftrightarrow\quad a(\lambda) =
\sum_{\nu\searrow\lambda} b(\nu)\, (\xi_\nu,\xi_\lambda) \quad
\forall\,\lambda\in\Bbb Y_n.
$$
Note that these relations imply that
$$
\sum_{\lambda\in\Bbb Y_n} |a(\lambda)|^2 = \Vert v_n\Vert^2 \le \Vert
v_{n+1}\Vert^2 = \sum_{\nu\in\Bbb Y_{n+1}} |b(\lambda)|^2.
$$

Assume now that $v\in V_z$ and $v_n=P_nv$ for $n=1,2,\ldots$.
Then
$$
f_v(\lambda) = (v,\xi_\lambda) = (v_n,\xi_\lambda), \qquad
n = |\lambda|.
$$
The above argument shows that the function $f_v$ satisfies the
pseudoharmonicity condition and
$$
\Vert f_v\Vert^2 := \sup_n \Vert v_n\Vert^2 = \lim_{n\to\infty} \Vert
v_n\Vert^2 = \Vert v\Vert^2 < \infty,
$$
hence $f_v\in \Cal F_z$. Therefore, the function $v\mapsto f_v$ provides an
isometric embedding $V_z\to \Cal F_z$.

Let us show now that the map is surjective. Given an arbitrary $f\in \Cal F_z$,
we set
$$
v_n = \sum_{\lambda\in\Bbb Y_n} f(\lambda)\,\xi_\lambda,
\qquad n = 1,2,\ldots\;.
$$
Then $v_n\in V^n$, $v_n=P_nv_{n+1}$ and
$$
\lim_{n\to\infty} \Vert v_n\Vert^2 = \sup_n \Vert v_n\Vert^2 = \Vert f\Vert^2.
$$
Therefore, the vectors $v_n$ converge to a vector $v\in H(T_z)$.
For every $m$ the vector $v$ is $K(m)$--invariant, because the
vectors $v_n$ possess this property for all $n\ge m$. Hence,
$v\in V_z$ and it follows that $f=f_v$. \qed\enddemo

\subhead 6.2. Two theorems
\endsubhead
We assume from this point to the end of \S6 that $z$ is an
integer and we write $z=k$. Recall that in \S5 we have
introduced invariant subspaces $H(T_{pq})\subset H(T_k)$ (with
the indices $p,q$ subject to the condition $p-q=k$), called the
blocks. Denote by $V_{pq}=V\cap H(T_{pq})$ the subspace of all
$K$--invariant vectors in the block $H(T_{pq})$. Let $\Cal
F_{pq}$ be the subspace in $\Cal F_k$ formed by the functions
supported by $\Bbb Y(p,q)\subset\Bbb Y$. The decomposition
$$
\Cal F_k = \bigoplus_{p-q=k}\Cal F_{pq}
$$
is parallel to
$$
V_k = \bigoplus_{p-q=k} V_{pq},
$$
and the latter corresponds to the decomposition
$$
H(T_k) = \bigoplus_{p-q=k} H(T_{pq}),
$$
of \S5.

Our goal is to produce a certain vector $v_{pq}\in V_{pq}$,
which will be described in terms of the realization $V_{pq}=\Cal
F_{pq}$. The following encoding of a Young diagram
$\lambda\in\Bbb Y(p,q)$ will be convenient. Recall that
$\lambda$ belongs to $\Bbb Y(p,q)$ if and only if $\lambda$
contains the rectangle of shape $p\times q$ ($p$ rows and $q$
columns) but not the box $(p+1,q+1)$. One can represent such a
diagram as union of three parts: the $p\times q$--rectangle, a
diagram $\lambda^+$ to its right, and a diagram below the
rectangle. The transpose of the latter diagram will be denoted
by $\lambda^-$. Formally:
$$
\lambda^+ = (\lambda^+_1,\,\ldots,\lambda^+_p), \qquad \text{where $\lambda^+_i
= \lambda_i - q,\quad 1\le i\le p$};
$$
$$
\lambda^- = (\lambda^-_1,\,\ldots,\lambda^-_q), \qquad \text{where $\lambda^-_j
= (\lambda')_j - p,\quad 1\le j\le q$}.
$$
If $q=0$ then $\lambda^-=\varnothing$ and $\lambda=\lambda^+$. Likewise, if
$p=0$ then $\lambda^+=\varnothing$ and $\lambda=(\lambda^-)'$.

If both $p,q$ are nonzero then the correspondence
$\lambda\mapsto(\lambda^+,\lambda^-)$ establishes a bijection between the sets
$\Bbb Y(p,q)$ and $\Bbb Y(p,0)\times\Bbb Y(q,0)$.

Below we set $n=|\lambda|$.

\proclaim{Theorem 6.2.1} Let $z=k$ be an integer and let $T_{pq}$ be an
arbitrary block of the representation $T_k$ {\rm(}recall that $p,q$ are
nonnegative integers such that $p-q=k$ and $(p,q)\ne(0,0)$ if $k=0${\rm)}.

In the space  $H(T_{pq})$ there is a $K$-invariant vector
$v_{pq}$, such that the corresponding function
$f_{pq}(\lambda)=f_{v_{pq}}(\lambda)$ in the space $\Cal
F_{pq}\subset \Cal F_z$ of functions on the graph $\Bbb Y(p,q)$
has the following form
$$
\gathered f_{pq}(\lambda) = (-1)^{|\lambda^-|}\;\;
\frac{\sqrt{((p-q)^2+n-1)!\, n!}}{(p^2+q^2-pq+n-1)!} \\ \times
\prod_{1\le i<j\le p} (\lambda^+_i-\lambda^+_j+j-i) \prod_{1\le
r<s\le q} (\lambda^-_r-\lambda^-_s+s-r),
\endgathered
$$
where $n=|\la|$.
\endproclaim

In the next theorem we use the concept of  {\it spectral
measure\/}; it is explained in \S9.7. We also need the
finite--dimensional faces $\Om(p,q)\subset\Om$, which are
defined in in \S9.6 ($\Om(p,q)$ is a simplex of dimension
$p+q-1$).

\proclaim{Theorem 6.2.2} Let $v_{pq}$ be the vector defined in Theorem 6.2.1,
let
$$
\varphi_{pq}(g) = \frac{(T_z(g) v_{pq}, v_{pq})}{\Vert v_{pq}\Vert^2} =
\frac{(T_{pq}(g) v_{pq}, v_{pq})}{\Vert v_{pq}\Vert^2}, \qquad g\in G,
$$
be the corresponding spherical function on the group $G$, and let $\si_{pq}$ be
the spectral measure on $\Om$, corresponding to $\varphi_{pq}$.

The  measure $\si_{pq}$ is supported by
$\Omega(p,q)\subset\Omega$ and has the density
$$
\Vert v_{pq}\Vert^{-2}\,\prod_{1\le i<j\le p} (\alpha_i-\alpha_j)^2 \prod_{1\le
r<s\le q} (\beta_r-\beta_s)^2
$$
with respect to Lebesgue measure on $\Omega(p,q)$.
\endproclaim

In order to prove these two theorems, we need a few lemmas.

\subhead 6.3. Preparatory Lemmas \endsubhead

\proclaim{Lemma 6.3.1} Let $p=1,2,\dots$ be fixed. The function
$g_p(\la)$ on the graph $\Y(p,0)$, determined by
$$
g_p(\lambda) =\prod_{1\le i<j\le p} (\lambda_i-\lambda_j+j-i),
$$
satisfies the relation
$$
(p^2+|\lambda|)\, g_p(\lambda) = \sum_{\nu\searrow\lambda} (p+c_{\lambda\nu})\,
g_p(\nu).
$$
\endproclaim

\demo{Proof} By setting
$$
\gather l_i=\lambda_i+p-i \quad (i=1,\dots,p), \qquad |l| =
l_1+\ldots+l_p,\\
V(l_1,\,\ldots,l_p) = \prod_{1\le i<j\le p} (l_i - l_j),
\endgather
$$
we transform the relation to the form
$$
\left(p^2 - {p(p-1) \over 2} + |l|\right)\, V(l_1,\,\ldots,l_p) = \sum_{i=1}^p
(l_i + 1)\, V(l_1,\ldots,l_i+1,\ldots,l_p).
$$
The sum extends over all indices $i=1,\,\ldots,p$, not just
those with $l_{i-1}>l_i+1$, because in case of $l_{i-1}=l_i+1$
the term $V(l_1,\ldots,l_i+1,\ldots,l_p)$ vanishes.

Let us check the latter relation. The right hand side, being a
skew symmetric polynomial in $l_1,\,\ldots,l_p$, is divisible by
$V(l_1,\,\ldots,l_p)$. Since its highest term is
$|l|\,V(l_1,\,\ldots,l_p)$, the right hand side has the form
$$
(|l| + \const) V(l_1,\,\ldots,l_p).
$$
In order to determine the constant, we specialize the identity
$$
\sum_{i=1}^p (l_i +
1)\, V(l_1,\ldots,l_i+1,\ldots,l_p) = (|l| + \const) V(l_1,\,\ldots,l_p).
$$
at
$$
(l_1,l_2,\,\ldots,l_p) = (p-1,p-2,\,\ldots,0),
$$
Then we obtain
$$
p\,V(p,p-2,p-3,\ldots,0) = \left({p(p-1)\over 2} + \const\right)
V(p-1,p-2,\,\ldots,0),
$$
whence
$$
\const = p^2 - {p(p-1) \over 2} = {p(p+1) \over 2},
$$
and we are done.
\qed\enddemo

\example{Remark 6.3.2} One can suggest another proof of Lemma 6.3.1. It is
seemingly more round about, but better explains the origin of the result. The
idea is that in case of $p=k$, $q=0$ we already dispose of an example of a
pseudoharmonic function on the Young lattice $\Bbb Y(p,0)$: the function
$$
\tilde g_p(\lambda) := (\xi_0,\xi_\lambda) =
\left({1\over(p^2)_n}\right)^{1/2}\; {\dim\lambda \over \sqrt{n!}}\;
\prod_{b\in\lambda} (p+c(b)), \qquad n = |\lambda|
$$
(cf. the proof of Theorem 5.5.1, part a); in this setup
$t=p^2$). Since
$$
{\dim\lambda \over n!} = \frac{V(l_1,\,\ldots,l_p)}{l_1! \ldots
l_p!}
$$
and
$$
\prod_{b\in\lambda} (p+c(b))=\frac{l_1!\dots l_p!}{(p-1)!\dots
0!}\,,
$$
we obtain that
$$
\tilde g_p(\lambda) = \frac1{(p-1)!\dots
0!}\,\Big(\frac{n!}{(p^2)_n}\Big)^{1/2}\, g_p(\lambda),
$$
and the required relation for $g_p(\lambda)$ follows from the pseudoharmonicity
condition for $\tilde g_p(\lambda)$. \qed
\endexample

Using Lemma 6.3.1 and the fact that the graph $\Y(p,q)$ is
isomorphic to the product $\Y(p,0)\times\Y(q,0)$ we will produce
a pseudoharmonic function on $\Y(p,q)$.

Let $(x_n)_{n\ge pq}$ be a sequence of positive real numbers
satisfying the recurrence relation
$$
{x_{n+1} \over x_n} =
\frac{n+p^2+q^2-pq}{\sqrt{(n+(p-q)^2)(n+1)}}.
$$
For instance, one can set
$$
x_n = \frac{\Gamma(p^2+q^2-pq+n)} {\sqrt{\Gamma((p-q)^2+n)\Gamma(n+1)}}. \tag
6-3-1
$$

Below we use the correspondence $\la\mapsto(\la^+,\la^-)$
introduced in \S6.2.

\proclaim{Lemma 6.3.3} Let $p,q$ be nonnegative integers,
$(p,q)\ne(0,0)$. The function $f_{pq}(\lambda)$ on the graph
$\Bbb Y(p,q)$ determined by the formula
$$
f_{pq}(\lambda) = {(-1)^{|\lambda^-|} \over x_n}
\prod_{1\le i<j\le p} (\lambda^+_i-\lambda^+_j+j-i)
\prod_{1\le r<s\le p} (\lambda^-_r-\lambda^-_s+s-r),
$$
satisfies the pseudoharmonicity condition with the parameter
$k=p-q$, i.e.,
$$
f_{pq}(\lambda) = \sum_{\nu\searrow\lambda}
\frac{p-q+c_{\lambda\nu}}{\sqrt{((p-q)^2+n)(n+1)}} f_{pq}(\nu)
$$
for any $\la\in\Y(p,q)\cap\Y_n$.
\endproclaim

In this formula we assume that $\nu$ belongs to $\Y(p,q)$.
However, the formula remains true without this assumption,
because if $\la\nearrow\nu$, $\la\in\Y(p,q)$, and
$\nu\notin\Y(p,q)$ then the factor $p-q+c_{\la\nu}$ vanishes.

\demo{Proof} Assume first that both $p$ and $q$ are positive. In the notation
of Lemma 6.3.1,
$$
f_{pq}(\lambda) = {(-1)^{|\lambda^-|} \over x_n}\, g_p(\lambda^+)\,
g_q(\lambda^-).
$$
Let $\nu\in\Bbb Y(p,q)$ denote a Young diagram with $n+1$ boxes. The condition
$\lambda\nearrow\nu$ implies one of the following two conditions:

(i)  $\lambda^+\nearrow\nu^+$, \qquad $\lambda^-=\nu^-$;

(ii)  $\lambda^+=\nu^+$, \qquad $\lambda^-\nearrow\nu^-$.

Note that
$$
c_{\lambda\nu} = \cases
q + c_{\lambda^+\nu^+}, & \text{ in case (i)}, \\
-p - c_{\lambda^-\nu^-}, & \text{ in case (ii)}.
\endcases
$$
Therefore,
$$
p - q + c_{\lambda\nu} = \cases
p + c_{\lambda^+\nu^+}, & \text{ in case (i)}, \\
-(q + c_{\lambda^-\nu^-}), & \text{ in case (ii)}.
\endcases
$$
It follows that
$$
\gathered \sum_{\nu\searrow\lambda} (p-q+c_{\lambda\nu})\, f_{pq}(\nu) = \\
= \sum_{\nu^+\searrow\lambda^+} {(-1)^{|\nu^-|} \over x_{n+1}}\,
(p+c_{\lambda^+\nu^+})\, g_p(\nu^+)\, g_q(\lambda^-)  \\ +
\sum_{\nu^-\searrow\lambda^-} {(-1)^{|\nu^-|+1} \over x_{n+1}}\,
(q+c_{\lambda^-\nu^-})\, g_p(\lambda^+)\, g_q(\nu^-).
\endgathered
$$
In the first sum $|\nu^-|=|\lambda^-|$, and in the second one
$|\nu^-|=|\lambda^-|+1$. Hence, the sign in both formulas is
$(-1)^{|\lambda^-|}$.

Applying the Lemma 6.3.1, we derive that the last expression equals
$$
\gathered
{(-1)^{|\lambda^-|} \over x_{n+1}}
(p^2 + |\lambda^+| + q^2 + |\lambda^-|)
g_p(\lambda^+)\, g_q(\lambda^-) = \\ =
{x_n \over x_{n+1}} (p^2+q^2+n-pq)\, f_{pq}(\lambda).
\endgathered
$$
We have shown that
$$
\sum_{\nu\searrow\lambda} (p-q+c_{\lambda\nu})\, f_{pq}(\nu) = {x_n \over
x_{n+1}} (p^2+q^2+n-pq)\, f_{pq}(\lambda).
$$
Dividing both sides by $\sqrt{(n+(p-q)^2)(n+1)}$ we obtain that
the coefficient of $f_{pq}(\lambda)$ in the right--hand side
equals
$$
{x_n \over x_{n+1}} \frac{p^2+q^2+n-pq}
{\sqrt{(n+(p-q)^2)(n+1)}}\,,
$$
which is $1$ by definition of the numbers $x_n$. The lemma follows.

If $q=0$ or $p=0$, our relation is simply equivalent to Lemma
6.3.1.
\qed\enddemo

\proclaim{Lemma 6.3.4} The function $f_{pq}(\lambda)$ introduced
in Lemma 6.3.3 satisfies the Hardy type condition of Definition
6.1.2. Hence, $f_{pq}\in \Cal F_{pq}$.
\endproclaim

\demo{Proof} In order to simplify the notation, set
$a=\lambda^+$, $b=\lambda^-$, i.e.,
$$
a_1=\lambda^+_1, \ldots, a_p=\lambda^+_p; \qquad
b_1=\lambda^-_1, \ldots, b_q=\lambda^-_q.
$$
When $\lambda$ ranges over $Y(p,q)\cap\Y_n$, the couple $(a,b)$
ranges over the subset in $\Z_+^{p+q}$ determined by the
conditions
$$
a_1\ge\ldots\ge a_p,\qquad
b_1\ge\ldots\ge b_q,\qquad
\sum a_i + \sum b_j = n - pq.
$$

We now have
$$
\underset{|\lambda|=n}\to{\sum_{\lambda\in\Bbb Y(p,q)}}
|f_{pq}(\lambda)|^2 = {1\over x_n^2} \sum_{a,b}
\prod_{1\le i<j\le p} (a_i-a_j+j-i)^2 \,
\prod_{1\le r<s\le p} (b_r-b_s+s-r)^2.
$$
Assuming here $n>pq$, we set
$$
\bar\iota_n(\la)=(\al;\be)=(\al_1,\,\dots,\al_p;\be_1,\,\dots,\be_q),
$$
where
$$
\al_1={a_1\over n-pq}, \dots, \al_p={a_p\over n-pq};\qquad \be_1={b_1\over
n-pq}, \dots, \be_q={b_q\over n-pq}.
$$
The vector $(\al;\be)$ belongs to the lattice
$\frac1{n-pq}\,\Z_+^{p+q}$ and satisfies the conditions
$$
\al_1\ge\dots\ge\al_p,\qquad \be_1\ge\dots\ge\be_q,\qquad \sum
\al_i + \sum \be_j = 1.
$$
In terms of $(\al;\be)$, the expression for
$\sum|f_{pq}(\lambda)|^2$ can be written as
$$
\frac{n^{p(p-1)+q(q-1)}}{x_n^2} \sum_{(\al;\be)} \prod_{1\le
i<j\le p} (\al_i-\al_j)^2 \prod_{1\le r<s\le q}
(\be_r-\be_s)^2\Big) \; \Big(1 + O\Big(\frac1n\Big)\Big).
$$
summed over all $(\al;\be)$ subject to the conditions stated
above, that is, over the set
$$
\Omega(p,q) \cap \Big({1 \over n-pq} \Z^{p+q}\Big),
$$
where $\Om(p,q)$ is the simplex introduced in \S9.6.

{}From the well-known formula
$$
{\Gamma(n+\const) \over \Gamma(n)} \sim n^{\const}, \qquad n\to\infty,
$$
we get
$$
x_n \sim n^{p^2+q^2-pq-{(p-q)^2+1 \over 2}} =
n^{p^2+q^2-1 \over 2},
$$
hence the coefficient in front of the sum has the asymptotics (as
$n\to\infty$)
$$
n^{p(p-1)+q(q-1)-(p^2+q^2-1)} = n^{-(p+q-1)}.
$$
But $p+q-1$ is exactly the dimension of  $\Omega(p,q)$, hence we
have
$$
\lim_{n\to\infty}\sum_{\la\in\Y(p,q)\cap\Y_n}
|f_{pq}(\lambda)|^2 =\int\limits_{\om=(\al;\be)\in\Omega(p,q)}
\prod_{1\le i<j\le p} (\alpha_i-\alpha_j)^2\; \prod_{1\le r<s\le
p} (\beta_r-\beta_s)^2\; d\omega,
$$
where $d\omega$ is  Lebesgue measure on $\Omega(p,q)$. Since the
integral in the right--hand side is finite, we conclude that the
Hardy type condition is satisfied. \qed\enddemo

\subhead 6.4. Proof of Theorems 6.2.1 and 6.2.2\endsubhead The definition of
the function $f_{pq}(\lambda)$ in Theorem 6.2.1 is identical with that given in
Lemma 6.3.3. By virtue of Lemmas 6.3.3 and 6.3.4, the function $f_{pq}$ belongs
to $\Cal F_{pq}$. Therefore, it determines a $K$--invariant vector $v_{pq}\in
H(T_{pq})$, which concludes the proof of the theorem.

Let us turn to Theorem 6.2.2. Compare two probability measures on $\Y_n$.

The first measure, denoted as $M^{(n)}_{pq}$, comes from the
coherent system corresponding to the spherical function
$\varphi_{pq}$ introduced in the statement of Theorem 6.2.2,
$$
M^{(n)}_{pq}(\la) = \frac{\Vert Q(\lambda)\,v_{pq}\Vert^2} {\Vert
v_{pq}\Vert^2}, \qquad \la\in\Y_n.
$$
By the general theory (see Theorem 9.7.3), the spectral measure $\si_{pq}$ is
the weak limit of the measures $\iota_n(M^{(n)}_{pq})$, where $\iota_n:
\Y_n\to\Om$ is the embedding defined just before Theorem 9.7.3.

The second measure, which we denote as $\overline M^{(n)}_{pq}$, has the form
$$
\overline M^{(n)}_{pq}(\la) =\frac{\Vert
Q_\la\,v^{(n)}_{pq}\Vert^2}{\Vert v^{(n)}_{pq}\Vert^2}
=\frac{|f_{pq}(\lambda)|^2}{\Vert v^{(n)}_{pq}\Vert^2} , \qquad
\la\in\Y_n,
$$
where $v^{(n)}_{pq}=P_nv_{pq}$ is the projection of $v_{pq}$
onto the subspace $H^n$, and the projection $Q_\la$ was defined
in \S6.1. We know that $\overline M^{(n)}_{pq}$ is concentrated
on the subset $\Y_n(p,q)=\Y(p,q)\cap\Y_n$.

Below we use the standard norm on signed Borel measures: given
such a measure $\mu$, its norm $\Vert \mu\Vert$ is defined as
the variance of $\mu_+$ plus the variance of $\mu_-$, where
$\mu=\mu_+-\mu_-$ stands for the Jordan decomposition of $\mu$.

\proclaim{Lemma 6.4.1} We have
$$
\Vert M^{(n)}_{pq}-\overline M^{(n)}_{pq}\Vert \to0, \quad n\to\infty.
$$
\endproclaim

\demo{Proof} Indeed, as is well known, for any two {\it
probability\/} Borel measures $\mu_1,\mu_2$, defined on one and
the same Borel space,
$$
\Vert\mu_1-\mu_2\Vert\le 2\sup_{X}|\mu_1(X)-\mu_2(X)|,
$$
where the supremum is taken over arbitrary Borel subsets. Let us apply this
inequality to $\mu_1=M^{(n)}_{pq}$, $\mu_2=\overline M^{(n)}_{pq}$. To simplify
the notation, set
$$
\xi_1=\Vert v_{pq}\Vert^{-1}\,v_{pq}, \qquad \xi_2=\Vert v^{(n)}_{pq}\Vert^{-1}
\,v_{pq}^{(n)}
$$
and, for any subset $X\subset\Y_n$, set
$$
Q_X=\sum_{\la\in X}Q_\la.
$$
The operators $Q_\la$, with $\la$ ranging over $\Y_n$, are
pairwise orthogonal projectors whose sum equal 1, whence $\Vert
Q_X\Vert\le1$ for any $X$. It follows
$$
\Vert M^{(n)}_{pq}-\overline M^{(n)}_{pq}\Vert \le
2\sup_{X\subset\Y_n} |(Q_X\xi_1,\xi_1)-(Q_X\xi_2,\xi_2))|\le
4\Vert\xi_1-\xi_2\Vert.
$$
But $\Vert\xi_1-\xi_2\Vert\to0$ as $n\to\infty$, because $\Vert
v_{pq}-v_{pq}^{(n)}\Vert\to0$. This completes the argument. \qed
\enddemo

On the other hand, it follows from the proof of Lemma 6.3.4 that the measures
$\bar\iota_n(\overline M^{(n)}_{pq})$ on $\Om(p,q)$ weakly converge to a
probability measure, which is absolutely continuous with respect to Lebesgue
measure on $\Om(p,q)$ and whose density is exactly as required in the statement
of the theorem.

\proclaim{Lemma 6.4.2} The measures $\iota_n(\overline M^{(n)}_{pq})$ on $\Om$
weakly converge, as $n\to\infty$, to the same limit measure concentrated on
$\Om(p,q)$.
\endproclaim

\demo{Proof} Indeed, let us compare the two embeddings,
$$
\bar\iota_n: \Y_n(p,q)\to\Om(p,q)\subset\Om \quad\text{and}\quad
\iota_n:\Y_n(p,q)\to\Om.
$$
For any $\la\in\Y_n(p,q)$, we can write (viewing $\bar\iota_n(\la)$ as an
element of $\Om$)
$$
\bar\iota_n(\la)=(\bar\al_1,\bar\al_2,\dots; \bar\be_1,\bar\be_2,\dots),
$$
where
$$
\bar\al_{p+1}=\bar\al_{p+2}=\dots=\bar\be_{q+1}=\bar\be_{q+2}=\dots=0,
$$
and similarly
$$
\iota_n(\la)=(\al_1,\al_2,\dots; \be_1,\be_2,\dots),
$$
where
$$
\al_{d+1}=\al_{d+2}=\dots=\be_{d+1}=\be_{d+2}\dots=0, \quad d=\max(p,q).
$$
Further, from the definition of $\bar\iota_n(\la)$ and  $\iota_n(\la)$ it
follows that for any fixed $i,j$,
$$
|\bar\al_i-\al_i|\le \frac{\const}n, \qquad |\bar\be_j-\be_j|\le \frac{\const}n
$$
where the constant does not depend on $\la$. Consequently, if $F$ is an
arbitrary bounded continuous function on $\R^\infty\times\R^\infty$ depending
on finitely many coordinates only, then the result of coupling between $f$ and
$\bar\iota_n(\overline M^{(n)}_{pq})-\iota_n(\overline M^{(n)}_{pq})$ is
$O(1/n)$. Since $\Om$ is a compact subset of $\R^\infty\times\R^\infty$, the
same is true for any continuous test function $F$ on $\Om$, which completes the
proof. \qed
\enddemo

By Lemma 6.4.1,
$$
\Vert \bar\iota_n(M^{(n)}_{pq}) -\bar\iota_n(\overline M^{(n)}_{pq})\Vert \to0,
\quad n\to\infty,
$$
so that the measures $\iota_n(M^{(n)}_{pq})$ must have the same limit as the
measures $\iota_n(\overline M^{(n)}_{pq})$. Together with Lemma 6.4.2 this
completes the proof of Theorem 6.2.2.

\head 7. The commutant of representation $T_{pq}$ \endhead

\subhead 7.1. The results \endsubhead Recall that if $z$ is an integer,
$z=k\in\Z$, then the representation $T_z=T_k$ splits into a direct sum of
subrepresentations called  the blocks,
$$
T_k\sim\bigoplus\limits_{p-q=k} T_{pq}
$$
where $p,q$ are assumed to be nonnegative integers, and $(p,q)\ne(0,0)$ if
$k=0$ (see \S5.7).

The main goal of this Section is to prove the following result.

\proclaim{Theorem 7.1.1} Let $z=k$ be an integer and let $T_{pq}$ be an
arbitrary block of the representation $T_k$. Then the vector $v_{pq}\in
H(T_{pq})$ constructed in \S6 is a cyclic vector of $T_{pq}$.
\endproclaim

Along with the Theorem 6.2.2 this implies our main result on representations
$T_k$ at integer points $z=k$:

\proclaim{Theorem 7.1.2} The block $T_{pq}\subset T_z$, where $z=k$ is an
integer and $p-q=k$, is equivalent to the direct integral of irreducible
spherical representations labelled by the points $\omega$ of the finite
dimensional face $\Omega(p,q)$ of the Thoma simplex, with respect to Lebesgue
measure $d\omega$.
\endproclaim

\proclaim{Corollary 7.1.3} The representations $T_{pq}$, where $p,q\in\Z$, are
pairwise disjoint. In particular, the representations $T_k$, $k\in\Z$, are
pairwise disjoint.
\endproclaim

\subhead 7.2. The commutant $\Cal A_{pq}$\endsubhead Recall that the
representation $T_{pq}$ is an inductive limit of finite dimensional
representations of the groups $G(n)$ in the spaces $H^n\cap H(T_{pq})\subset
\Reg^n$. Since
$$
(v_{pq},\xi_\lambda) = f_{pq}(\lambda) \ne 0 \qquad \text{ for }
\lambda \in \Bbb Y(p,q),
$$
the projection $v^{(n)}_{pq}$ of the vector $v_{pq}$ onto
subspace $H^n\cap H(T_{pq})$ is a cyclic vector, for every $n$.
If one of the numbers $p$, $q$ vanishes, we have
$v^{(n)}_{pq}\equiv v_{pq}$ and hence $v_{pq}$ is obviously a
cyclic vector. But if $v^{(n)}_{pq}\ne v_{pq}$, the fact that
each vector $v^{(n)}_{pq}$ is cyclic in the corresponding
representation of the subgroup $G(n)$ does not formally imply
that the limiting vector is cyclic, too.  See subsection 7.8.
for a counterexample.

In order to prove that the vector $v_{pq}$ is cyclic we study
the commutant of $T_{pq}$, making use of Proposition 5.3.1 and
Theorem 5.5.1. We shall show that the matrix element
$((\,\cdot\,)v_{pq},v_{pq})$ provides a {\it faithful\/} state
on the commutant, which implies the cyclicity immediately.
Unfortunately, the proof turns out to be rather long.

Fix a block $T_{pq}\subset T_k$, where $k=p-q$, and assume that $p\ge1$,
$q\ge1$ (if one of the numbers $p$, $q$ vanishes, the claim of the Theorem
7.1.1 is trivial).

Set $\Bbb Y_n(p,q)=\Bbb Y_n\cap\Bbb Y(p,q)$; this is the $n$th level of the
graph $\Bbb Y(p,q)$. It is not empty, starting with $n=pq>0$. The graph $\Bbb
Y(p,q)$ can be identified with the direct product of the graphs $\Bbb Y(p,0)$
and $\Bbb Y(0,q)$.

Denote by $\Cal A_{pq}$ the commutant of the representation $T_{pq}$. According
to Proposition 5.3.1, there is a linear isomorphism $\Cal A_{pq}\to\widetilde
{\Cal A}_{pq}$, where $\widetilde {\Cal A}_{pq}$ is the Banach space of
complex--valued bounded functions $A(\lambda)$ on $\Bbb Y(p,q)$, satisfying the
condition
$$
A(\lambda) = \sum_{\nu\searrow\lambda} \wt p_z(\lambda,\nu)\;
A(\nu), \qquad \lambda \in \Bbb Y(p,q),
$$
with the norm
$$
\Vert A\Vert = \sup_\lambda |A(\lambda)|.
$$
According to Theorem 5.5.1, the function $\wt p_z(\lambda,\nu)$
has the form
$$
\wt p_z(\lambda,\nu) =p_z(\lambda,\nu)=
\frac{|z+c_{\lambda,\nu}|^2}{|z|^2+n}\cdot
\frac{\dim\nu}{(n+1)\dim\lambda}.
$$
It will be important for us that $p_z(\lambda,\nu)\ne0$ for all
couples $\lambda\nearrow\nu$ in the graph $\Bbb Y(p,q)$.

Denote by $\Cal A^+_{pq}$ the cone of nonnegative Hermitian
operators in $\Cal A_{pq}$. Its image with respect to the
isomorphism $\Cal A_{pq}\to\widetilde {\Cal A}_{pq}$ is
contained in the cone $\widetilde {\Cal
A}^+_{pq}\subset\widetilde {\Cal A}_{pq}$ generated by the
functions $A\in\widetilde {\Cal A}_{pq}$ with nonnegative
values. One can  show that the image actually coincides with
$\widetilde{\Cal A}^+_{pq}$, but we do not need this fact.

\subhead 7.3. A faithful state on the algebra $\Cal A_{pq}$
\endsubhead
Denote by $\lambda_{\min}$ the rectangular diagram of size
$p\times q$; this is the only vertex of the graph $\Bbb Y(p,q)$
at the level $pq$.

\proclaim{Lemma 7.3.1} The linear functional
$$
\eta: A \mapsto A(\lambda_{\min}), \qquad A \in \widetilde{\Cal A}_{pq}
$$
provides a faithful trace on the algebra $\Cal A_{pq}$.
\endproclaim

\demo{Proof} Clearly, the functional $\eta$ equals $1$ at the
function $A(\lambda)\equiv1$ (which determines the identity of
$\Cal A_{pq}$), has norm $1$ and is nonnegative on the cone
$\Cal A^+_{pq}$. Hence, it defines a state on the algebra $\Cal
A_{pq}$. It remains to prove that this state is faithful. To do
this we shall show that $\eta(A)>0$ for every non--zero function
$A(\lambda)$ in $\widetilde{\Cal A}^+_{pq}$.

Define a ``weight function'' $\frak M_{pq}(\lambda)$ on the graph $\Bbb Y(p,q)$
by the recurrence relation
$$
\frak M_{pq}(\nu) = \sum_{\lambda:\lambda\nearrow\nu} \frak
M_{pq}(\lambda)\,  p_z(\lambda,\nu), \qquad |\nu| > pq
$$
and the initial condition
$$
\frak M_{pq}(\lambda_{\min}) = 1.
$$

Denote by $\frak M^{(n)}_{pq}$ the measure on $\Y_n(p,q)$, $n\ge pq$, with the
weights $\frak M_{pq}(\lambda)$ at the vertices $\lambda\in\Bbb Y_n(p,q)$.
Since $p_z(\lambda,\nu)$ has the property
$$
\sum_{\nu\searrow\lambda} p_z(\lambda,\nu) = 1,
$$
all $\frak M^{(n)}_{pq}$ are probability measures.

By virtue of the main relation for the functions
$A(\,\cdot\,)\in \wt{\Cal A}_{pq}$, the expression
$$
\eta_n(A) = \langle \frak M^{(n)}_{pq},A \rangle = \sum_{\lambda\in\Bbb
Y_n(p,q)} \frak M_{pq}(\lambda)\, A(\lambda), \qquad n \ge pq,
$$
does not depend on $n$. Hence, it coincides with
$\eta_{pq}(A)=\eta(A)$.

Since $p_z(\lambda,\nu)>0$, all the weights $\frak M_{pq}(\lambda)$ are
strictly positive. Let now $A(\,\cdot\,)$ be a non--zero function from
$\widetilde {\Cal A}^+_{pq}$. Then there exists $\lambda\in\Bbb Y(p,q)$, such
that $A(\lambda)>0$. If $n=|\lambda|$, then $\langle \frak M^{(n)}_{pq},A
\rangle>0$, hence we conclude that $\eta(A)=\eta_n(A)>0$. \qed\enddemo

\subhead 7.4. Evaluation of $\frak M_{pq}$ \endsubhead As in \S6
above, we associate with a diagram $\lambda\in\Bbb Y(p,q)$ a
couple of diagrams $(\lambda^+$, $\lambda^-)$, and in order to
simplify the notation we set $a_i=\lambda^+_i$,
$b_r=\lambda^-_r$ for $1\le i\le p$, $1\le r\le q$.

\proclaim{Lemma 7.4.1} The ``weight function'' $\frak M_{pq}(\lambda)$
introduced in the proof of Lemma 7.2.1 can be given by the explicit formula
$$
\gathered
\frak M_{pq}(\lambda) = \frac{C(p,q)}{(n-pq+1)_{p^2+q^2-pq-1}} \times
\\ \times \frac{\underset{1\le i<j\le p}\to{\prod} (a_i-a_j+j-i)^2
\underset{1\le r<s\le q}\to{\prod} (b_r-b_s+s-r)^2} {\underset{i,r}\to{\prod}
(a_i+b_r+p-i+q-r+1)},
\endgathered
$$
where $n=|\lambda|$ and
$$
C(p,q) =
\frac{(p^2+q^2-pq-1)! \prod_{i,r} (p-i+q-r+1)}
{\big(\prod_i (p-i)!\; \prod_r (q-r)! \big)^2}.
$$
\endproclaim

\demo{Proof} Recall that $\frak M_{pq}(\lambda)$ satisfies the recurrence
relation
$$
\frak M_{pq}(\nu) = \sum_{\lambda\nearrow\nu}
\frac{(p-q+c_{\lambda\nu})^2}{(p-q)^2+n}\, \frac{\dim\nu}{(n+1)\dim\lambda}\,
\frak M_{pq}(\lambda),
$$
where $|\nu|=n+1>pq$.

Set in this formula
$$
\frak M_{pq}(\lambda) = \frac{\underset{b\in\lambda\setminus\lambda_{\min}}
\to{\prod}(p-q+c(b))^2} {((p-q)^2 + n-1)!}\; \frac{\dim\lambda}{n!}\; \frak
M'_{pq}(\lambda)
$$
with a new unknown function $\frak M_{pq}'(\lambda)$,
$n=|\lambda|$. Then for $\frak M'_{pq}(\lambda)$ we get the
recurrence relation
$$
\frak M'_{pq}(\nu) = \sum_{\lambda\nearrow\nu} \frak M'_{pq}(\lambda).
$$
Its solution has the form
$$
\frak M'_{pq}(\lambda) = \const\cdot \operatorname{Dim}\lambda,
$$
where $\operatorname{Dim}\la$ stands for the number of paths, in the graph
$\Y(p,q)$, going from $\la_{\min}$ to $\la$ (the dimension function of the
graph $\Bbb Y(p,q)$).

Since the graph $\Bbb Y(p,q)$ is isomorphic to the product of
the graphs $\Y(p,0)\times\Y(q,0)$ (with the empty diagrams
$\varnothing$ added to $\Y(p,0)$ and to $\Y(0,q)$), we get
$$
\operatorname{Dim}\lambda = \frac{(n-pq)!}{|\lambda^+|!\,|\lambda^-|!}\;
\dim\lambda^+\; \dim\lambda^-.
$$
It follows that
$$
\frak M_{pq}(\lambda) = \const
\frac{\underset{b\in\lambda\setminus\lambda_{\min}} \to{\prod}(p-q+c(b))^2}
{((p-q)^2 + n-1)!}\; \frac{\dim\lambda}{|\lambda|!}\;
\frac{\dim\lambda^+}{|\lambda^+|!}\; \frac{\dim\lambda^-}{|\lambda^-|!}\;
(n-pq)!.
$$
Substitute here the explicit expressions:
$$
\underset{b\in\lambda/\lambda_{\min}}\to{\prod}(p-q+c(b))^2 =
\left(\prod_{i=1}^p \frac{(a_i+p-i)!}{(p-i)!} \prod_{r=1}^q
\frac{(b_r+q-r)!}{(q-r)!}\right)^2
$$
(this is readily checked),
$$
\frac{\dim\lambda}{|\lambda|!} =
\frac{\underset{1\le i<j\le p}\to{\prod} (a_i-a_j+j-i)\,
\underset{1\le r<s\le q}\to{\prod} (b_r-b_s+s-r)}
{\underset{1\le i\le p}\to{\prod} (a_i+p-i)!
\underset{1\le r\le q}\to{\prod} (b_r+q-r)!\;
\underset{i,r}\to{\prod}(a_i+b_r+p-i+q-r+1)},
$$
(this can be derived from the hook formula in the same way as
formula \tht{2.7} in Olshanski \cite{Ol5}), and
$$
\frac{\dim\lambda^+}{|\lambda^+|!} = \frac{\underset{1\le i<j\le
p}\to{\prod} (a_i-a_j+j-i)} {\underset{1\le i\le p}\to{\prod}
(a_i+p-i)!},\qquad \frac{\dim\lambda^-}{|\lambda^-|!} =
\frac{\underset{1\le r<s\le q}\to{\prod} (b_r-b_s+s-r)}
{\underset{1\le r\le q}\to{\prod} (b_r+q-r)!}
$$
(this is the well--known Frobenius formula, see, e.g., Macdonald
\cite{Mac, Ex. I.7.6}). Then we get the desired expression,
where the constant $\const$ can be found from the condition
$\frak M_{pq}(\lambda_{\min})=1$. \qed\enddemo

\subhead 7.5. The limit of the measures $\frak M^{(n)}_{pq}$
\endsubhead
We shall need two technical lemmas.

\proclaim{Lemma 7.5.1} For any $\varepsilon>0$ there exists
$\delta>0$ such that
$$
\underset{\lambda\in\Bbb Y_n(p,q);\; a_p+b_q\le\delta n}\to {\sum} \frak
M_{pq}(\lambda) \le \varepsilon,
$$
for all $n$ large enough.
\endproclaim

Proof will be given in \S7.7.

\proclaim{Lemma 7.5.2} The function
$$
\frac{\underset{1\le i<j\le p}\to{\prod}(\alpha_i-\alpha_j)^2
\underset{1\le r<s\le q}\to{\prod}(\beta_r-\beta_s)^2}
{\underset{i,r}\to{\prod}(\alpha_i+\beta_r)}
$$
on the simplex $\Omega(p,q)$ is integrable with respect to
Lebesgue measure $d\omega$.
\endproclaim

Proof will be given in \S7.8.

Consider now the embedding $\bar\iota_n:
\Y_n(p,q)\to\Omega(p,q)$ introduced in \S6.3. Let $\wh{\frak
M}^{(n)}_{pq}=\bar\iota_n(\frak M^{(n)}_{pq})$ be the
push--forward of the probability measure $\frak M^{(n)}_{pq}$;
this is a probability measure on $\Om(p,q)$.

\proclaim{Lemma 7.5.3} As $n\to\infty$, the measures $\wh{\frak
M}^{(n)}_{pq}$ on $\Om(p,q)$ weakly converge to a certain
probability measure $\wh{\frak M}^{(\infty)}_{pq}$. The measure
$\wh{\frak M}^{(\infty)}_{pq}$ is absolutely continuous with
respect to Lebesgue measure $d\omega$ on $\Om(p,q)$, and has the
density
$$
C(p,q)\; \frac{\underset{1\le i<j\le p}\to{\prod}(\alpha_i-\alpha_j)^2
\underset{1\le r<s\le q}\to{\prod}(\beta_r-\beta_s)^2}
{\underset{i,r}\to{\prod}(\alpha_i+\beta_r)},
$$
where the constant is the same as in Lemma 7.4.1.
\endproclaim

\demo{Proof} Let $\lambda\in\Y_n(p,q)$ and
$\alpha;\beta)=\bar\iota_n(\la)$. It follows from the expression
for $\frak M_{pq}(\lambda)$ in Lemma 7.4.1 that
$$
\gathered
\frak M^{(n)}_{pq}(\lambda) = C(p,q)\; n^{-(p+q-1)}\; (1+O(1/n)) \times \\
\times \frac{\underset{1\le i<j\le p}\to{\prod}
(\alpha_i-\alpha_j+O(1/n))^2 \underset{1\le r<s\le q}\to{\prod}
(\beta_r-\beta_s+O(1/n))^2}
{\underset{i,r}\to{\prod}(\alpha_i+\beta_r+O(1/n))}\,.
\endgathered
$$
Given $\de>0$, we split the simplex $\Omega(p,q)$ into the union
of two subsets,
$$
\Omega(p,q) = \Omega_{\ge\delta}(p,q) \cup \Omega_{<\delta}(p,q),
$$
determined by the conditions $\alpha_p+\beta_q\ge\delta$ and
$\alpha_p+\beta_q<\delta$, respectively. On the set
$\Omega_{\ge\delta}(p,q)$ the denominator of the expression for
$\frak M^{(n)}_{pq}(\lambda)$ is bounded from below, uniformly
in $n$. This implies the weak convergence of measures
$$
\wh{\frak M}^{(n)}_{pq}\bigg|_{\Omega_{\ge\delta}(p,q)} \longrightarrow
\wh{\frak M}^{(\infty)}_{pq}\bigg|_{\Omega_{\ge\delta}(p,q)}.
$$

By Lemma 7.5.1, the total mass of the set
$\Omega_{<\delta}(p,q)$ with respect to $\wh{\frak
M}^{(n)}_{pq}$ can be made arbitrarily small (uniformly in $n$),
provided that $\delta$ is chosen sufficiently small, and $n$ is
large enough. By Lemma 7.5.2, the mass of the set
$\Omega_{<\delta}(p,q)$ with respect to $\wh{\frak
M}^{(\infty)}_{pq}$ also tends to $0$, together with $\delta$.
This implies that $\wh{\frak M}^{(n)}_{pq}$ weakly converges to
$\wh{\frak M}^{(\infty)}_{pq}$ on the entire simplex
$\Omega(p,q)$. \qed\enddemo

\proclaim{Lemma 7.5.4} The state
$$
A \mapsto {(Av_{pq},v_{pq}) \over (v_{pq},v_{pq})}
$$
on the algebra $\Cal A_{pq}$ is faithful.
\endproclaim

\demo{Proof} We have to prove that if $A\in \Cal A_{pq}^+$,
$A\ne0$, then $(Av_{pq},v_{pq})>0$. Set
$$
A^{(n)}=P_nAP_n, \qquad v_{pq}^{(n)} = P_n v_{pq}.
$$
If $n\to\infty$, the operator $P_nAP_n$ converges to $A$
strongly, hence
$$
(Av_{pq},v_{pq}) = \lim_{n\to\infty} (P_nAP_n v_{pq},v_{pq}) =
\lim_{n\to\infty} (A^{(n)} v^{(n)}_{pq},v^{(n)}_{pq}).
$$

According to Proposition 5.3.1,
$$
A^{(n)} = \underset{\lambda\in\Bbb Y_n(p,q)}\to{\sum}
A(\lambda)P(\la),
$$
where $A(\lambda)$ is the function in $\widetilde{\Cal A}_{pq}^+$ corresponding
to the operator $A$. On the other hand,
$$
v^{(n)}_{pq} = \underset{\lambda\in\Bbb Y_n(p,q)}\to{\sum}
f_{pq}(\lambda)\, \xi_\lambda
$$
by definition of the vector $v_{pq}$. Therefore,
$$
(A^{(n)} v^{(n)}_{pq},v^{(n)}_{pq}) = \underset{\lambda\in\Bbb
Y_n(p,q)}\to{\sum} |f_{pq}(\lambda)|^2 A(\lambda) = \underset{\lambda\in\Bbb
Y_n(p,q)}\to{\sum} A(\lambda)\,\frak M_{pq}(\lambda) {|f_{pq}(\lambda)|^2 \over
\frak M_{pq}(\lambda)}.
$$

We have explicit expressions for both $f_{pq}(\lambda)$ and $\frak
M_{pq}(\lambda)$, see Theorem 6.2.1 and Lemma 7.4.1. These expressions imply
that
$$
{|f_{pq}(\lambda)|^2 \over \frak M_{pq}(\lambda)} = c_{pqn} \prod_{i,r}
(a_i+b_r+p-i+q-r+1),
$$
where
$$
c_{pqn} = \frac{(n-pq+1)_{p^2+q^2-pq-1}}
{x_n^2\,\cdot\,C(p,q)}
$$
(we have employed here an equivalent expression for
$f_{pq}(\lambda)$, see Lemma 6.3.3). Since
$$
x_n \sim n^{(p^2+q^2-1)/2},
$$
we have
$$
c_{pqn} \sim {1 \over C(p,q)\; n^{pq}},
$$
whence
$$
\multline {|f_{pq}(\lambda)|^2 \over \frak M_{pq}(\lambda)} = \big(\const +
o(1)\big)
\prod_{i,r} {a_i+b_r+p-i+q-r+1 \over n}\\
= \const\; \prod_{i,r}(\alpha_i+\beta_r) + o(1),
\endmultline
$$
where we assume that $(\alpha;\beta)=\bar\iota_n(\lambda)$.

Let $\frak N^{(n)}_{pq}$ denote the measure on $\Y_n(p,q)$
defined by
$$
\frak N^{(n)}_{pq}(\la)=A(\lambda)\frak M_{pq}(\lambda), \qquad
\la\in\Y_n(p,q).
$$
In the notation introduced in the proof of Lemma 7.3.1,
$$
\sum_{\la\in\Y_n(p,q)}A(\la)\frak M_{pq}(\la)=\eta_n(A).
$$
We have noted there that $\eta_n(A)$ does not depend on $n$, and
is strictly positive. Therefore, the mass of the set $\Y_n(p,q)$
with respect to the measure $\frak N^{(n)}_{pq}$ is strictly
positive and does not depend on $n$. Consequently, passing to an
appropriate subsequence of indices $n_1<n_2<\dots$, we may
conclude that the measures $\wh{\frak
N}^{(n)}_{pq}:=\bar\iota_n(\frak N^{(n)}_{pq})$ converge weakly
on $\Omega(p,q)$ to a certain nonzero measure $\wh{\frak
N}^{(\infty)}_{pq}$.

On the other hand, for the weights of the measures $\frak N^{(n)}_{pq}$ and
$\frak M^{(n)}_{pq}$ there is an estimate
$$
\frak N^{(n)}_{pq}(\lambda) = A(\lambda)\,\frak M^{(n)}_{pq}(\lambda) \le \Vert
A\Vert\, \frak M^{(n)}_{pq}(\lambda),
$$
hence
$$
\wh{\frak N}^{(n)}_{pq} \le \Vert A\Vert\, \wh{\frak M}^{(n)}_{pq}.
$$

By Lemma 7.5.3, the measures $\wh{\frak M}^{(n)}_{pq}$ converge
weakly to a measure $\wh{\frak M}^{(\infty)}_{pq}$ which is
absolutely continuous with respect to the Lebesgue measure
$d\omega$ on the simplex $\Omega(p,q)$. It follows that
$$
\wh{\frak N}^{(\infty)}_{pq} \le \Vert A\Vert\, \wh{\frak
M}^{(\infty)}_{pq}\,,
$$
which implies that $\wh{\frak N}^{(\infty)}_{pq}$ is absolutely
continuous, too.

Return now to the quantity $(Av_{pq},v_{pq})$. We have shown
that it can be represented as a limit
$$
(Av_{pq},v_{pq}) = \lim_{n\to\infty} (A^{(n)}v^{(n)}_{pq},v^{(n)}_{pq}) =
\lim_{n\to\infty}\, \const\, \bigg\langle \wh{\frak N}^{(n)}_{pq},\big(
\prod_{i,r}(\alpha_i+\beta_r) + o(1)\big)\bigg\rangle,
$$
where $\const>0$. Hence,
$$
(Av_{pq},v_{pq}) = \const\, \big\langle \wh{\frak N}^{(\infty)}_{pq},
\prod_{i,r}(\alpha_i+\beta_r)\big\rangle.
$$
Since the measure $\wh{\frak N}^{(\infty)}_{pq}$ is absolutely continuous and
nonzero, and the function $\prod(\alpha_i+\beta_r)$ on the simplex
$\Omega(p,q)$ is continuous and positive almost everywhere, the result is a
strictly positive number. \qed\enddemo

\subhead 7.6. The Cauchy determinant \endsubhead We shall need a generalization
of the classical Cauchy identity
$$
\det \left[{1 \over x_i+y_j} \right]_{i,j=1}^m =
{V(x)\, V(y) \over \prod_{i,j} (x_i+y_j)},
$$
where
$$
V(x)=\prod_{1\le i<j\le m} (x_i - x_j),
$$
to the case when  the numbers of $x_i$'s and $y_j$'s are not
necessarily equal.

Set
$$
x=(x_1,\,\ldots,x_p), \qquad
y=(y_1,\,\ldots,y_q),
$$
and assume (to be definite) that $p\le q$. We shall denote by
the symbol $y=y'\sqcup y''$ an arbitrary decomposition of the
set of variables $y=(y_1,\,\ldots,y_q)$ into disjoint subsets of
cardinalities $q-p$ and $p$ respectively:
$$
y' = (y'_1,\,\ldots,y'_{q-p}) = (y_{i_1},\,\ldots,y_{i_{q-p}}),
\qquad
y'' = (y''_1,\,\ldots,y''_p) = (y_{j_1},\,\ldots,y_{j_p}),
$$
where
$$
i_1<\ldots<i_{q-p}, \quad j_1<\ldots<j_p, \qquad
\{i_1,\,\ldots,i_{q-p}\}\cup
\{j_1,\,\ldots,j_p\} = \{1,\,\ldots,q\}.
$$
Set $\sgn(y',y'')=\pm1$, where the sign plus/minus is taken
depending on parity/imparity  of the number of inversions in the
permutation $(i_1,\,\ldots,i_{q-p},j_1,\,\ldots,j_p)$ of the
numbers $1$, $\ldots$, $q$.

\proclaim{Lemma 7.6.5} The following formula generalizes the Cauchy
determinant:
$$
\multline \frac{V(x_1,\,\ldots,x_p)\, V(y_1,\,\ldots,y_q)} {\prod_{i,j} (x_i +
y_j)}\\
= \underset{y'\sqcup y''=y}\to{\sum} \sgn(y',y'') V(y'_1,\,\ldots,y'_{q-p})
\det\left[\frac1{x_i+y''_j} \right]_{i,j=1}^p.
\endmultline
$$
\endproclaim

\demo{Proof} One can easily derive this identity from the
classical Cauchy identity by replacing the variables
$x_1,\,\ldots,x_p$ with $x''_1,\,\ldots,x''_p$, adding a new
group of variables $x'_1,\,\ldots,x'_{q-p}$, writing down the
Cauchy identity in terms of the variables $(\varepsilon x'_1$,
$\ldots$, $\varepsilon x'_{q-p}$; $x''_1$, $\ldots,x''_p)$ and
$y_1,\,\ldots,y_q$, and then applying the Laplace rule while
$\varepsilon$ goes to $0$. \qed\enddemo

\demo{Another proof} The argument given below is similar to the
well--known proof of the classical Cauchy identity.

Denote the right--hand side of the identity in question by
$A(x,y)$. Then the product $A(x,y)\Pi(x_i+y_j)$ is a polynomial.
It suffices to check the following three claims:

(i) $A(x,y)$ is skew symmetric, separately in $x$ and in $y$.

(ii) $\deg A(x,y)=\deg V(x)+\deg V(y)-\deg\Pi(x_i+y_j)$.

(ii) Let us order the variables as $x_1>\dots>x_p>y_1>\dots>y_q$
and consider the corresponding order on the monomials; then the
leading term  in the expansion of $A(x,y)\Pi(x_i+y_j)$ is the
same as that for $V(x)\,V(y)$.

Now we have:

(i) The skew symmetry with respect to $x$ is evident. Let us
show that $A(x,y)$ changes the sign upon the elementary
transposition $y_j\leftrightarrow y_{j+1}$, where
$j=1,\,\ldots,q-1$. If the variables $y_j,y_{j+1}$ enter the
same group, $y'$ or $y''$, then the corresponding term is
already skew symmetric. Consider now the terms for which $y_j$
and $y_{j+1}$ belong to different groups. Those terms split into
pairs: in each pair the terms are switched by the transposition
$y_j\leftrightarrow y_{j+1}$ and the corresponding signs are
opposite. Thus, the skew symmetry follows.

(ii) This is trivial.

(iii) It is readily verified that the leading monomial comes
from the partition $y'=(y_1,\,\ldots,y_{q-p})$,
$y''=(y_{q-p+1},\,\ldots,y_q)$. \qed\enddemo

\subhead 7.7. Proof of Lemma 7.5.1 \endsubhead Without loss of
generality we may assume that $p\le q$. Set
$$
x_i = a_i+p-i, \quad 1 \le i \le p; \qquad
y_j = b_j+q-j, \quad 1 \le j \le q.
$$
Taking into account the explicit formula for $\frak M_{pq}(\lambda)$ (Lemma
7.4.1) we have to prove the estimate
$$
{1 \over n^{p^2-pq+q^2-1}} \sum_{x,y}
\frac{V^2(x)\, V^2(y)}{\prod_{i,j} (x_i+y_j+1)} = O(\delta),
$$
summed over the integer vectors $x\in\Z_+^p$, $y\in\Bbb Z_+^q$ satisfying the
conditions
$$
\gathered
x_1>\ldots>x_p\ge0, \qquad y_1>\ldots>y_q\ge0, \qquad
x_p+y_q\le\delta n \\
x_1 + \ldots + x_p + y_1 + \ldots + y_q = n + \const,
\endgathered
\tag 7-6-1
$$
where ``$\const$'' is an integer not depending on $n$.

Apply the identity of Lemma 7.6.1. Since all variables are of
order $O(n)$, we have
$$
V(x)\,V(y)\,V(y') = O\left(
n^{{1\over2}(p(p-1)+q(q-1)+(q-p)(q-p-1))}\right) =
O\left(n^{p^2-pq+q^2-q}\right).
$$
Using this and expanding the determinant in the right hand side
of the identity in question, we reduce the problem to the
following bound:
$$
\frac1{n^{q-1}} \sum_{x,y}
\frac1{(x_1+y''_{\sigma(1)}+1)\ldots(x_p+y''_{\sigma(p)}+1)} =
O(\delta),
$$
where a splitting of variables $y'\sqcup y''=y$ is fixed,
$\sigma$ is a fixed permutation of the indices $1,\,\ldots,p$,
and summation is again taken under the conditions \tht{7-6-1}.

There are three possibilities:

(i) $y''$ contains $y_q$, and $y''_{\sigma(p)}=y_q$;

(ii)  $y''$ contains $y_q$, but $y''_{\sigma(p)}\ne y_q$;

(iii) $y''$ does not contain $y_q$.

The case (iii) reduces to that of (i), since replacing
$y''_{\sigma(p)}$ with $y_q<y''_{\sigma(p)}$ only increases the
sum.

The case (ii) can also be reduced to (i). Indeed, there exists
$i<p$ such that $y''_{\sigma(i)}=y_q$. We can now change
$\sigma$ by switching $y''_{\sigma(i)}$ and $y''_{\sigma(p)}$.
This can be done by virtue of the inequality
$$
\frac1{(A_1+B_2)(A_2+B_1)} < \frac1{(A_1+B_1)(A_2+B_2)},
$$
(correct for $A_1>A_2>0$, $B_1>B_2>0$) which we apply to
$$
A_1 = x_i +{1\over2}, \quad A_2=x_p+{1\over2}, \qquad
B_1 = y''_{\sigma(p)}+{1\over2}, \quad
B_2=y''_{\sigma(i)}+{1\over2}.
$$

Hence, we are left with the case (i), where
$y''_{\sigma(p)}=y_q$. Let us relax the system of restrictions
\tht{7-6-1} by removing from it the inequalities
$x_1>\ldots>x_p$, $y_1>\ldots>y_q$ (but we still assume that all
variables are nonnegative integers). This will result in a
bigger amount of arrays $x,y$, hence the sum will increase, too.

Set
$$
w_1=x_1+y''_{\sigma(1)},\,\quad\dots,\quad
w_{p-1}=x_{p-1}+y''_{\sigma(p-1)}\,, \quad
w_p=x_p+y''_{\sigma(p)}=x_p+y_p.
$$
Note that for every $w_i$ there are precisely $w_i+1$ ways to
represent it as a sum of two nonnegative terms. Therefore, our
bound reduces to the following one: the number of vectors
$(y'_1,\dots, y'_{q-p},w_1,\dots,w_p)\in\Z_+^q$ such that
$$
y'_1+\dots+y'_{q-p}+w_1+\dots+w_p = n+\const, \qquad w_p \le
\delta n
$$
is of order $O(\de)n^{q-1}$. This bound is readily verified.

\subhead 7.8. Proof of Lemma 7.5.2 \endsubhead We shall prove
the analogous claim for the same density on the larger simplex
$\widetilde{\Omega}(p,q)$, obtained by removing the restrictions
$\alpha_1\ge\ldots\ge\alpha_p$, $\beta_1\ge\ldots\ge\beta_q$. In
other words, $\widetilde{\Omega}(p,q)$ is just the standard
$(p+q-1)$--dimensional simplex. Without loss of generality we
may assume that $p\le q$. By the Cauchy type identity of Lemma
7.6.5 and taking into account the trivial bounds $\alpha_i\le1$,
$\beta_j\le1$ we are reduced to the following claim:
$$
\int \frac{dw}{(\alpha_1+\beta_1)\ldots(\alpha_p+\beta_p)} <
\infty,
$$
where $dw$ is the Lebesgue measure on the standard $(p+q-1)$--dimensional
simplex $\sum\alpha_i+\sum\beta_j=1$, $\alpha_i\ge0$, $\beta_j\ge0$, the domain
of integration. Project our simplex onto a $(q-1)$--dimensional simplex by
applying the map
$$
(\alpha_1,\,\ldots,\alpha_p;\beta_1,\,\ldots,\beta_q)
\longmapsto
(\alpha_1+\beta_1,\,\ldots,\alpha_p+\beta_p;
\beta_{p+1},\,\ldots,\beta_q).
$$
Under this projection, the push--forward of the measure with
density
$$
{1 \over (\alpha_1+\beta_1)\ldots(\alpha_p+\beta_p)}
$$
is simply the Lebesgue measure (this is a particular case of a
more general fact on the behavior of Dirichlet measures under
projections of simplices), which completes the proof of the
lemma.

This  completes the proofs of Theorems 7.1.1, 7.1.2.

\subhead 7.9. Example of a noncyclic vector with cyclic
projections
\endsubhead One can think that the proof of the cyclicity of the vectors $v_{pq}$
given above is a little bit too involved. Surely, one would try
to simplify it. But the claim is not entirely trivial, and the
goal of the present addendum is show that by an example.

We shall provide a representation $T=\varinjlim T_n$ of a group
$G=\varinjlim G_n$, and a vector $\xi\in H(T)$, such that, for
each $n$, the projection $\xi_n$ of the vector $\xi$ onto the
subspace $H(T_n)$ is cyclic for $T_n$, though the vector $\xi$
itself is not.

Set $G_n=\Z_2^n$, so that $G=\bigoplus_1^\infty\Z_2$. The dual group to $G_n$
is again $\Z_2^n$, while the dual group to $G$ is a compact group,
$\widehat{G}=\prod_1^\infty\Z_2$. With arbitrary finite Borel measure $\sigma$
on $\widehat{G}$ one can associate a unitary representation $T$ of the group
$G$ acting in the Hilbert space $H=L^2(\widehat{G},\sigma)$ by the formula
$$
(T(g)f)(\chi) = \chi(g)\,f(\chi), \qquad
g\in G,\quad \chi\in \widehat{G},\quad f\in H.
$$
Denote by $H_n\subset H$ the subspace of functions $f(\chi)$
depending on $\chi_n:=\chi\big|_{G_n}$ only, and by $T_n$ the
natural representation of the group $G_n$ in $H_n$. The
representation $T$ coincides with $\varinjlim T_n$. The space
$H_n$ can be identified with $L^2(\Z_2^n,\sigma^{(n)})$, where
$\sigma^{(n)}$ is the image of the measure $\sigma$ under the
projection $\widehat{G}\to\Z_2^n$ (taking of the first $n$
coordinates). If $\xi=f(\,\cdot\,)\in H$, then the vector
$\xi_n\in H_n=L^2(\Bbb Z_2^n,\sigma^{(n)})$ can be obtained by
integrating $f$ along the fibers of the projection $G\to\Z_2^n$.

For $p\in(0,1)$ denote by $\sigma_p$ the Bernoulli measure on
$\widehat{G}=\prod_1^\infty \Z_2=\prod_1^\infty\{0,1\}$ with the
weights $p$ and $1-p$ at the points $0$ and $1$ respectively. By
virtue of the law of large numbers, $\sigma_p$ is supported by
the set $X_p\subset\prod_1^\infty\{0,1\}$ formed by 0--1
sequences with the limiting frequency of 0's equal to $p$.

Take two distinct numbers $p,p'\in(0,1)$ and set
$\sigma=\sigma_p+\sigma_{p'}$. The measure $\sigma$ is supported
by the union of two disjoint sets $X_p$, $X_{p'}$ each of which
has measure $1$. Take for $\xi$ the characteristic function of
the set $X_p$. Clearly, $\xi$ is not cyclic, since its cyclic
span is a proper subspace generated by the functions in
$L^2(\widehat{G},\sigma)$ supported by $X_p$.

On the other hand, the vector $\xi_n$, as a function on
$\Z_2^n$, coincides with the Radon--Nikodym derivative
$\sigma_p^{(n)}/(\sigma_p^{(n)}+\sigma_{p'}^{(n)})$, hence is a
strictly positive function. Therefore, $\xi_n$ is a cyclic
vector for $T_n$ for any $n$.

\head 8. Disjointness \endhead

\subhead 8.1. The results \endsubhead Our aim is to prove the following result

\proclaim{Theorem 8.1.1} If $z$ ranges over the upper half--plane $\Im z\ge0$,
then the representations $T_z$ are pairwise disjoint.
\endproclaim

The assumption $\Im z\ge0$ is introduced because $T_z\sim T_{\bar z}$. The
definition of disjoint representations is given in \S9...

Let $\sigma_z$ the spectral measure of the character $\chi_z$ (\S4.1). We shall
deduce Theorem 8.1.1 from the following result.

\proclaim{Theorem 8.1.2} Assume $z$ ranges over the set $\{z\in\C\setminus\Z,
\Im z\ge0\}$.

{\rm(i)} The measures $\sigma_z$ are pairwise disjoint.

{\rm(ii)} Each of the faces $\Omega(p,q)\subset\Omega$ is a null set with
respect to $\si_z$.
\endproclaim

\demo{Derivation of Theorem 8.1.1 from Theorem 8.1.2} Let $z_1$, $z_2$ be two
distinct complex numbers from the upper half--plane $\Im z\ge0$. We have to
prove that $T_{z_1}$ and $T_{z_2}$ are disjoint. Assume first that both $z_1$
and $z_2$ are not integers. We know that if $z\in\Bbb C\setminus\Z$ then the
distinguished spherical vector $\xi_0$ is a cyclic vector; hence the measure
$\sigma_z$ determines the representation $T_z$ entirely. According to claim (i)
of Theorem 8.1.2, the measures $\si_{z-1}$ and $\si_{z_2}$ are disjoint;
therefore the representations $T_{z_1}$ and $T_{z_2}$ are disjoint, too. Assume
now that $z_1\in\C\setminus\Z$ while $z_2\in\Bbb Z$. According to Theorem
7.1.2, the representation $T_{z_2}$ decomposes into a direct sum of
representations labelled by the faces $\Omega(p,q)$, $p-q=z_2$. By virtue of
claim (ii) of Theorem 8.1.2, $T_{z_1}$ and $T_{z_2}$ are disjoint. Finally,
when $z_1$, $z_2$ are two distinct integers, the disjointness of the
representations was pointed out in Corollary 7.1.3. \qed\enddemo

Now we proceed to the proof of Theorem 8.1.2.

\subhead 8.2. Reduction to central measures \endsubhead There is a one--to--one
correspondence $\si\leftrightarrow\wt M$ between probability  measures $\si$ on
the Thoma simplex $\Om$ and central probability measures $\wt M$ on the path
space $\Cal T$, see \S9...

\proclaim{Lemma 8.2.1} Let $\sigma_1$ and $\sigma_2$ be two probability
measures on $\Omega$, and let $\wt M_1$ and $\wt M_2$ be  the corresponding
central measures on $\Cal T$. Then  $\sigma_1$, $\sigma_2$ are disjoint if and
only if $\wt M_1$ and $\wt M_2$ are disjoint.
\endproclaim

\demo{Proof} First, introduce a notation. Given two finite (not necessarily
normalized) Borel measures $\nu_1$, $\nu_2$ on a Borel space, let us denote by
$\nu_1\wedge\nu_2$ their greatest lower bound. Its existence can be verified as
follows.  Let $f_1$ and $f_2$ be the Radon--Nikodym derivatives of $\nu_1$ and
$\nu_2$ with respect to $\nu_1+\nu_2$, then we set
$\nu_1\wedge\nu_2=\min(f_1,f_2)(\nu_1+\nu_2)$. Observe that $\nu_1$ and $\nu_2$
are disjoint if and only if $\nu_1\wedge\nu_2=0$.

Next, observe that the correspondence $\si\leftrightarrow\wt M$ can be extended
to not necessarily normalized measures.

Now we can proceed to the proof. In one direction the implication is trivial.
Namely, if $\sigma_1$ and $\sigma_2$ are not disjoint, then $\si_1\wedge\si_2$
is a nonzero measure. Let $\wt M$ be the corresponding central measure. From
the integral representation of Theorem 9.7.2 it follows that $\wt M\le \wt M_1$
and $\wt M\le \wt M_2$, so that $\wt M_1\wedge\wt M_2\ne0$, whence $\wt M_1$
and $\wt M_2$ are not disjoint.

In the opposite direction, assume that $\wt M_1$ and $\wt M_2$ are not
disjoint, so that $\wt M_1\wedge \wt M_2$ is nonzero. We claim that $\wt
M_1\wedge \wt M_2$ is a central measure. Indeed, this follows from the
characterization of central measures as invariant measures with respect to a
countable group action, as explained in Proposition 9.4.1. Now, let $\si$ be
the measure on $\Om$ corresponding to $\wt M_1\wedge \wt M_2$. It is a nonzero
measure. Next, since $\wt M_1\wedge \wt M_2\le \wt M_1$ and $\wt M_1\wedge \wt
M_2\le \wt M_2$, we also have $\si\le\si_1$, $\si\le\si_2$. (Indeed, this claim
can be restated as follows: if $\wt M$, $\wt M'$ are two central probability
measures such that $\wt M\le \const\wt M'$ then the same inequality holds for
the corresponding spectral measures on $\Om$, and the latter claim can be
checked by applying Theorem 9.7.3 or Theorem 9.7.4).) Therefore, $\sigma_1$ and
$\sigma_2$ are not disjoint. \qed
\enddemo

\subhead 8.3. Proof of claim (ii) of Theorem 8.1.2 \endsubhead Let $p,q\in\Z_+$
be not equal to 0 simultaneously. Denote by $\Gamma(p,q)$ the set
$$
\Gamma(p,q) = \{(i,j)\mid 1\le i\le p,\quad j=1,2,\dots\} \cup \{(i,j)\mid 1\le
j\le q,\quad i=1,2,\dots\},
$$
(a ``fat hook'') and by $\Cal T(p,q)$ the set of paths $\tau=(\tau_n)\in\Cal T$
such that $\tau_n\subset\Gamma(p,q)$ for all $n$.

\proclaim{Lemma 8.3.1} If a measure $\sigma$ is supported by
$\Omega(p,q)\subset\Omega$, then the corresponding central measure $\wt M$ is
supported by $\Cal T(p,q)\subset\Cal T$.
\endproclaim

\demo{Proof} Let $(M^{(n)})$ be the coherent system corresponding to $\si$. We
have (see \S9...)
$$
M^{(n)}(\lambda) = \dim\lambda\, \int_{\om\in\Omega(p,q)} \wt
s_\lambda(\omega)\, \sigma(d\omega),
$$
where $\wt s_\lambda(\omega)$ is the supersymmetric Schur
polynomial
$s_\lambda(\alpha_1,\dots,\alpha_p;\beta_1,\dots,\beta_q)$. It
is known (Macdonald,  \cite{Mac, Example I.5.23 (a)}) that this
polynomial vanishes unless $\lambda$ belongs to $\Gamma(p,q)$.
Thus, for each $n$, $M^{(n)}$ vanishes outside
$\Ga(p,q)\cap\Y_n$. This is equivalent to the statement of the
lemma. \qed
\enddemo

\proclaim{Lemma 8.3.2} Let $z\in\C\setminus \Z$ and  $p,q\in \Z_+$ be fixed. We
assume that $p+q>0$ thus excluding the case $p=q=0$. Let $p_z(\lambda,\nu)$ be
the coherent system $M_z$.

There exists $\varepsilon>0$ depending on $z$, $p$, $q$ only, with the
following property. If $\lambda\subset\Gamma(p,q)$ is an arbitrary Young
diagram such that the set $\nu=\lambda\cup\{(p+1,q+1)\}$ is also a diagram,
then
$$
p_z(\lambda,\nu) \ge \varepsilon/n, \qquad n=|\lambda|.
$$
\endproclaim

\demo{Proof} Since the content of the box $\nu/\la$ is $q-p$, we have
$$
p_z(\lambda,\nu) = \frac{|z+q-p|^2}{|z|^2+n}\cdot
\frac{\dim\nu}{(n+1)\dim\lambda}.
$$
The first factor can be easily estimated from below: since $z\notin\Z$, there
exists $\varepsilon_1>0$ depending on $z$ only, such that
$$
\frac{|z+q-p|^2}{|z|^2+n} \ge {\varepsilon_1 \over n}.
$$
Now consider the second factor. It follows from the hook formula,
that
$$
\frac{\dim\nu}{(n+1)\dim\lambda} = \prod_b {h(b) \over h(b)+1},
$$
where the product is taken over the boxes $b\in\lambda$ such that either the
arm or the leg of $b$ contains the box $(p+1,q+1)$. There are exactly $p+q$
such boxes $b$, namely
$$
(p+1,j), \quad 1\le j\le q; \qquad (i,q+1), \quad 1\le i\le p.
$$
Therefore, there is a product of $p+q$ factors of the form
$k/(k+1)$, where $k=1,2,\ldots$. Each of the factors is greater
or equal to $1/2$, and the entire product is not less than
$2^{-(p+q)}$. This provides the required estimate.
\qed\enddemo

\proclaim{Lemma 8.3.3} Let $z\in\C\setminus\Z$ and let $\wt M_z$ be the central
measure on $\Cal T$ corresponding to the measure $\sigma_z$. We have $\wt
M_z(\Cal T(p,q))=0$ for all $p$, $q$.
\endproclaim

\demo{Proof} Denote by $\Cal T'(p,q)$ the set of those paths $\tau\in\Cal
T(p,q)$ that are not contained the smaller sets $\Cal T(p-1,q)$ and $\Cal
T(p,q-1)$. It suffices to prove that $\Cal T'(p,q)$ has measure 0 with respect
to $\wt M_z$.

Let $\mu$ be an arbitrary diagram in $\Gamma(p,q)$ that contains those two
boxes, set $m=|\mu|$, and denote by $\Cal T(p,q;\mu)$ the set of paths
$\tau\in\Cal T(p,q)$ with $\tau_m=\mu$. For any path $\tau=(\tau_n)\in\Cal
T'(p,q)$ there exists a number $n$ such that the diagram $\tau_n$ contains the
boxes $(p+1,q)$ and $(p,q+1)$. Consequently the set $\Cal T'(p,q)$ is the
countable sum of the sets of the form $\Cal T(p,q;\mu)$. Thus, it remains to
prove that each set $\Cal T(p,q;\mu)$ has measure 0.

It will be convenient to look at the measure $\wt M_z$ as describing a Markov
process with the transition function $p_z(\lambda,\nu)$. Set
$$
p_n=\Prob\{\tau_{n+1}\mid \mu\subseteq\tau_n\subset\Gamma(p,q)\}.
$$
The measure of the set $\Cal T(p,q;\mu)$ coincides with the probability of the
event $\tau_m=\mu$, multiplied by the product of the conditional probabilities
$\prod_{n\ge m} p_n$.

By Lemma 8.3.2, we have
$$
p_n \le 1 - {\varepsilon \over n}
$$
so that $\prod_{n\ge m} p_n=0$. \qed
\enddemo

\demo{Proof of claim (ii) of Theorem 8.1.2} Let $\si$ be the restriction of the
measure $\si_z$ to $\Omega(p,q)$, and let us show that $\si=0$. Let $\wt M$ be
the central measure corresponding to $\sigma$. According to Lemma 8.2.1, $\wt
M$ is supported by $\Cal T(p,q)$. On the other hand, it follows from Lemma
8.3.3 that $\Cal T(p,q)$ is a null set for $\wt M_z$. Since $\wt M\le \wt M_z$,
we conclude that $\wt M=0$, hence $\si=0$. \qed
\enddemo

\subhead 8.4. Proof of claim (i) of Theorem 8.1.2 \endsubhead Recall that if
$z\in\C\setminus\Z$, then the measure $M_z^{(n)}$ has nonzero weights
$M_z^{(n)}(\lambda)$ for all $\lambda\in\Bbb Y_n$, and we have an explicit
formula for $M_z^{(n)}(\lambda)$, see Theorem 4.1.1.

Fix two distinct numbers $z_1$, $z_2$ in the upper half--plane $\Im z\ge0$,
which are not integers. We have to prove that the spectral measures $\si_{z_1}$
and $\si_{z_2}$ are disjoint. By virtue of Lemma 8.2.1, it suffices to prove
that the corresponding central measures $\wt M_{z_1}$ and $\wt M_{z_2}$ are
disjoint. To simplify the notation, we set $\wt M_1=\wt M_{z_1}$, $\wt M_2=\wt
M_{z_2}$. We also denote by $(M_1^{(n)})$ and $(M_2^{(n)})$ the corresponding
coherent systems.

Introduce a sequence $\varphi_n(\tau)$ of functions on $\Cal T$,
$$
\varphi_n(\tau) = \frac {M_2^{(n)}(\tau_n)} {M_1^{(n)}(\tau_n)}, \qquad
n=1,2,\dots, \quad \tau=(\tau_n)\in\Cal T.
$$
Let $X$ be the set of paths $\tau\in\Cal T$ such that the sequence
$(\varphi_n(\tau))_{n\ge1}$ converges, as $n\to\infty$, to a finite nonzero
limit. This is a Borel subset of $\Cal T$.

\proclaim{Lemma 8.4.1} We have $\wt M_1(X) = \wt M_2(X) = 0$.
\endproclaim

\demo{Proof} We shall show that $X$ is contained in the union of the sets $\Cal
T(p,q)$, so that the claim will follow from Lemma 8.3.3.

Denote by $c_k(\tau)$ the content of the $k$th box $\tau_k\setminus\tau_{k-1}$.
By virtue of the explicit formula for $M_z^{(n)}$ (Theorem 4.1.1), we get
$$
\varphi_n(\tau) = \prod_{k=1}^n
\left|\frac
{z_2+c_k(\tau)}{z_1+c_k(\tau)}\right|^2
\frac{|z_1|^2+k-1}{|z_2|^2+k-1}.
$$
Therefore, $X$ consists of those paths $\tau$ for which the
infinite product
$$
\varphi_\infty(\tau) = \prod_{k=1}^\infty
\left|\frac
{z_2+c_k(\tau)}{z_1+c_k(\tau)}\right|^2
\frac{|z_1|^2+k-1}{|z_2|^2+k-1}
$$
converges. In particular, the $k$th factor in the product should go to $1$.
Since the second fraction in right--hand side converges to $1$, as
$k\to\infty$, we conclude that
$$
\lim_{k\to\infty}\left|\frac{z_2+c_k(\tau)}{z_1+c_k(\tau)}\right|^2 = 1, \qquad
\tau \in X.
$$

It follows from our assumptions on $z_1$, $z_2$ that the equality
$|z_2+c|^2=|z_1+c|^2$ may hold for at most one real number $c$. Indeed, this
equation on $c$ describes the set of points $c$ that are equidistant from
$-z_1$ and $-z_2$. Since $z_1\ne z_2$, this set is a line in the complex plane
$\C$, which cannot coincide with the real axis $\R$, because $z_1$, $z_2$ are
both in the upper half--plane. Thus, the line is either parallel to $\R$ (then
there is no real $c$ at all) or intersects $\R$ at a single point.

Now, we fix an arbitrary integer $c$ such that
$$
\left|\frac{z_2+c}{z_1+c}\right|^2 \ne 1.
$$
For any $\tau\in X$, the existence of the limit above implies that there is
only a finite number of integers $k$ such that $c_k(\tau)=c$. This means that
any path $\tau\in X$ may contain only a finite number of boxes $(p,q)$ on the
diagonal $q-p=c$. Therefore, $\tau$ is contained in some subset of type $\Cal
T(p,q)$, which concludes the proof. \qed
\enddemo

\proclaim{Lemma 8.4.2} Let $\wt A$ and $\wt B$ be two central probability
measures on $\Cal T$ and $(A^{(n)})$, $(B^{(n)})$ be the corresponding coherent
systems. Assume $\wt A\le \const\wt B$ and let $f(\tau)$ denote the
Radon--Nikodym derivative of $\wt A$ with respect to $\wt B$. Assume further
that $B^{(n)}(\lambda)\ne0$ for all $n$ and all $\lambda\in\Bbb Y_n$. Then
$$
\lim_{n\to\infty} \frac {A^{(n)}(\tau_n)}{B^{(n)}(\tau_n)} = f(\tau)
$$
for almost all paths $\tau=(\tau_n)\in \Cal T$ with respect to $\wt B$.
\endproclaim

\demo{Proof} Let $\Cal T^{[n]}$ denote the set of finite paths in $\Y$ going
from $\varnothing$ to a vertex in $\Y_n$. There is a natural projection $\Cal
T\to\Cal T^{[n]}$ assigning to a path $\tau$ its finite part
$\tau^{[n]}=(\tau_0,\dots,\tau_n)$. Notice that the infinite path space $\Cal
T$ can be identified with the projective limit space $\varprojlim\Cal T^{[n]}$.

Denote by $\Sigma^{[n]}$ the finite algebra of cylinder subsets with the bases
in $\Cal T^{[n]}$. The algebras $\Sigma^{[n]}$ form an increasing family, and
the union $\Sigma=\bigcup\Sigma^{[n]}$ coincides with the algebra of all Borel
sets with respect to the topology of $\Cal T$.

Consider the probability space $(\Cal T,\Sigma, \wt B)$. The
function $f$ is bounded and $\wt B$--measurable. Hence, by the
martingale theorem (cf., e.g., Shiryaev \cite{Shir, Ch.~VII,
Section 4, Theorem 3}
$$
\lim_{n\to\infty} \Bbb E(f\mid \Sigma^{[n]})=f.
$$
almost everywhere.

On the other hand, let $A^{[n]}$ and $B^{[n]}$ be the push--forwards of the
measures $\wt A$ and $\wt B$ taken with respect to the projection $\Cal
T\to\Cal T^{[n]}$. The conditional expectation $\Bbb E(f\mid \Sigma^{[n]})$ is
nothing but the function
$$
f_n(\tau) = \frac{A^{[n]}(\tau^{[n]})}{B^{[n]}(\tau^{[n]})}
$$
Since $\wt A$ is a central measure, we have
$$
A^{[n]}(\tau^{[n]})=\frac1{\dim\tau_n} \, A^{(n)}(\tau_n),
$$
and similarly
$$
B^{[n]}(\tau^{[n]})=\frac1{\dim\tau_n} \, B^{(n)}(\tau_n),
$$
It follows
$$
f_n(\tau) = \frac{A^{(n)}(\tau_n)}{ B^{(n)}(\tau_n)},
$$
and the proof is completed. \qed
\enddemo

\demo{Proof of claim (i) of Theorem 8.1.2} We have to show that the measures
$\wt M_1$ and $\wt M_2$ are disjoint. Set $\wt A=\wt M_1$, $\wt B=(\wt M_1+\wt
M_2)/2$. Then $\wt A\le 2\wt B$ and hence there exists the Radon--Nikodym
derivative of $\wt A$ with respect to $\wt B$. Denote it by $f(\tau)$. We have
$0\le f(\tau)\le 2$. The measures $\wt M_1$ and $\wt M_2$ are disjoint if and
only if  $f(\tau)$ takes only two values $0$ and $1$, almost surely with
respect to the measure $\wt B$.

On the other hand, by virtue of  Lemma 8.4.2 above, $f(\tau)$ is $\wt
B$--almost surely the limit of the functions $f_n(\tau)$. Let $Y$ be the set of
those paths $\tau$ for which the limit of $f_n(\tau)$ exists and is distinct
from 0 and 2. Observe that
$$
\gather f_n(\tau) = \frac{A^{(n)}(\tau_n)}{B^{(n)}(\tau_n)} =
2\frac{M_1^{(n)}(\tau_n)}{M_1^{(n)}(\tau_n)+M_2^{(n)}(\tau_n)} \\
= 2\left(1 + \frac{M_2^{(n)}(\tau_n)}{M_1^{(n)}(\tau_n)}\right)^{-1} = \frac2{1
+ \varphi_n(\tau)},
\endgather
$$
in the notation introduced before Lemma 8.4.1. Consequently, $Y$ coincides with
the set of those paths $\tau$ for which  $\varphi_n(\tau)$ has a finite nonzero
limit, that is,  $Y=X$. But $\wt M_1(X)=\wt M_2(X)=0$ by virtue of Lemma 8.4.1.
Hence, $\wt B(Y)=B(X)=0$, so that $f(\tau)$ is $0$ or $2$ almost surely with
respect to $\wt B$. \qed
\enddemo

The proof of Theorem 8.1.2 is completed.

\head 9. Appendix \endhead

\subhead 9.1. Young diagrams and representations of finite symmetric groups
\endsubhead
We identify partitions with Young diagrams, and denote by $\Y_n$ the set of
Young diagrams with $n$ boxes. Given a Young diagram $\la$, we denote by
$|\la|$ the number of its boxes, and by $\la'$  the transposed diagram.

Recall that $\Y_n$ is a natural set of labels for irreducible representations
of the finite symmetric group $S(n)$. Given $\la\in\Y_n$, we denote by
$\pi^\la$ the corresponding irreducible representation of $S(n)$ and by
$\chi^\la$ its character. Let
$$
\dim\la=\dim\pi^\la=\chi^\la(e).
$$

The quantity $\dim\la$ is called the {\it dimension\/} of $\la$. There are
several explicit formulas for the dimension. For example, the {\it hook
formula\/} says
$$
\dim \lambda  = \frac {n!} {\prod_{b \in \lambda} h(b)}, \qquad \lambda \vdash
n.
$$
Here the symbol $b\in\lambda$ means that $b$ is a box of $\la$, and if $(i,j)$
are its coordinates than $h(b)$, the {\it hook length} of $b$, is defined by
$$
h(b) = \lambda_i + \lambda_j' -i - j + 1.
$$

Given two Young diagrams $\la$ and $\mu$, we write $\mu\nearrow\la$ if
$\mu\subset\la$ and $|\la|=|\mu|+1$, that is, $\la$ is obtained from $\mu$ by
adding a box.

The classical {\it Young branching rule\/} says that for any $\la\in\Y_n$
$$
\left.\pi^\la\right|_{S(n-1)}\sim\sum_{\mu\in\Y_{n-1}\,:\, \mu\nearrow\la}
\pi^\mu\,.
$$
This implies
$$
\left.\chi^\la\right|_{S(n-1)}=\sum_{\mu\in\Y_{n-1}\,:\, \mu\nearrow\la}
\chi^\mu\,.
$$

\subhead 9.2. The Young graph and cotransition probabilities\endsubhead Let
$\Y$ be the set of all Young diagrams: the disjoint union of the sets $\Y_n$,
where $n=0,1,\dots$ (we agree that $\Y_0$ consists of a single element, the
empty diagram $\varnothing$). We view $\Y$ as the set of vertices of a graph,
called the {\it Young graph\/} and denoted also by $\Y$. By definition, the
edges of the Young graph are arbitrary couples $\mu\nearrow\la$. The Young
graph is a convenient way to encode the Young branching rule.

For $\mu\in\Y_{n-1}$ and $\la\in\Y_n$ set
$$
q(\mu,\la)=\cases \dim\mu/\dim\la, & \text{if $\mu\nearrow\la$} \\ 0, &
\text{otherwise.} \endcases
$$
By convention, $q(\varnothing, \la)=1$ for the single element $\la\in\Y_1$ (the
one--box diagram).

By the Young branching rule,
$$
\sum_{\mu\in\Y_{n-1}:\, \mu\nearrow\la} \dim\mu=\dim\la,
$$
so that
$$
\sum_{\mu\in\Y_{n-1},\la} q(\mu,\la)=1.
$$
The numbers $q(\mu,\la)$ are called the {\it cotransition
probabilities\/} (see Kerov \cite{Ker2}, \cite{Ker4}). They
constitute a probability distribution for any fixed $\la$
--- the {\it cotransition distribution.\/}

\subhead 9.3. Coherent systems of distributions on the Young
graph
\endsubhead Let $\De_n$ be the set of probability distributions on the finite
set $\Y_n$. This is a finite--dimensional simplex whose vertices are Dirac
measures $\de_\la$ with $\la\in\Y_n$.

For any $n=1,2,\dots$, define an affine map $\De_n\to\De_{n-1}$ by
$$
\de_\la\to\sum_{\mu\in\Y_{n-1}}q(\mu,\la)\,\de_\mu.
$$
Let
$$
\De=\varprojlim\De_n
$$
be the projective limit of the simplices taken with respect to these maps.

By the very definition, an element of $\De$ is a sequence
$M=\{M^{(n)}\}_{n=0,1,\dots}$ such that $M^{(n)}$ is a probability distribution
on $\Y_n$ and any two measures $M^{(n-1)}$, $M^{(n)}$ with consecutive indices
fulfill the {\it coherency relation}
$$
M^{(n-1)}(\mu)=\sum_{\la\in\Y_n} q(\mu,\la)M^{(n)}(\la), \qquad
\forall\mu\in\Y_{n-1}.
$$
Elements $M\in\De$ will be called {\it coherent systems of distributions.}

Alternatively, we may regard any $M\in\De$ as a real function on
$\Y$ such that

$\bullet$ First, $M(\la)\ge0$ for any $\la\in\Y$.

$\bullet$ Second, for any $n=1,2,\dots$ and any $\mu\in\Y_{n-1}$
$$
M(\mu)=\sum_{\la\in\Y_n} q(\mu,\la)M(\la).
$$

$\bullet$ Third, $M(\varnothing)=1$.

Indeed, the only point to be checked is that if $M$  satisfies the above
conditions then its restriction $M^{(n)}=M\mid_{Y_n}$ is a probability
distribution on $\Y_n$, that is, $\sum_{\la}M^{(n)}(\la)=1$. But this is
readily proved by induction on $n$.

A function $f(\la)$ on $\Y$ is called {\it harmonic\/} if it
satisfies the relation
$$
f(\la)=\sum_{\nu\searrow\la}f(\nu) \qquad \forall\la\in\Y.
$$
The coherency relation for a function $M(\la)$ is equivalent to
the harmonicity relation for the function
$f(\la)=M(\la)/\dim\la$.

For further details, see Vershik--Kerov \cite{VK2}, Kerov
\cite{Ker2}, \cite{Ker4}.

\subhead 9.4. Central measures on paths and transition
probabilities\endsubhead By a {\it path\/} in the Young graph we
mean a sequence of vertices
$$
\tau=(\tau_n\nearrow\tau_{n+1}\nearrow\dots), \qquad \tau_i\in\Y_i,
$$
which may be finite or infinite. Let $\Cal T$ be the set of all infinite paths
starting at $\varnothing$. This is a subset of the infinite product set
$\prod_{n=0}^\infty \Y_n$. We endow $\Cal T$ with the induced topology. Then it
becomes a compact totally disconnected topological space. Given a probability
measure on the path space $\Cal T$, we may speak about {\it random\/} paths.

To any coherent system $M$ on $\Bbb Y$ we assign a probability measure $\wt M$
on $\Cal T$ with cylinder probabilities defined as follows. Let $\la\in\Y_n$ be
an arbitrary vertex and $\tau$ be an arbitrary finite path going from
$\varnothing$ to $\la$, then the probability that the $\wt M$--random path goes
along $\tau$ (up to $\la$) equals $M(\la)/\dim\la$.

The measure $\wt M$ is {\it central} in the sense that the cylinder
probabilities  depend only on the final vertices $\la$ but not on the paths
$\tau_n$ chosen. Conversely, any central probability measure on the path space
comes from a (unique) coherent system.

There is a useful characterizations of central measures as invariant measures
with respect to a countable group of transformations of $\Cal T$. This group is
defined as follows. First, for each $n$ we let $\Cal G(n)$ be the group of the
transformations $g: \Cal T\to\Cal T$ such that for any path
$\tau=(\tau_n)\in\Cal T$, we have $\tau_m=(g(\tau))_m$ for all $m\ge n$.
Clearly, this is a finite group and we have $\Cal G(n)\subset\Cal G(n+1)$.
Next, we define the group $\Cal G$ as the union of the groups $\Cal G(n)$.

\proclaim{Proposition 9.4.1} A measure on $\Cal T$ is central if and only if it
is invariant under the action of $\Cal G$.
\endproclaim

Define the {\it support\/} of a coherent system $M$ as the subset
$$
\supp(M) = \{\lambda\in\Bbb Y:\; M(\lambda) \ne 0\}\subset\Y.
$$
The measure $\wt M$ is concentrated on the subspace of paths entirely contained
in $\supp(M)$. We may view $\wt M$ as a Markov chain on the state set
$\supp(M)$, with the {\it transition probabilities}
$$
p(\la,\nu) = \Prob\{\tau_{n+1} = \nu \mid \tau_n = \la\}, \qquad \la\in\Y_n,
\quad \nu\in\Y_{n+1},
$$
where $\tau=(\tau_n)$ is the random path. The transition probabilities
$p(\lambda,\nu)$ are unambiguously defined for all $\la\in\supp(M)$ by
$$
p(\lambda,\nu) = \frac {M(\nu)} {M(\lambda)} \, \cdot \, \frac {\dim \lambda}
{\dim \nu}, \quad \la \in \supp(M).
$$

The system of transition probabilities uniquely determines the initial central
measure, so that distinct central measures have distinct transition
probabilities. On the other hand, all central measures have one and the same
system of cotransition probabilities, which are nothing but the quantities
$q(\mu,\la)$ introduced in \S9.2. That is, we have
$$
q(\mu,\la) = \Prob\{\tau_{n-1} =\mu\mid \tau_n=\la\}, \qquad \mu\in\Y_{n-1},
\quad \la\in\Y_n.
$$

For further details, see Vershik--Kerov \cite{VK2}, Kerov
\cite{Ker2}, \cite{Ker4}.

\subhead 9.5. Characters of the group $S(\infty)$ \endsubhead Recall that we
have defined the {\it infinite symmetric group\/} $S(\infty)$ as the inductive
limit of the finite symmetric group $S(n)$ as $n\to\infty$.  A function $\chi:
S(\infty)\to\C$ will be called a {\it character\/} if it is central (i.e.
constant on conjugacy classes), positive definite, and normalized at the unity
(i.e. $\chi(e)=1$). The set of all characters of $S(\infty)$ will be denoted by
$\Cal X$.

\proclaim{Proposition 9.5.1} There is a natural bijective correspondence
$$
\Cal X\ni\chi\leftrightarrow M\in\De
$$
between characters of $S(\infty)$ and coherent systems of
probability measures on the Young graph.
\endproclaim

\demo{Proof} Let $\chi\in\Cal X$. For any $n$, set
$\chi_n=\left.\chi\right|_{S(n)}$. This is a central, positive definite,
normalized  function on $S(n)$. As readily verified, such functions are exactly
the convex combinations of {\it normalized\/} irreducible characters
$\chi^\la/\dim\la$. Thus,
$$
\chi_n=\sum_{\la\in\Y_n}M^{(n)}(\la)\,\frac{\chi^\la}{\dim\la}
$$
with certain nonnegative coefficients $M^{(n)}(\la)$,
$$
\sum_{\la\in\Y_n}M^{(n)}(\la)=1.
$$
These coefficients  may be viewed as the {\it Fourier coefficients\/} of the
function $\chi_n$. They form a probability distribution on $\Y_n$; denote it by
$M^{(n)}$. Let us check that the distributions $M^{(n)}$ obey the coherency
relation. Indeed, by virtue of the Young branching rule, this is simply
equivalent to saying that the function $\chi_{n-1}$ coincides with the
restriction of $\chi_n$ to $S(n-1)$. Thus, starting from a character $\chi$ we
obtain a coherent system $M=\{M^{(n)}\}$.

Conversely, let $M$ be a coherent system. Then, for any $n$, we may define a
function $\chi_n$ on $S(n)$ as above. These functions are pairwise compatible
and hence define a function $\chi$ on the group $S(\infty)$. It is readily
verified that $\chi$ is a character. \qed
\enddemo

Both $\Cal X$ and $\De$ are certain sets of functions (on $S(\infty)$ and $\Y$,
respectively). These sets are convex, and the bijection $\Cal
X\leftrightarrow\De$ is an isomorphism of convex sets.

Further, both $\Cal X$ and $\De$ are compact topological spaces with respect to
the topology of pointwise convergence, and our bijection is a homeomorphism
with respect to this topology.

\subhead 9.6. Thoma's theorem \endsubhead Let $\Ex\Cal X$ denote the set of
extreme points of the convex set $\Cal X$. Elements of $\Ex\Cal X$ will be
called {\it extremal characters\/} of the group $S(\infty)$.

The first examples of extremal characters are as follows. Let
$\Om(p,q)\subset\R^p\times\R^q$ be the set of couples $(\al,\be)$, where
$\al=(\al_1\ge\dots\ge\al_p\ge0)$ and $\be=(\be_1\ge\dots\ge\be_q\ge0)$ be two
collections of numbers such that
$$
\sum_{i=1}^p \al_i+\sum_{j=1}^q\be_j=1.
$$
Here one of the numbers $p,q$ may be zero (then the corresponding collection
$\al$ or $\be$ disappears).

For $(\al,\be)\in\Om(p,q)$ and $k=1,2,\dots$, set
$$
\wt p_k(\al,\be)=\sum_{i=1}^p \al_i^k+(-1)^{k-1}\sum_{j=1}^q\be_j^k.
$$
Notice that
$$
\wt p_1(\al,\be)\equiv1.
$$
Given $s\in S(\infty)$, we denote by $m_k(s)$ the number of $k$--cycles in $s$.
Since $s$ is a finite permutation, we have
$$
m_1(s)=\infty, \qquad \text{$m_k(s)<\infty$ for $k\ge2$}, \qquad
\text{$m_k(s)=0$ for $k$ large enough}.
$$

In this notation, we define a function on $S(\infty)$ by
$$
\chi^{(\al,\be)}(s)=\prod_{k=1}^\infty (\wt
p_k(\al,\be))^{m_k(s)}=\prod_{k=2}^\infty (\wt p_k(\al,\be))^{m_k(s)}, \qquad
s\in S(\infty),
$$
where we agree that $1^\infty=1$ and $0^0=1$. Any such function turns out to be
an extremal character: this claim is a particular case of a more general result
stated below.

If $p=1$ and $q=0$ (i.e., $\al_1=1$ and all other parameters disappear) then we
get the trivial character, which equals 1 identically. If $p=0$ and $q=1$ then
we get the alternate character $\sgn(s)=\pm1$, where the plus--minus sign is
chosen according to the parity of the permutation. More generally, we have
$$
\chi^{(\al,\be)}\cdot\sgn=\chi^{(\be,\al)}.
$$

Let $\R^\infty$ denote the direct product of countably many copies of $\R$. We
equip $\R^\infty$ with the product topology. Let $\Om$ be the subset of
$\R^\infty\times\R^\infty$ formed by couples $\al\in\R^\infty$,
$\be\in\R^\infty$ such that
$$
\al=(\al_1\ge\al_2\ge\dots\ge0), \qquad \be=(\be_1\ge\be_2\ge\dots\ge0), \qquad
\sum_{i=1}^\infty \al_i+\sum_{j=1}^\infty\be_j\le1.
$$

We call $\Om$ the {\it Thoma simplex.\/} As affine coordinates of the simplex
one can take the numbers
$$
\al_1-\al_2,\,\dots,\, \al_{p-1}-\al_p, \,\al_p, \, \be_1-\be_2,\,\dots,
\,\be_{q-1}-\be_q, \, \be_q
$$
but we will not use these coordinates. We equip $\Om$ with topology induced
from that of the space $\R^\infty\times\R^\infty$. It is readily seen that
$\Om$ is a compact space. Clearly, each set $\Om(p,q)$ may be viewed as a
subset of $\Om$ (this is one of finite--dimensional faces of $\Om$).

Notice that the union of the simplices $\Om(p,q)$ is dense in $\Om$. For
instance, the point $(\underline0,\underline0)=(\al\equiv0,\be\equiv0)\in\Om$
can be approximated by points of the simplices $\Om_{p0}$ as $p\to\infty$,
$$
(\underline0,\underline0)=\lim_{p\to\infty} ((\underbrace{1/p,\dots,1/p}_p),
\underline0).
$$

Now we extend by continuity the definition of the functions $\chi^{(\al,\be)}$
given above. First, for any $k=2,3,\dots$ we define the function $\wt p_k$ on
$\Om$ as follows. If $\om=(\al,\be)\in\Om$ then
$$
\wt p_k(\om)=\wt p_k(\al,\be) =\sum_{i=1}^\infty
\al_i^k+(-1)^{k-1}\sum_{j=1}^\infty\be_j^k.
$$
Note that $\wt p_k$ is a continuous function on $\Om$. It should be emphasized
that the condition $k\ge2$ is necessary here: the similar expression with $k=1$
(that is, the sum of all coordinates) {\it is not\/} continuous.

Next, for any $\om=(\al,\be)\in\Om$ we set
$$
\chi^{(\om)}(s)=\chi^{(\al,\be)}(s)=\prod_{k=2}^\infty (\wt
p_k(\al,\be))^{m_k(s)}, \qquad s\in S(\infty),
$$

\proclaim{Theorem 9.6.1 (Thoma's theorem)} The functions
$\chi^{(\om)}$ are precisely the extremal characters of the
group $S(\infty)$.
\endproclaim

That is, for any $\om\in\Om$ the function $\chi^{(\om)}$  is an extremal
character, each extremal character is obtained in this way, and different
points $\om\in\Om$ define different characters.

In particular, the character $\chi^{(\underline0,\underline0)}$ is the delta
function at $e\in S(\infty)$. It corresponds to the biregular representation
defined in \S1.3.

Notice that the set $\Cal X$ carries a natural topology -- that
of pointwise convergence on the group. Endow the subset $\Ex\Cal
X\subset \Cal X$ with the induced topology. Then the
correspondence $\Ex\Cal X\leftrightarrow\Om$ given by Thoma's
theorem becomes a homeomorphism of topological spaces.

This implies, in particular, that  characters $\chi^{(\om)}$
with parameters $\om\in\cup_{p,q}\Om(p,q)$ are dense in the
whole set $\Ex\Cal X$ with respect to the topology of pointwise
convergence on the group $S(\infty)$.

\demo{Comments to Thoma's theorem} 1. The original proof of
Thoma was given in his paper \cite{Tho1} published in 1964.
Thoma first proved (Satz 1 in \cite{Tho1}) that a character
$\chi$ is extremal if and only if it is a multiplicative class
function, that is,
$$
\chi(s)=\prod_{k=2}^\infty p_k^{m_k(s)}, \qquad s\in S(\infty),
$$
with certain real numbers $p_2,p_3,\dots$. This reduced the
problem to the following one: find all sequences
$(p_2,p_3,\dots)$ such that the expression above is a positive
definite function on the group $S(\infty)$. An equivalent
condition on $(p_2, p_3,\dots)$ is as follows: let
$h_1,h_2,\dots$ be defined by
$$
1+h_1u+h_2u^2+h_3u^3+\dots=\exp\Big(u+\sum_{k=2}^\infty
p_k\frac{u^k}k \Big)
$$
and set $h_0=1$, $h_{-1}=h_{-2}=\dots=0$; then
$$
\det_{1\le i,j\le \ell(\la)}[h_{\la_i-i+j}]\ge0 \qquad \text{for
any $\la\in\Y$}.
$$

Then Thoma succeeded to prove that the sequences
$(h_1=1,h_2,h_3,\dots)$ with this property  are exactly those
given by the formula
$$
1+h_1u+h_2u^2+h_3u^3+\dots=e^{\ga u}\prod_{i=1}^\infty
\frac{1+\be_i u}{1-\al_i u}\,,
$$
where
$$
(\al,\be)\in\Om, \quad \ga=1-\sum_{i=1}^\infty(\al_i+\be_i),
$$
which implies the theorem. Actually, this result is equivalent
to Edrei's classification (1952) of one--sided, totally positive
sequences in the sense of Schoenberg, see Edrei \cite{Edr}.
Since the paper \cite{Tho1} contains no reference to Schoenberg
or Edrei, one may conclude that Thoma was unaware of their work.

2. Quite a different proof of Thoma's theorem was given by
Vershik and Kerov (\cite{VK1}, \cite{VK2}, 1981). Instead of
function--theoretic arguments of Edrei and Thoma, Vershik and
Kerov used an asymptotic method (whose general idea was
suggested by Vershik's paper \cite{Ver}): approximation of
extremal characters of $S(\infty)$ by irreducible normalized
characters $\chi^\la/\dim\la$ of finite groups $S(n)$, as
$n\to\infty$. The asymptotic method explains the origin of
Thoma's parameters $\al_i,\be_i$: they arise as limits of
normalized Frobenius coordinates of the growing diagram $\la$.

3. An important combinatorial lemma, stated in \cite{VK2}
without proof (see \cite{VK2, \S5, Lemma 1}), was proved in
Kerov--Olshanski \cite{KO}. A particular case of it (which is
sufficient for completing the proof of Thoma's theorem) was
proved in Wassermann's thesis (\cite{Was}). For more detail, see
also Okounkov--Olshanski \cite{OkOl, \S8},
Olshanski--Regev--Vershik \cite{ORV}.

4. Okounkov's work \cite{Ok1, Ok2} provides one more approach to
Thoma's theorem. In particular, Okounkov showed that a crucial
step in Thoma's proof can be replaced by a simple
representation--theoretic argument.

5. Finally, the paper Kerov--Okounkov--Olshanski \cite{KOO}
contains a far generalization of Thoma's theorem obtained by the
asymptotic method of \cite{VK2}.

\enddemo

\subhead 9.7. Spectral decomposition of characters \endsubhead

\proclaim{Theorem 9.7.1} {\rm(i)} For any character $\chi\in\Cal X$, there
exists a probability measure $\si$ on the Thoma simplex $\Om$ such that
$$
\chi(s)=\int_{\Om} \chi^{(\om)}(s)\si(d\om), \qquad s\in S(\infty).
$$

{\rm(ii)} Such a measure is unique.

{\rm(iii)} Conversely, for any probability measure $\si$ on $\Om$, the function
$\chi$ defined by the above formula is a character of $S(\infty)$.

Thus, $\Cal X$ is isomorphic, as a convex set, to the set of all probability
measures on the compact space $\Om$.
\endproclaim

We call this integral representation the {\it spectral decomposition\/} of a
character. The measure $\si$ will be called the {\it spectral measure\/} of
$\chi$. If $\chi$ is extremal, i.e., $\chi=\chi^{(\om)}$, then its spectral
measure reduces to the Dirac mass at $\om$.

Theorem 9.7.1 admits an equivalent formulation in terms of
coherent systems. To state it we need to extend the above
definition of the functions $\wt p_k(\om)$ to arbitrary
symmetric functions. Let $\La$ be the algebra of symmetric
functions (see Macdonald \cite{Mac}). The power--sums $p_k$ are
algebraically independent generators of $\La$, so that the
assignment $p_k\mapsto \wt p_k$ can be extended to a
homomorphism of the algebra $\La$ into the algebra $C(\Om)$ of
continuous functions on the Thoma simplex $\Om$. Given
$f\in\La$, we denote by $\wt f$ its image in $C(\Om)$. In
particular, we apply this to the Schur functions $s_\la$: the
corresponding functions $\wt s_\la(\om)$ are called the {\it
extended Schur functions\/} (see Vershik--Kerov \cite{VK3, \S6},
Kerov--Okounkov--Olshanski \cite{KOO, Appendix}). Notice that
the restriction of $\wt s_\la$ to $\Om(p,q)\subset\Om$ is
nothing but the {\it supersymmetric\/} Schur polynomial
$s_\la(\al_1,\dots,\al_p;\be_1,\dots,\be_q))$ (see, e.g.,
Berele--Regev \cite{BR}, Macdonald \cite{Mac, Example I.3.23};
in Macdonald's notation, this is
$s_\la(\al_1,\dots,\al_p;-\be_1,\dots,-\be_q)$).

\proclaim{Theorem 9.7.2} Let $\chi$ be a character of $S(\infty)$,
$M=(M^{(n)})$ be the corresponding coherent system, and $\si$ be the spectral
measure of $\chi$. For any $n=0,1,2,\dots$ and any $\la\in\Y_n$ we have
$$
M^{(n)}(\la)=\int_\Om \dim \la \cdot \wt s_\la(\om)\,\si(d\om).
$$
\endproclaim

In particular, the coherent system corresponding to an extremal character
$\chi^{(\om)}$ has the form $M^{(n)}(\la)=\dim \la \cdot \wt s_\la(\om)$.

For a proof of Theorem 9.7.2, see Kerov--Okounkov--Olshanski \cite{KOO} (that
paper actually contained a more general result).

The claim of the theorem  is similar to the Poisson integral
representation of harmonic functions (see Kerov \cite{Ker2},
\cite{Ker4}, Kerov--Okounkov--Olshanski \cite{KOO}). Notice that
the role of the Poisson kernel is played here by the function
$(\la,\om)\mapsto \wt s_\la(\om)$.

Theorems 9.7.1 and 9.7.2 involve the claim that a certain convex
set (that of characters or coherent systems, or yet
equivalently, that of central measures) is a {\it Choquet
simplex\/}. That is, each point of the convex set in question is
{\it uniquely\/} representable by a probability measure on the
subset of extreme points. This fact does not rely on the
specific nature of the group $S(\infty)$ or the Young graph, and
can be derived from some very general theorems. See, e.g.,
Diaconis--Freedman \cite{DS, section 4}, Olshanski \cite{Ol7,
\S9}.

The next result can be viewed as a kind of Fatou's theorem on boundary values
of harmonic functions. To state it we need to introduce some important notation
and definitions.

Recall the definition of the {\it Frobenius coordinates\/} of a nonempty
diagram $\la\in\Y$: these are the integers $p_1>\dots> p_d\ge0$, $q_1>\dots>
q_d\ge0$, where $d$ is the number of boxes on the main diagonal of $\la$ and
$$
p_i=\la_i-i, \quad q_i=\la'_i-i, \qquad i=1,\dots,d.
$$
Any collection of integers $p_1>\dots> p_d\ge0$, $q_1>\dots> q_d\ge0$
corresponds to a Young diagram. The transposition $\la\mapsto\la'$ corresponds
to interchanging $p_i\leftrightarrow q_i$.

We also need the so called {\it modified\/} Frobenius coordinates of a diagram
$\la$, which are defined by
$$
\wt p_i=p_i+\tfrac12, \quad \wt q_i=q_i+\tfrac12, \qquad i=1,\dots,d.
$$
It is convenient to agree that
$$
\wt p_i=\wt q_i=0, \qquad i>d.
$$
Notice that
$$
\sum_{i=1}^\infty (\wt p_i+\wt q_i)=|\la|.
$$

For each $n$ we define an embedding $\iota_n: \Y_n\to\Om$ by
$$
\iota_n(\la)=(\al,\be), \qquad \al_i=\frac{\wt p_i}n, \quad \be_i=\frac{\wt
q_i}n, \quad i=1,2,\dots, \quad \la\in\Y_n.
$$
Notice that the union of the finite sets $\iota_n(\Y_n)$ is dense in $\Om$.

\proclaim{Theorem 9.7.3} Let $\chi$ be a character of $S(\infty)$, $\si$ be its
spectral measure, and $M=(M^{(n)})$ be the coherent system of distributions
corresponding to $\chi$. Further, let $\iota_n(M^{(n)})$ be the push--forward
of the measure $M^{(n)}$ under the embedding $\iota_n:\Y_n\to\Om$.

As $n\to\infty$, the measures $\iota_n(M^{(n)})$ converge to $\si$ in the weak
topology of measures on the compact space $\Om$.
\endproclaim

For a proof, see Kerov--Okounkov--Olshanski \cite{KOO}.

Finally, we shall state a related result concerning central measures. An
infinite path $\tau=(\tau_n)\in\Cal T$ is called {\it regular\/} if the points
$\iota_n(\tau_n)\in\Om$ converge to a limit as $n\to\infty$; then the limit
point is called the {\it end\/} of the path. Let $\Cal T'$ be the set of
regular paths; this is a Borel subset of $\Cal T$. Assigning to a regular path
its end we get a Borel map $\Cal T'\to\Om$.

\proclaim{Theorem 9.7.4} Any central probability measure $\wt M$ is supported
by $\Cal T'$. The push-forward of $\wt M$ under the map $\Cal T'\to\Om$
coincides with the spectral measure $\si$ of the character $\chi\leftrightarrow
\wt M$.
\endproclaim

The proof is similar to (and actually simpler than) the proof of
Theorem 10.2 in Olshanski \cite{Olsh7}.

\subhead 9.8. Spherical representations and spherical functions \endsubhead Let
$G$ be the bisymmetric group $S(\infty)\times S(\infty)$ and $K$ be the
diagonal subgroup in $G$, canonically isomorphic to $S(\infty)$.

Assume we are given a unitary representation $T$ of the group $G$ in a Hilbert
space $\Cal H$. A vector $\xi\in\Cal H$ is said to be a {\it cyclic vector\/}
if the linear span of the vectors $T(g)\xi$, where $g$ ranges over $G$, is
dense in $\Cal H$. Suppose that $\xi$ is cyclic, invariant under the action of
the subgroup $K$, and $\Vert\xi\Vert=1$. In such a case we say that the couple
$(T,\xi)$ is a {\it spherical representation\/}. We will call $\xi$ the {\it
spherical vector.}

Denote by $\Phi$ the set of all functions on $G$ that are positive definite,
$K$--bi--invariant, and normalized at the unity. If $(T,\xi)$ is a spherical
representation of $(G,K)$, then  the matrix coefficient corresponding to the
spherical vector,
$$
\varphi(g)=(T(g)\xi,\xi), \qquad g\in G,
$$
is an element of $\Phi$. We call $\varphi$ the {\it spherical function\/} of
$(T,\xi)$. The  couple $(T,\xi)$ is uniquely (up to a natural equivalence)
reconstructed from its spherical function, by use of the Gelfand--Naimark
construction. Moreover, any $\varphi\in\Phi$ comes from a certain $(T,\xi)$.
Thus, there is a one--to--one correspondence between functions $\varphi\in\Phi$
and (equivalence classes of) spherical representations $(T,\xi)$.

Assume $T$ is an irreducible unitary representation of $G$. Then
the space of its $K$--invariant vectors has dimension 0 or 1
(indeed, this follows from the fact that $(G,K)$ is a Gelfand
pair, see Olshanski \cite{Ol3, \S1}). Thus, if $T$ possesses a
nonzero $K$--invariant vector $\xi$ then $\xi$ is unique, within
a scalar multiple. Observe that $\xi$ is automatically cyclic,
because any nonzero vector in an irreducible representation is
cyclic. Thus, assuming $\Vert\xi\Vert=1$, we see that  $(T,\xi)$
is a spherical representation. The only lack of uniqueness in
the choice of $\xi$ is reduced to multiplying $\xi$ by a complex
scalar of absolute value 1, which does not affect the spherical
function $\varphi(g)=(T(g)\xi,\xi)$. Notice that $\varphi$ is an
extreme point of the convex set $\Phi$. Conversely, if
$\varphi\in\Phi$ is extreme then the corresponding spherical
representation is irreducible.

\proclaim{Proposition 9.8.1} There is a natural bijective correspondence
$\chi\leftrightarrow\varphi$ between characters $\chi\in\Cal X$ and spherical
functions $\varphi\in\Phi$.
\endproclaim

\demo{Proof} Indeed, the relation between $\chi$ and $\varphi$ has the form
$$
\varphi(g_1,g_2)=\chi(g_1g_2^{-1}), \qquad \chi(s)=\varphi(s,e),
$$
where $g_1,g_2,s$ are elements of $S(\infty)$. Clearly, the normalization
$\chi(e)=1$ is equivalent to the normalization $\varphi(e)=1$. It is readily
verified that $\chi$ is constant on conjugacy classes if and only if $\varphi$
is constant on double cosets. Next, let $\{g_i\}=\{(g_{i1},g_{i2})\}$ be a
finite collection of element of the group $G$, and let $s_i=g_{i1}g_{i2}^{-1}$
be the corresponding elements in $S(\infty)$. Remark that $g_j^{-1}g_i$ lies in
the same double coset modulo $K$ as $(s_j^{-1}s_i,e)$. It follows that
$\varphi(g_j^{-1}g_i)=\chi(s_j^{-1}s_i)$, so that $\varphi$ is positive
definite if and only if $\chi$ is. Thus, $\chi\leftrightarrow\varphi$ is indeed
a bijective correspondence between the two sets. \qed\enddemo

Clearly, the bijection $\chi\leftrightarrow\varphi$ is an isomorphism of convex
set. Therefore, irreducible spherical representations of $(G,K)$ are
parametrized by extremal characters, and finally by points $\om\in\Om$.

More generally, combining Proposition 9.8.1 with the description of characters
given in Theorem 9.7.1 we obtain a general description of spherical
representations $(T,\xi)$. Specifically, any such $(T,\xi)$ is determined by a
probability measure $\si$ on $\Om$. It is worth noting that, as long we are
dealing with reducible spherical representations, a given $T$ may well possess
a lot of $K$--invariant cyclic vectors. If $\xi$ is replaced by another
spherical vector $\xi'$ then $\si$ is replaced by an equivalent probability
measure $\si'$. Thus, the equivalence class of $T$ is determined by the
equivalence class of $\si$. This equivalence class of measures on $\Om$ will be
called the {\it spectral type of $T$.}

Using the abstract machinery of direct integrals of Hilbert spaces one can show
that any spherical representation $(T,\xi)$ can be decomposed into a
multiplicity free direct integral of irreducible spherical representations
$(T^{(\om)},\xi^{\om})$. This decomposition, which may be called the {\it
spectral decomposition of $(T,\xi)$,\/} is unique, and it is governed by the
spectral measure $\si$:
$$
T=\int_\om T^{(\om)}\si(d\om), \qquad \xi=\int_\om \xi^{(\om)}\si(d\om).
$$

Assume $(T',\xi')$ is another spherical representation and $\si'$ is the
corresponding spectral measure. The measures $\si$ and $\si'$ are said to be
{\it disjoint\/} if they are singular with respect to each other (then there
exist two disjoint Borel sets supporting them). The representations $T$ and
$T'$ are said to be {\it disjoint\/} if they have no equivalent nonzero
subrepresentations. Disjointness of  $\si$ and $\si'$ is equivalent to
disjointness of  $T$ and $T'$.

\subhead 9.9. Admissible representations \endsubhead Spherical representations
enter a wider class of representations that will be defined now.

For $m \le n$, let $S_m(n)\subseteq S(n)$ denote the subgroup fixing the points
$1,2,\,\ldots,m$. Set
$$
S_m(\infty) = \bigcup_{n \ge m} S_m(n) \subset S(\infty)
$$
and denote by $K_m\subset K$ the corresponding subgroup of $K\cong S(\infty)$.
Recall that by $G(m)$ we denote the subgroup $S(m)\times S(m)$ in the
bi--symmetric group $G$. An important fact is that the subgroups $K_m$ and
$G(m)$ commute to each other.

Given a unitary representation $T$ of the group $G$ in a Hilbert space $\Cal
H$, let $\Cal H_m = \Cal H^{K_m}$ be the subspace of $K_m$--invariant vectors,
and set
$$
\Cal H_\infty = \bigcup_m \Cal H_m\,.
$$

We remark that $\Cal H_\infty$ is a $G$--invariant (algebraic) subspace in
$\Cal H$. Indeed, for any $m$, $\Cal H_m$ is invariant under $G(m)$, because
$G(m)$ and $K_m$ commute. Since $G$ is the union of $G(m)$'s, it follows that
$\Cal H_\infty$ is invariant under $G$. Thus, the closure of $\Cal H_\infty$ is
an invariant subspace of the representation $T$.

We say that $T$ is an {\it admissible} representation of the pair $(G,K)$ if
the subspace $\Cal H_\infty$ as defined above is dense in $\Cal H$.

For more detail about this definition, see \cite{Ol3} and also
\cite{Ol2}, \cite{Ol4}. Not all representations of $G$ are
admissible, for it may well happen that $\Cal H_\infty$ is
reduced to $\{0\}$. If $T$ is irreducible then either it is
admissible or $\Cal H_\infty=\{0\}$.

Any spherical representation $(T,\xi)$ is admissible. Indeed, the spherical
vector $\xi$ belongs to the subspace $\Cal H_\infty$. Therefore, all vectors
$T(g)\xi$, where $g\in G$, are also in $\Cal H_\infty$. Since these vectors
generate a dense algebraic subspace, $\Cal H_\infty$ is dense in $\Cal H$, so
that $T$ is admissible.

As shown in \cite{Ol3}, admissible representations are exactly
those unitary representations that can be continuously extended
to the topological group $\overline G\supset G$. Thus,
admissible representations of $(G,K)$ are in essence the same as
continuous unitary representations of the group $\overline G$.
However, for technical reasons, admissible representations are
more convenient to deal with than representations of $\overline
G$.

Any admissible representation is a type I representation, i.e.,
the von Neumann algebra generated by it is of type I (\cite{Ol3,
Theorem 4.1}). This means that inside the class of admissible
representation there is no pathologies occurring for general
representations of non--tame groups. Irreducible admissible
representations admit a complete classification (Okounkov
\cite{Ok1}, \cite{Ok2}), their explicit realization is described
in \cite{Ol3, \S5}.

\enddocument
\end